\documentclass[12pt,twocolumn,a4paper]{article}

\usepackage{amsmath}
\usepackage{amsthm}
\usepackage{amssymb}
\usepackage[utf8]{inputenc}
\usepackage{tikz}
\usepackage{color}
\usepackage{caption}
\usepackage{subcaption}
\usepackage{cite}

\usepackage{hyperref}
\hypersetup{pdftex,colorlinks=true,allcolors=blue}
\usepackage{hypcap}
\usepackage{pgfplots}

\usetikzlibrary{decorations.markings}
\tikzset{->-/.style={decoration={
  markings,
  mark=at position #1 with {\arrow{>}}},postaction={decorate}}}
\tikzset{-<-/.style={decoration={
  markings,
  mark=at position #1 with {\arrow{<}}},postaction={decorate}}}

\def\defn#1{{\bf #1}}
\def\multi#1{{\begin{tabular}{@{}l@{}}#1\end{tabular}}}

\def\innerprod(#1,#2){{\left<#1\,,\,#2\right>}}

\usepackage{fancyhdr}
\fancypagestyle{firststyle}
{
   \fancyhf{}
   
\fancyfoot[C]{1}
\fancyfoot[R]{\footnotesize\bf File: \jobname.tex}
}
\def\VE{{\boldsymbol E}}
\def\VB{{\boldsymbol B}}
\def\VD{{\boldsymbol D}}
\def\VH{{\boldsymbol H}}
\def\VJ{{\boldsymbol J}}
\def\Jcurr{{\cal J}}
\def\Hform{{\cal H}}

\newenvironment{jgitemize}{\begin{list}{$\bullet$}
{
\setlength{\itemindent}{1.5 em}
\setlength{\itemsep}{1pt}
\setlength{\labelsep}{0.5em}
\setlength{\labelwidth}{15em} 
\setlength{\leftmargin}{1em}
\setlength{\parsep}{0em}
\setlength{\listparindent}{1em}
\setlength{\parskip}{0em}
\setlength{\partopsep}{0pt}
\setlength{\topsep}{0pt}
}}
{\end{list}}

\newlength{\figsize}
\title{A pictorial introduction to differential geometry, \\leading to
Maxwell's equations as three pictures. 
}

\author{Dr Jonathan Gratus
\\\normalsize Physics Department, Lancaster University
  and the Cockcroft Institute of accelorator Science.
\\\normalsize \texttt{j.gratus@lancaster.ac.uk}
\\\normalsize \texttt{http://www.lancaster.ac.uk/physics/about-us/people/jonathan-gratus}}

\begin{document}


\maketitle
\thispagestyle{firststyle} 

\begin{abstract}
In this article we present pictorially the foundation of differential
geometry which is a crucial tool for multiple areas of physics, notably
general and special relativity, but also mechanics, thermodynamics and solving
differential equations. As all the concepts are presented as pictures,
there are no equations in this article. As such this article may be
read by pre-university students who enjoy physics, mathematics and
geometry. However it will also greatly aid the intuition of an
undergraduate and masters students, learning general relativity and
similar courses. It concentrates on the tools needed to understand
Maxwell's equations thus leading to the goal of presenting Maxwell's
equations as 3 pictures. 
\end{abstract}

\subsection*{Prefix}

When I was young, somewhere around 12, I was given a book on relativity,
gravitation and cosmology. Being dyslexic I found reading the text
torturous. However I really enjoyed the pictures. To me, even at that
age, understanding spacetime diagrams was natural. It was obvious, once
it was explained, that a rocket ship could not travel more in space
than in time and hence more horizontal than $45^\circ$. Thus I could
easily understand why, having entered a stationary, uncharged black
hole it was impossible to leave and that you were doomed to reach the
singularity. Likewise for charged black holes you could escape to
another universe. 

I hope that this document may give anyone enthusiastic enough to get a
feel for differential geometry with only a minimal mathematical or
physics education. However it helps having a good imagination, to
picture things in 3 dimension (and possibly 4 dimension) and a good supply of
pipe cleaners.

I teach a masters course in differential geometry to physicists and
this document should help them to get some intuition before embarking
on the heavy symbol bashing.

{\small 
\tableofcontents

\listoffigures
\listoftables
}

\section{Introduction}
\label{ch_Intro}

Differential geometry is a incredibly useful tool in both physics and
mathematics. Due to the heavy technical mathematics need to define all the
objects precisely it is not usually introduced until the final years in
undergraduate study. Physics undergraduates usually only see it in the context
of general relativity, where spacetime is introduced. This is unfortunate as
differential geometry can be applied to so many objects in physics. These
include: 
\begin{jgitemize}
\item
The 3-dimensional space of Newtonian physics. 
\item
The 4-dimensional spacetime of special and general relativity.
\item
2-dimensional shapes, like spheres and torii.
\item
6-dimensional phase space in Newtonian mechanics.
\item
7-dimensional phase-time space in  general relativity.
\item
The state space in thermodynamics.
\item
The configuration space of  physical system.
\item
The solution space in differential equations.
\end{jgitemize}

There are many texts on differential geometry, many of which use diagrams to
illustrate the concepts they are trying to portray, for example
\cite{warnick2014differential,burke1985applied,jancewicz1989multivectors,tonti2001finite,warnick1996electromagnetics,misner1957classical,hehl2000gentle,schouten1954tensor}. They
probably date back to Schouten \cite{schouten1954tensor}. The novel approach adopted
here is to present as much differential geometry as possible solely by using
pictures. As such there are almost no equations in this article.
Therefore this document may be used by first year undergraduates, or even keen
school students to gain some intuition of differential geometry. Students
formally studying differential geometry may use this text in conjunction with a
lecture course or standard text book.

In geometry there is always a tension between drawing pictures and
manipulating algebra.
Whereas the former can give you intuition and some simple results in
low dimensions, only by expressing geometry in terms of mathematical
symbols and manipulating them can one derive deeper results. In
addition, symbols do not care whether one is in 1, 3, 4 dimensions or
even 101 dimensions (the one with the black and white spots), whereas
pictures can only clearly represent objects in 1, 2, or 3
dimensions. This is unfortunate as we live in a 4-dimensional universe
(spacetime) and one of the goals of this document is to represent
electrodynamics and Maxwell's equations.

As well as presenting differential geometry without equations, this document is
novel because it concentrates on two aspects of differential geometry which
students often find difficult: Exterior differential forms and
orientations. This is done even before introducing the metric.  

Exterior differential forms, which herein we will simply call \defn{forms},
are important both in physics and mathematics. In electromagnetism
they are the natural way to define the electromagnetic fields. In
mathematics they are the principal objects that can be integrated. In this
article we represent forms pictorially using curves and surfaces.

Likewise orientations are important for both integration and
electromagnetism. There are two types of orientation, internal and
external, also called twisted and untwisted, and these are most
clearly demonstrated pictorially.  Some physicists will be familiar
with the statement that the magnetic field is a twisted vector field,
also known as a pseudo-vectors or an axial vector. Here the
twistedness is extended to forms and submanifolds. 

The metric to late is one of the most important objects in
differential geometry. In fact without it some authors say you are not in fact
studying differential geometry but instead another subject called differential
topology. Without a metric you can think of a manifold as made of a rubber,
which is infinity deformable, whereas with a metric the manifold is made of
concrete.  The metric gives you the length of vectors and the angle between two
vectors. It can also tell you the distance between two points. In higher
dimensions the metric tells you the area, volume or $n$-volume of your manifold
and submanifolds within it. It gives you a measure so that you can integrate
scalar fields. 

General relativity is principally the study of 4-dimensional manifolds
with a spacetime metric. It is the metric which defines curvature and
therefore gravity. General relativity therefore is the study of how
the metric affects particles and fields on the manifold and how these
fields and particles effect the metric. There are many popular science
books and programs which introduce general relativity by showing how
light and other object move in a curved space, usually by showing
balls moving on a curved surface. Therefore one of the things we do not
cover in this document are the spacetime diagrams in special
relativity and the diagrams associated with curvature, general
relativity and gravity.

Since the metric is so important one may ask why one should study
manifolds without metrics. There are a number of reasons for doing
this. One is that many manifolds do not possess a metric. To be more
precise, since on every manifold one may construct a metric, one
should say it has no single preferred metric, i.e. one with physical
significance. For example phase space does not possess a
metric, neither do the natural manifolds for studying differential
equations. Even in general relativity, one may consider varying the
metric in order to derive Einstein's equations and the
stress-energy-momentum tensor. Transformation optics, by contrast, may
be considered as the study of spacetimes with two metrics, the real
and the optical metric.  As a result knowing about objects and
operations which do not require a metric is very useful. The other
reason for introducing differential forms before the metric is that,
if one has a metric, then one can pass between 1-forms and vectors. Thus
they both look very similar and it is difficult for the student to get
an intuition about differential forms as distinct from vectors.

One of the goals of this article is to introduce the equations of
electromagnetism and electrodynamics pictorially. Electromagnetism
lends itself to this approach as the objects of electromagnetism are
exterior differential forms, and most of electromagnetism does not
need a metric. Indeed the four macroscopic Maxwell's equations can be
written (and therefore drawn) without a metric. It is only when we wish
to give the additional equations which prescribe vacuum or the medium
that we need the metric.  This has lead to school of thought
\cite{hehl2012foundations} which suggest that the metric may be a
consequence of electromagnetism.

The challenge when depicting Maxwell's equations pictorially is that they are
equations on four dimensional spacetime and we lack 4-dimensional paper.
Unfortunately, were we blessed with such 4-dimensional paper, we would have to
live in a universe with 5 spatial dimensions and therefore 6-dimensional
spacetime. Whatever field the 5-dimensional equivalent of me would be attempting
to represent in pictures he/she/it would presumably be bemoaning the lack of
6-dimensional paper.

\vspace{1em}

This article is organised as follows. We start in section \ref{ch_DG}
with basic differential geometry, discussing manifolds, submanifolds,
scalar and vector fields. In section \ref{ch_Orientation} we introduce
both untwisted and twisted orientation, dealing first with
orientations on submanifolds, since this is the easiest to transfer to
other objects like forms and vectors. We then introduce closed
exterior differential forms, section \ref{ch_CForms}, showing how one
can add them and take the wedge product. We then talk about
integration and hence conservation laws. This leads to a discussion
about manifolds with ``holes'' in them. Finally in this section we
discuss the challenge of representing forms in four dimensions. In
section \ref{ch_NonClosed_Forms} we discuss non-closed forms and how
to visualise and integrate them. In section \ref{ch_other_ops} we
summarise a number of additional operation one can do with vectors and
forms.  The metric is
introduced in section \ref{ch_Metric} where we show how we can use
it to measure the length of a vector, angles between vectors as
well the the metric and Hodge dual. Section \ref{ch_Maps} is about
smooth maps which shows how forms can be pulled from one manifold to
another.  Finally we reach our goal in section \ref{ch_Maxwell} which
shows how to picture Maxwell's equations and the Lorentz force
equation. We start with static electromagnetic fields as these are the
easiest to visualise and then attempt to draw the full four
dimensional pictures. We then return to three dimensional pictures but
this time one time and two spatial dimensions. We then conclude in
section \ref{ch_Conclusion} with a discussion of further pictures one
may contemplate.

All figures are available in colour online.

\section{Manifolds, scalar and vector fields}
\label{ch_DG}


\subsection{Manifolds.}
\label{ch_DG_Man}

\begin{figure}[tb]
\begin{tikzpicture}[scale=2]
\draw  (-1,0) node {(a)} ;
\shade [ball color=green] (0,-.7) circle (1) ;
\draw  [very thick,yscale=0.4] (.5,0.4) arc (-50:-130:.77) ;  
\draw  (1.4,0) node {(b)} ;
\shade [ball color=green] (1.4,-.3) arc (120:60:2) 
(1.4,-.3) arc (-105:-75:3.8) ;
\draw  [very thick] (1.4,-.3) arc (-105:-75:3.8) ;  
\draw  [very thick,red] (1.7,0.14-.3) arc (-100:-80:4) ;  
\draw (1.4,-1.2) node {(c)} ;
\shade [ball color=green] (1.4,-1.5) arc (110:70:3) 
(1.4,-1.5) arc (-100:-80:6) ;
\draw  [very thick,red] (1.4,-1.5) arc (-100:-80:6) ;
\end{tikzpicture}
\caption[Zooming in on a sphere.]{Zooming in on a sphere. As we go
  from (a) to (b) to (c), the patch of the sphere becomes flatter.}
\label{fig_Man_zoom_sphere}
\end{figure}
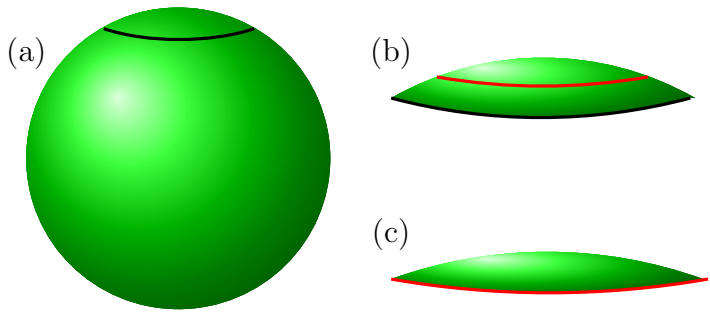
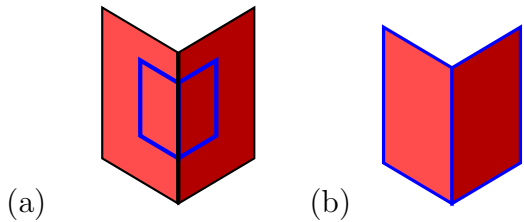
\begin{figure}[tb]
\begin{tikzpicture}[scale=2]
\draw  (-1,0) node {(a)} ;
\filldraw [thick,fill=red!70!white] (0,0) -- (-.5,.3) -- (-.5,1.3) --
(0,1) -- cycle ; 
\filldraw [thick,fill=red!70!black] (0,0) -- (.5,.3) -- (.5,1.3) --
(0,1) -- cycle ;
\draw [ultra thick,blue,shift={(0,.3)},scale=0.5] (0,0) -- (.5,.3) --
(.5,1.3) -- (0,1) -- cycle -- (-.5,.3) -- (-.5,1.3) -- (0,1) -- cycle;
\draw [ultra thick] (0,0) -- (0,1) ;
\draw  (1,0) node {(b)} ;
\filldraw [very thick,fill=red!70!white,draw=blue,scale=0.9] (2,0) --
+(-.5,.3) -- +(-.5,1.3) -- +(0,1) -- cycle ; 
\filldraw [very thick,fill=red!70!black,draw=blue,scale=0.9] (2,0) --
+(.5,.3) -- +(.5,1.3) -- +(0,1) -- cycle ;

\end{tikzpicture}
\caption[Zooming in on a corner.]{Zooming in on a corner. As we go
  from (a) to (b), the patch does not get any flatter. This is not a
  manifold.}
\label{fig_Man_zoom_corner}
\end{figure}


\begin{figure}[tb]
\begin{subfigure}[b]{0.48\figsize}
\centering
\begin{tikzpicture}
\fill (1,1) circle [radius=0.05] ;
\fill (0,0) circle [radius=0.0] ;
\end{tikzpicture}
\caption{A 0-dimensional manifold: a point.}
\end{subfigure}
\quad
\begin{subfigure}[b]{0.48\figsize}
\begin{tikzpicture}
\draw [ultra thick, blue] (0,0) .. controls (1,1) and (2,-1) .. (3,0) ;
\end{tikzpicture}
\caption{A 1-dimensional manifold: a curve.}
\end{subfigure}
\quad
\begin{subfigure}[b]{0.48\figsize}
\begin{tikzpicture}
\filldraw [ultra thick,fill=red!30] (0,0) .. controls (1,1) and (2,-1) .. (3,0) 
         .. controls (2.8,1.2) and (3.7,2) .. (3,3)
         .. controls (2,2.3) and (1,3.7) .. (0,3) 
        .. controls (.8,2)and (-0.7,1.2) .. (0,0)  ;
\draw (1,0.2) .. controls (0.2,1)  and (1.8,2) .. (1,3.13) ;
\draw (2,-0.2) .. controls (1.4,1) and (2.8,2) .. (2,2.85) ;
\draw (-.15,1) .. controls (1,2) and (2,0) .. (3.07,1) ;
\draw (0.2,2) .. controls (1,3) and (2,1) .. (3.3,2) ;
\end{tikzpicture}
\caption{A 2-dimensional manifold: a surface.}
\end{subfigure}
\quad
\begin{subfigure}[b]{0.48\figsize}
\begin{tikzpicture}[scale=0.7]
\draw [ultra thick] (0,0) .. controls (1,1) and (2,-1) .. (3,0) ;
\draw [ultra thick] (0,3) .. controls (1,3.7) and (2,2.3) .. (3,3) ;
\draw [ultra thick] (0,0) .. controls (-0.7,1.2) and (.8,2) .. (0,3) ;
\draw [ultra thick] (3,0) .. controls (2.8,1.2) and (3.7,2) .. (3,3) ;
\draw [ultra thick] (3,3) .. controls (3.2,3.7) and (4.,3.3) .. (4.5,4) ;
\draw [ultra thick] (0,3) .. controls (0.2,3.7) and (1.,3.3) .. (1.5,4) ;
\draw [ultra thick] (3,0) .. controls (3.2,0.7) and (4.,0.3) .. (4.5,1) ;
\draw [ultra thick] (1.5,4) .. controls (2.5,4.3) and (3.5,3.3) .. (4.5,4) ;
\draw [ultra thick] (4.5,1) .. controls (4.4,2.2) and (5.2,2) .. (4.5,4) ;
\draw [dashed] (0,0) .. controls (0.2,0.7) and (1.,0.3) .. (1.5,1) ;
\draw [dashed] (1.5,1) .. controls (2.5,1.3) and (3.5,0.3) .. (4.5,1) ;
\draw [dashed] (1.5,1) .. controls (1.4,2.2) and (2.2,2) .. (1.5,4) ;
\end{tikzpicture}
\caption{A 3-dimensional manifold: a volume.}
\end{subfigure}
\caption{0, 1, 2 and 3-dimensional manifolds.}
\label{fig_Manifolds}
\end{figure}
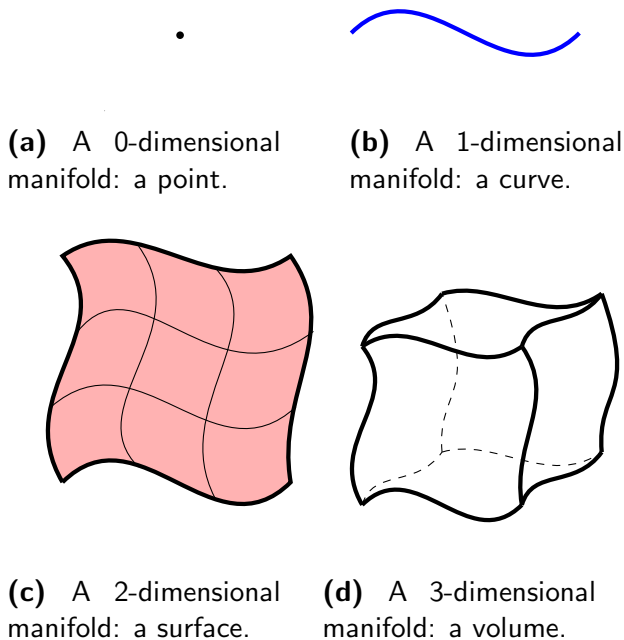


Manifolds are the fundamental object in differential geometry on which
all other object exist. These are particularly difficult to define
rigorously and it usually take the best part of an undergraduate
mathematics lecture course to define. This means that the useful tools
discussed here are relegated to the last few lectures.

Intuitively we should think of a manifold as a shape, like a sphere or
torus, which is smooth. That is, it does not have any corners. If we
take a small patch of a sphere and enlarge it, then the result looks
flat, like a small patch of the plane, figure
\ref{fig_Man_zoom_sphere}. By contrast if we take a small patch of a
corner then, no mater how large we made it, it would still never be
flat, figure \ref{fig_Man_zoom_corner}.

In figure \ref{fig_Manifolds} are
sketched manifolds with 0, 1, 2 and 3 dimensions.
0-dimensional manifolds are simply points, 1-dimensional
manifolds are curves, 2-dimensional manifolds are surfaces and
3-dimensional manifolds are volumes. 

Topologically, there are only
two types of 1 dimensional manifolds, those which are lines and
those which are loops. Higher dimensional manifolds are more interesting.

\begin{figure*}[tb]
\centering
\begin{subfigure}[t]{0.22\textwidth}
\centering
\begin{tikzpicture}
\draw [thick] (0,0) to[out=70,in=-100] (0,1) 
   to[out=80,in=-80]  (0,2) to[out=100,in=-90] (0,3) ;
\draw [ultra thick,green!70!black] (0,1) 
   to[out=80,in=-80]  (0,2)  ;
\fill [blue!70!black] (0,1) circle (.1) (0,2) circle (.1) ;
\end{tikzpicture}
\caption{The interval (green) is a compact 1-dim submanifold of a 1-dim
  manifold. The boundary of the submanifold are the two points, this
  is a compact 0-dim  submanifold (blue).}
\label{fig_subman_01}
\end{subfigure}
\quad
\begin{subfigure}[t]{0.22\textwidth}
\begin{tikzpicture}
\fill [green!70!black] (0,0) 
  to[out=30,in=210]  (3,0)
  .. controls (2.95,.5) .. (3.05,1) 
  .. controls (2,0) and  (1,2).. (-.15,1)
  .. controls (-.23,.5) .. (0,0) ; 
\draw [thick] (0,0) to[out=30,in=210] (3,0) ;
\draw [thick] (0,3) .. controls (1,3.7) and (2,2.3) .. (3,3) ;
\draw [thick] (0,0) .. controls (-0.7,1.2) and (.8,2) .. (0,3) ;
\draw [thick] (3,0) .. controls (2.8,1.2) and (3.7,2) .. (3,3) ;
\draw [color=blue!70!black, ultra thick] 
   (-.15,1) .. controls (1,2) and (2,0) .. (3.05,1) ;
\end{tikzpicture}
\vspace{-2.5em}\ 
\caption{The lower portion is a non compact 2-dim submanifold (green) 
of a 2-dim manifold. The boundary of this submanifold is the 
non compact 1-dim submanifold (blue).}
\label{fig_subman_12_open}
\end{subfigure}
\quad
\begin{subfigure}[t]{0.22\textwidth}
\begin{tikzpicture}
\draw [thick](0,0) .. controls (1,1) and (2,-1) .. (3,0) ;
\draw [thick](0,3) .. controls (1,3.7) and (2,2.3) .. (3,3) ;
\draw [thick](0,0) .. controls (-0.7,1.2) and (.8,2) .. (0,3) ;
\draw [thick](3,0) .. controls (2.8,1.2) and (3.7,2) .. (3,3) ;
\filldraw [draw=blue,ultra thick,fill=green!70!black] 
(1.5,1.5) circle [radius=0.8] ;
\end{tikzpicture}
\vspace{-2.5em}\ \caption{The disc (green) is a compact 2-dim
  submanifold of a 2-dim manifold. The boundary of the disk is the
  1-dim circle (blue), a compact submanifold.}
\label{fig_subman_12_closed}
\end{subfigure}
\quad
\begin{subfigure}[t]{0.22\textwidth}
\begin{tikzpicture}[scale=.8]
\fill [green!70!black] (2.0,0.5) 
.. controls (2.0,2) .. (2.5,3.5) 
-- (3.2,3.3) to[out=-55,in=95] (3.2,.3) 
--  (2.0,0.5) ;

\draw [color=blue,very thick] (2.0,0.5) .. controls (2.0,2) .. (2.5,3.5) ;
\draw [thick](0,0) .. controls (1,1) and (2,-1) .. (3,0) ;
\draw [thick](0,3) .. controls (1,3.7) and (2,2.3) .. (3,3) ;
\draw [thick](0,0) .. controls (-0.7,1.2) and (.8,2) .. (0,3) ;
\draw [thick](3,0) to[out=95,in=-55] (3,3) ;
\draw [thick](3,3) .. controls (3.2,3.7) and (4.,3.3) .. (4.5,4) ;

\draw [thick](0,3) .. controls (0.2,3.7) and (1.,3.3) .. (1.5,4) ;
\draw [thick](3,0) .. controls (3.2,0.7) and (4.,0.3) .. (4.5,1) ;
\draw [thick](1.5,4) .. controls (2.5,4.3) and (3.5,3.3) .. (4.5,4) ;
\draw [thick](4.5,1) .. controls (4.4,2.2) and (5.2,2) .. (4.5,4) ;
\draw [dashed] (0,0) .. controls (0.2,0.7) and (1.,0.3) .. (1.5,1) ;
\draw [dashed] (1.5,1) .. controls (2.5,1.3) and (3.5,0.3) .. (4.5,1) ;
\draw [dashed] (1.5,1) .. controls (1.4,2.2) and (2.2,2) .. (1.5,4) ;
\end{tikzpicture}
\vspace{-2.5em}\ 
\caption{The half plane is a non compact 2-dim submanifold  (green) 
of a 3-dim manifold. Its boundary is the non compact 1-dim curve (blue).}
\label{fig_subman_13_open}
\end{subfigure}
\quad
\begin{subfigure}[t]{0.22\textwidth}
\begin{tikzpicture}[scale=0.8]
\draw (0,0) .. controls (1,1) and (2,-1) .. (3,0) ;
\draw (0,3) .. controls (1,3.7) and (2,2.3) .. (3,3) ;
\draw (0,0) .. controls (-0.7,1.2) and (.8,2) .. (0,3) ;
\draw (3,0) .. controls (2.8,1.2) and (3.7,2) .. (3,3) ;
\draw (3,3) .. controls (3.2,3.7) and (4.,3.3) .. (4.5,4) ;
\draw (0,3) .. controls (0.2,3.7) and (1.,3.3) .. (1.5,4) ;
\draw (3,0) .. controls (3.2,0.7) and (4.,0.3) .. (4.5,1) ;
\draw (1.5,4) .. controls (2.5,4.3) and (3.5,3.3) .. (4.5,4) ;
\draw (4.5,1) .. controls (4.4,2.2) and (5.2,2) .. (4.5,4) ;
\draw [dashed] (0,0) .. controls (0.2,0.7) and (1.,0.3) .. (1.5,1) ;
\draw [dashed] (1.5,1) .. controls (2.5,1.3) and (3.5,0.3) .. (4.5,1) ;
\draw [dashed] (1.5,1) .. controls (1.4,2.2) and (2.2,2) .. (1.5,4) ;
\fill [color=blue,opacity=0.5] (0.0,1.5) .. controls (1,1) and (2,2) .. (3,1) 
          .. controls (3.5,1.2) and (3.9,2) .. (4.65,2) 
          .. controls (3.8,2.5) and (2.5,2) .. (1.8,2.5) 
          .. controls (1.2,2) and (0.4,2.2) .. cycle ;
\end{tikzpicture}
\vspace{-2.5em} \caption{A non compact 2-dim submanifold (blue) of a
  3-dim manifold. It is the boundary of the non compact 3-dim
  submanifold consisting of it and all points above it.}
\label{fig_subman_23_open}
\end{subfigure}
\quad
\begin{subfigure}[t]{0.22\textwidth}
\begin{tikzpicture}[scale=0.8]
\draw (0,0) .. controls (1,1) and (2,-1) .. (3,0) ;
\draw (0,3) .. controls (1,3.7) and (2,2.3) .. (3,3) ;
\draw (0,0) .. controls (-0.7,1.2) and (.8,2) .. (0,3) ;
\draw (3,0) .. controls (2.8,1.2) and (3.7,2) .. (3,3) ;
\draw (3,3) .. controls (3.2,3.7) and (4.,3.3) .. (4.5,4) ;
\draw (0,3) .. controls (0.2,3.7) and (1.,3.3) .. (1.5,4) ;
\draw (3,0) .. controls (3.2,0.7) and (4.,0.3) .. (4.5,1) ;
\draw (1.5,4) .. controls (2.5,4.3) and (3.5,3.3) .. (4.5,4) ;
\draw (4.5,1) .. controls (4.4,2.2) and (5.2,2) .. (4.5,4) ;
\draw [dashed] (0,0) .. controls (0.2,0.7) and (1.,0.3) .. (1.5,1) ;
\draw [dashed] (1.5,1) .. controls (2.5,1.3) and (3.5,0.3) .. (4.5,1) ;
\draw [dashed] (1.5,1) .. controls (1.4,2.2) and (2.2,2) .. (1.5,4) ;
\shade [ball color=blue] (2,2) circle [radius=1cm];
\end{tikzpicture}
\vspace{-2.5em}
\caption{A compact 2-dim submanifold, the surface of the sphere,
  of a 3-dim manifold. It is the boundary of the compact 3-dim ball.}
\label{fig_subman_23_sphere}
\end{subfigure}
\quad
\begin{subfigure}[t]{0.25\textwidth}
\begin{tikzpicture}[scale=.8]
\draw (0,0) .. controls (1,1) and (2,-1) .. (3,0) ;
\draw (0,3) .. controls (1,3.7) and (2,2.3) .. (3,3) ;
\draw (0,0) .. controls (-0.7,1.2) and (.8,2) .. (0,3) ;
\draw (3,0) .. controls (2.8,1.2) and (3.7,2) .. (3,3) ;
\draw (3,3) .. controls (3.2,3.7) and (4.,3.3) .. (4.5,4) ;
\draw (0,3) .. controls (0.2,3.7) and (1.,3.3) .. (1.5,4) ;
\draw (3,0) .. controls (3.2,0.7) and (4.,0.3) .. (4.5,1) ;
\draw (1.5,4) .. controls (2.5,4.3) and (3.5,3.3) .. (4.5,4) ;
\draw (4.5,1) .. controls (4.4,2.2) and (5.2,2) .. (4.5,4) ;
\draw [dashed] (0,0) .. controls (0.2,0.7) and (1.,0.3) .. (1.5,1) ;
\draw [dashed] (1.5,1) .. controls (2.5,1.3) and (3.5,0.3) .. (4.5,1) ;
\draw [dashed] (1.5,1) .. controls (1.4,2.2) and (2.2,2) .. (1.5,4) ;

\draw (2,2) node {\tikz[scale=0.5] {
\shadedraw[ball color=blue,yscale=0.7] (0,0) circle (2) ;
\filldraw[fill=white,yscale=0.7] (0.7,0.2) arc (20:160:0.7) ;
\fill[white,yscale=0.7] (0.7,0.2) arc (-50:-130:1) ;
\draw[yscale=0.7] (1,0.6) arc (-20:-160:1) ;
}} ;
\end{tikzpicture}
\vspace{-2.em}
\caption{A compact 2-dim submanifold, the surface of a torus.
It is the boundary of the compact 3-dim solid torus.}
\label{fig_subman_23_torus}
\end{subfigure}

\caption[A selection of submanifolds and their boundaries]
{A selection of submanifolds and their boundaries. Observe that the
  boundaries is itself a submanifold, but it does not have a
  boundary. This  implies ate the boundary must extend to infinity or
close in on itself.}
\label{fig_Man_submanifolds}
\end{figure*}
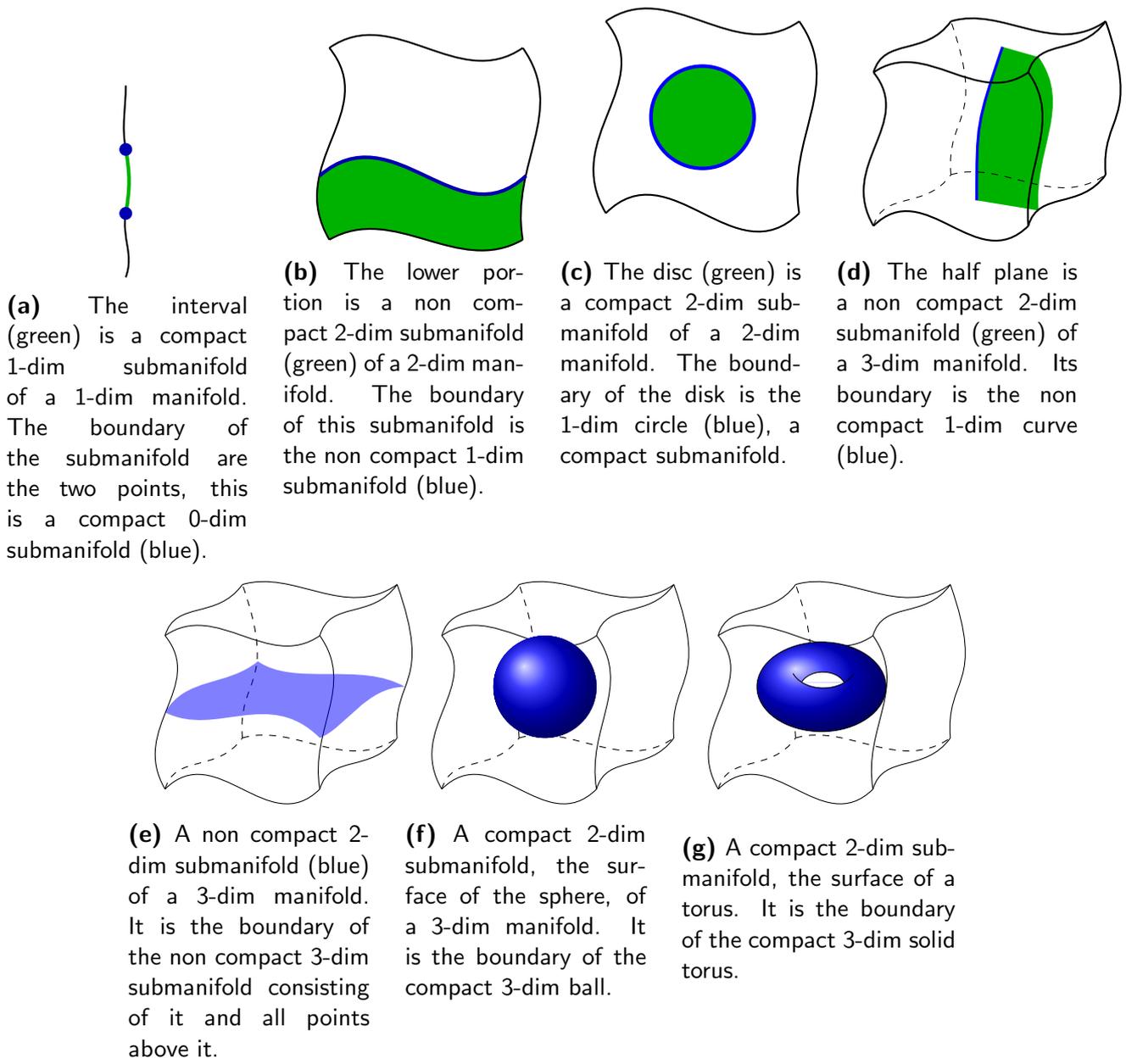

\subsection{Submanifolds.}
\label{ch_DG_SMan}

A submanifold is a manifold contained inside a larger manifold, called
the \defn{embedding manifold}. The
submanifold can be either the same dimension or a lower dimension than
the embedding manifold. Examples of embedded submanifolds are given in
figure \ref{fig_Man_submanifolds}

\begin{figure}
\centering
\begin{subfigure}[b]{0.48\figsize}
\centering
\begin{tikzpicture}[scale=0.8]
\draw [thick,color=blue] (0,0) to [out=90,in=135] (2,0) 
to [out=-45,in=-90] (4,0) to [out=90,in=-90] (0,0) ;
\end{tikzpicture}
\caption{A self intersecting curve.}
\label{fig_Man_patho_selfinter}
\end{subfigure}
\quad
\begin{subfigure}[b]{0.48\figsize}
\centering
\begin{tikzpicture}[scale=0.8]
\draw [->,thick,color=blue] (0,0) to [out=0,in=-90] (4,0) 
to [out=90,in=90] (2,-0.4) ;
\end{tikzpicture}
\caption{A curve approaching itself.}
\label{fig_Man_patho_approach}
\end{subfigure}
\caption{Pathological examples which are not submanifolds.}
\label{fig_Man_patho_submanifolds}
\end{figure}
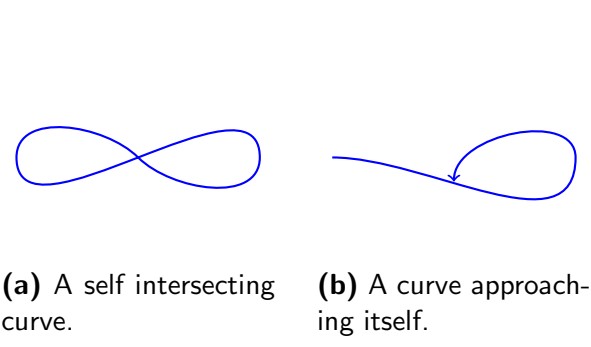

There are all kinds of bizarre pathologies which can take place when
one tries to embed one manifold into another, such as self
intersection, figure
\ref{fig_Man_patho_selfinter},
or the a manifold approaching itself
figure \ref{fig_Man_patho_approach}.
However we will avoid
these cases.

Submanifolds may either be compact or non compact. A compact
submanifold is bounded. By contrast a non compact submanifold extends
all the way to ``infinity'' at least in some direction. In figure
\ref{fig_Man_submanifolds} we see that examples \ref{fig_subman_01},
\ref{fig_subman_12_closed}, \ref{fig_subman_23_sphere} and
\ref{fig_subman_23_torus} are all compact, whereas the others are not
compact. A compact submanifold may either have a boundary, as in
section \ref{ch_DG_Bdd} below, or close in on itself like a sphere or torus.

\subsection{Boundaries of submanifold}
\label{ch_DG_Bdd}

Unlike the embedding manifold, it is also useful to consider
submanifolds with boundaries. See figure \ref{fig_Man_submanifolds}.
Both compact and non compact submanifolds can have boundaries.  The
importance of the boundary is when doing integration, in particular
the use of Stokes' theorem. The boundary of a submanifold is an
intuitive concept. For example the boundary of a 3-dimensional ball is
its surface, the 2-dimensional sphere, figure
\ref{fig_subman_23_sphere}.  We assume that the boundary of an
$n$-dimensional submanifold is an $(n-1)$-dimensional submanifold. We
observe that the boundary of the boundary vanishes.


\subsection{Scalar fields.}
\label{ch_DG_Scalar}

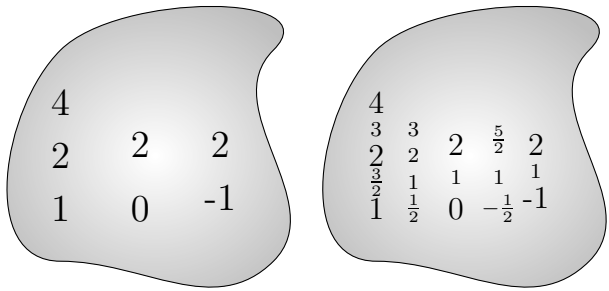
\begin{figure}
\begin{subfigure}[b]{0.48\figsize}
\begin{tikzpicture}[scale=0.7]
\shadedraw [inner color=white, outer color = lightgray] 
(0,0) to [out=0,in=-135] (4,0) to [in=-135] (4,4) to [in=45]
(0,4) to [out=-135,in=180] cycle ;
\draw (0,1) node {\large 1} ;
\draw (0,2) node {\large 2} ;
\draw (0,3) node {\large 4} ;
\draw (1.5,1) node {\large 0} ;
\draw (1.5,2.2) node {\large 2} ;
\draw (3,1.2) node {\large -1} ;
\draw (3,2.2) node {\large 2} ;
\end{tikzpicture}
\caption{Some of the values of the scalar field are shown.}
\end{subfigure}
\quad
\begin{subfigure}[b]{0.48\figsize}
\begin{tikzpicture}[scale=0.7]
\shadedraw [inner color=white, outer color = lightgray] 
(0,0) to [out=0,in=-135] (4,0) to [in=-135] (4,4) to [in=45]
(0,4) to [out=-135,in=180] cycle ;
\draw (0,1) node { 1} ;
\draw (0,2) node { 2} ;
\draw (0,3) node { 4} ;
\draw (1.5,1) node { 0} ;
\draw (1.5,2.2) node { 2} ;
\draw (3,1.2) node { -1} ;
\draw (3,2.2) node { 2} ;
\draw (0,1.5) node {\scriptsize $\tfrac32$} ;
\draw (0,2.5) node {\scriptsize 3} ;
\draw (0.7,1) node {\scriptsize $\tfrac12$} ;
\draw (0.7,1.5) node {\scriptsize 1} ;
\draw (0.7,2) node {\scriptsize 2} ;
\draw (0.7,2.5) node {\scriptsize 3} ;
\draw (1.5,1.6) node {\scriptsize 1} ;
\draw (2.3,1) node {\scriptsize $-\tfrac12$} ;
\draw (2.3,1.6) node {\scriptsize $1$} ;
\draw (2.3,2.3) node {\scriptsize $\tfrac52$} ;
\draw (3,1.7) node {\scriptsize 1} ;
\end{tikzpicture}
\caption{Intermediate values of the scalar field are also shown.}
\end{subfigure}
\caption{A representation of a scalar field.}
\label{fig_Scalar_field}
\end{figure}


\begin{figure}[tb]
\centering
\begin{tikzpicture}
\begin{axis}[grid=both,yticklabels={,,},xticklabels={,,},zticklabels={,,}]
  \addplot3[surf,shader=faceted] {1/(1+(x-1)*(x-1)+y*y)+1/(2+(x+4)*(x+4)+y*y)};
\end{axis}
\draw (-.3,2) node[rotate=90] {$f(x,y)$} ;
\draw (3,0) node[rotate=-10] {$x$} ;
\draw (6,.7) node[rotate=60] {$y$} ;
\end{tikzpicture}
\caption{A scalar field as a function of $x$ and $y$.}
\label{fig_Scalar_mountain}
\end{figure}
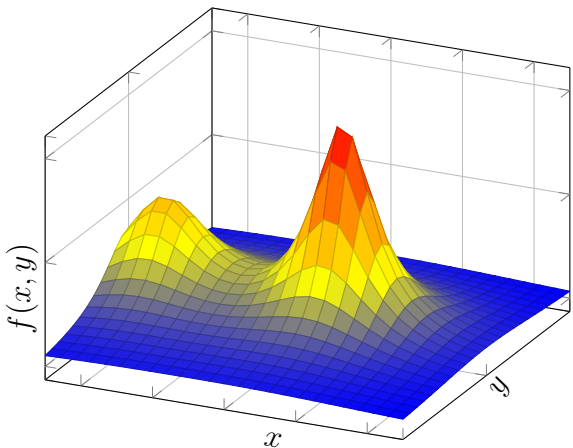


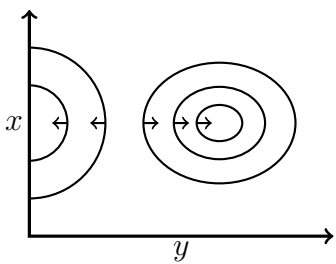
\begin{figure}[tb]
\centering
\begin{tikzpicture}
\draw [very thick,<->] (0,3) -- (0,0) -- (4,0) ; 
\draw (-.2,1.5) node {$x$} ;
\draw (2,-.2) node {$y$} ;
\draw [thick] (0,.5) arc (-90:90:1) (0,1) arc (-90:90:.5) ;
\draw [thick] (2.5,1.5/.8) [yscale=0.8] circle (1) circle (.3) circle (.6) ; 
\draw [thick,->] (.5,1.5) -- +(-.2,0) ; 
\draw [thick,->] (1,1.5) -- +(-.2,0) ; 
\draw [thick,->] (1.5,1.5) -- +(.2,0) ; 
\draw [thick,->] (1.9,1.5) -- +(.2,0) ; 
\draw [thick,->] (2.2,1.5) -- +(.2,0) ; 
\end{tikzpicture}
\caption[A scalar field depicted as contours.]{A scalar field depicted as contours. The arrows point up.}
\label{fig_Scalar_contours}
\end{figure}

Having defined manifolds, the simplest object that one may define on
a manifold is a scalar field.  This is simply a real number associated
with every point. See figure \ref{fig_Scalar_field}. The smoothness 
corresponds to nearby points having close values, also that these values
can be differentiated infinitely many times.  A scalar field on a
2-dimensional manifold may be pictured as function as in figure
\ref{fig_Scalar_mountain}. Alternatively one can depict a scalar field
in terms of its contours, figure \ref{fig_Scalar_contours}. We will
see below, section \ref{ch_CForms}, that these contours actually
depict the 1-form which is the the exterior derivative of the scalar field.


\subsection{Vector fields.}
\label{ch_DG_Vectors}

\begin{figure}[tb]
\centering
\begin{tikzpicture}[scale=0.8]
\shadedraw [inner color=white, outer color = lightgray] 
(0,0) to [out=0,in=-135] (4,0) to [in=-135] (4,4) to [in=45]
(0,4) to [out=-135,in=180] cycle ;
\draw (0,1) [blue,->]   -- +(1,0) ;
\draw (0,1.5) [blue,->] -- +(1,0.2) ;
\draw (0,2) [blue,->]   -- +(1,0.4) ;
\draw (0,2.5) [blue,->] -- +(1,0.6) ;
\draw (0,3) [blue,->]   -- +(1,0.7) ;
\draw (0,3.5) [blue,->] -- +(1,0.6) ;
\draw (0.7,1) [blue,->]   -- +(1,-0.2) ;
\draw (0.7,1.5) [blue,->] -- +(0.7,0.1) ;
\draw (0.7,2) [blue,->]   -- +(0.8,0.3) ;
\draw (0.7,2.5) [blue,->] -- +(0.9,0.6) ;
\draw (0.7,3) [blue,->]   -- +(1,0.8) ;
\draw (0.7,3.5) [blue,->] -- +(0.9,0.7) ;

\draw (1.5,1) [blue,->]   -- +(1,-0.4) ;
\draw (1.5,1.5) [blue,->] -- +(0.8,0.0) ;
\draw (1.5,2) [blue,->]   -- +(0.9,0.1) ;
\draw (1.5,2.5) [blue,->] -- +(1.0,0.4) ;
\draw (1.5,3) [blue,->]   -- +(1.1,0.5) ;
\draw (1.5,3.5) [blue,->] -- +(1.1,0.5) ;

\draw (2.3,1) [blue,->]   -- +(1,-0.4) ;
\draw (2.3,1.5) [blue,->] -- +(0.8,-0.2) ;
\draw (2.3,2) [blue,->]   -- +(0.9,-0.1) ;
\draw (2.3,2.5) [blue,->] -- +(1.0,0.0) ;
\draw (2.3,3) [blue,->]   -- +(1.1,0.3) ;
\draw (2.3,3.5) [blue,->] -- +(1.1,0.5) ;

\end{tikzpicture}
\caption{A representation of a vector field.}
\label{fig_Vector_field}
\end{figure}
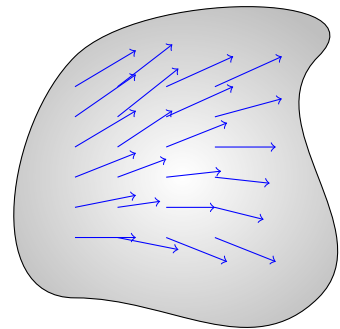

A vector field is a concept most familiar to people. A vector field
consists of an arrow, with a base point (bottom of the
stalk) at each point on the manifold, see figure  \ref{fig_Vector_field}.

Although we have only drawn a finite number of vectors, there is a
vector at each point. The smoothness of the vector field corresponds
to nearby vectors not differing too much.

To multiply a vector field by a scalar field, simply scale the length of
each vector by the value of the scalar field at the base point of
that vector. Of course if the scalar field is negative we reverse the
direction of the arrow. 
We are all familiar with the addition of vectors from school, via the
use of a parallelogram.

The vector presented in figure \ref{fig_Vector_field} are untwisted
vectors. We will see below in section \ref{ch_Orientation_Vectors}
that there also exist twisted vectors. 


\section{Orientation.}
\label{ch_Orientation}

We now discuss the tricky role that orientation has to play in differential
geometry, especially the difference between twisted and untwisted
orientations. However the
distinction is worth emphasising since the objects in electromagnetism are
necessarily twisted and others are necessarily untwisted.

The take home message is:
\begin{jgitemize}
\item
Orientation is important for integration and Stokes's theorem.
\item
The orientation usually only affects the overall sign. However it is
important when adding forms.
\item
The pictures are pretty.
\end{jgitemize}


\subsection{Orientation on Manifolds (internal)}


\begin{figure}[tb]
\centering
\begin{subfigure}[b]{0.22\figsize}
\centering
\begin{tikzpicture}
\fill [green!60!black] (0,0) circle (0.1) ; 
\draw (0,0.4) node {+} ;
\end{tikzpicture}
\quad
\begin{tikzpicture}
\fill [green!60!black] (0,0) circle (0.1) ; 
\draw (0,0.4) node {$-$} ;
\end{tikzpicture}
\caption{0-dim\\ manifold.}
\label{fig_ori_Man0}
\end{subfigure}
\quad
\begin{subfigure}[b]{0.22\figsize}
\centering
\begin{tikzpicture}
\draw[ultra thick,->-=.7,blue] (0,0) to[out=80,in=-110] (0,1) to[out=80,in=-110] (0,2) ; 
\end{tikzpicture}
\quad
\begin{tikzpicture}
\draw[ultra thick,-<-=.3,blue] (0,0) to[out=80,in=-110] (0,1) to[out=80,in=-110] (0,2) ; 
\end{tikzpicture}
\caption{1-dim\\ manifold.}
\label{fig_ori_Man1}
\end{subfigure}
\begin{subfigure}[b]{0.4\figsize}
\centering
\begin{tikzpicture}
\draw (0,0) to[out=80,in=-110] (0,2) 
   to[out=10,in=-170] (1,2) 
   to[out=-80,in=100] (1,0) 
   to[out=-170,in=10] cycle ;
\draw[very thick,red!60!black,->] (0.7,1) arc(0:130:0.25) ; 
\draw[very thick,red!60!black] (0.7,1) arc(0:260:0.25) ; 
\end{tikzpicture}
\ 
\begin{tikzpicture}
\draw (0,0) to[out=80,in=-110] (0,2) 
   to[out=10,in=-170] (1,2) 
   to[out=-80,in=100] (1,0) 
   to[out=-170,in=10] cycle ;
\draw[very thick,red!60!black,-<] (0.7,1) arc(0:130:0.25) ; 
\draw[very thick,red!60!black] (0.7,1) arc(0:260:0.25) ; 
\end{tikzpicture}
\caption{2-dim manifold.}
\label{fig_ori_Man2}
\end{subfigure}
\begin{subfigure}[b]{0.40\textwidth}
\centering
\begin{tikzpicture}[scale=0.5]
\draw[black!50] (0,0) .. controls (1,1) and (2,-1) .. (3,0) ;
\draw[black!50] (0,3) .. controls (1,3.7) and (2,2.3) .. (3,3) ;
\draw[black!50] (0,0) .. controls (-0.7,1.2) and (.8,2) .. (0,3) ;
\draw[black!50] (3,0) .. controls (2.8,1.2) and (3.7,2) .. (3,3) ;
\draw[black!50] (3,3) .. controls (3.2,3.7) and (4.,3.3) .. (4.5,4) ;
\draw[black!50] (0,3) .. controls (0.2,3.7) and (1.,3.3) .. (1.5,4) ;
\draw[black!50] (3,0) .. controls (3.2,0.7) and (4.,0.3) .. (4.5,1) ;
\draw[black!50] (1.5,4) .. controls (2.5,4.3) and (3.5,3.3) .. (4.5,4) ;
\draw[black!50] (4.5,1) .. controls (4.4,2.2) and (5.2,2) .. (4.5,4) ;
\draw[dashed,black!50] (0,0) .. controls (0.2,0.7) and (1.,0.3) .. (1.5,1) ;
\draw[dashed,black!50] (1.5,1) .. controls (2.5,1.3) and (3.5,0.3) .. (4.5,1) ;
\draw[dashed,black!50] (1.5,1) .. controls (1.4,2.2) and (2.2,2) .. (1.5,4) ;

\draw[ultra thick,draw=white,double=blue,xshift=2cm]
(0.5,1) \foreach \x in {0,10,...,720}
{ to ( {0.5*cos(\x)}, {0.0025*\x+0.5*sin(\x)+1} ) } ;
\draw[ultra thick,draw=white,double=blue,xshift=2cm]
(-.5,0.0025*180+1) \foreach \x in {180,190,...,360}
{ to ( {0.5*cos(\x)}, {0.0025*\x+0.5*sin(\x)+1} ) } ;
\draw[ultra thick,draw=white,double=blue,xshift=2cm]
(-.5,0.0025*540+1) \foreach \x in {540,550,...,720}
{ to ( {0.5*cos(\x)}, {0.0025*\x+0.5*sin(\x)+1} ) } ;
\draw (2,-1) node {Right handed} ;
\end{tikzpicture}
\begin{tikzpicture}[scale=0.5]
\draw[black!50] (0,0) .. controls (1,1) and (2,-1) .. (3,0) ;
\draw[black!50] (0,3) .. controls (1,3.7) and (2,2.3) .. (3,3) ;
\draw[black!50] (0,0) .. controls (-0.7,1.2) and (.8,2) .. (0,3) ;
\draw[black!50] (3,0) .. controls (2.8,1.2) and (3.7,2) .. (3,3) ;
\draw[black!50] (3,3) .. controls (3.2,3.7) and (4.,3.3) .. (4.5,4) ;
\draw[black!50] (0,3) .. controls (0.2,3.7) and (1.,3.3) .. (1.5,4) ;
\draw[black!50] (3,0) .. controls (3.2,0.7) and (4.,0.3) .. (4.5,1) ;
\draw[black!50] (1.5,4) .. controls (2.5,4.3) and (3.5,3.3) .. (4.5,4) ;
\draw[black!50] (4.5,1) .. controls (4.4,2.2) and (5.2,2) .. (4.5,4) ;
\draw[dashed,black!50] (0,0) .. controls (0.2,0.7) and (1.,0.3) .. (1.5,1) ;
\draw[dashed,black!50] (1.5,1) .. controls (2.5,1.3) and (3.5,0.3) .. (4.5,1) ;
\draw[dashed,black!50] (1.5,1) .. controls (1.4,2.2) and (2.2,2) .. (1.5,4) ;

\draw[ultra thick,draw=white,double=blue,xscale=-1,xshift=-2cm]
(0.5,1) \foreach \x in {0,10,...,720}
{ to ( {0.5*cos(\x)}, {0.0025*\x+0.5*sin(\x)+1} ) } ;
\draw[ultra thick,draw=white,double=blue,xscale=-1,xshift=-2cm]
(-.5,0.0025*180+1) \foreach \x in {180,190,...,360}
{ to ( {0.5*cos(\x)}, {0.0025*\x+0.5*sin(\x)+1} ) } ;
\draw[ultra thick,draw=white,double=blue,xscale=-1,xshift=-2cm]
(-.5,0.0025*540+1) \foreach \x in {540,550,...,720}
{ to ( {0.5*cos(\x)}, {0.0025*\x+0.5*sin(\x)+1} ) } ;
\draw (2,-1) node {Left handed} ;
\end{tikzpicture}
\caption{3-dim manifold.}
\label{fig_ori_Man3}
\end{subfigure}

\caption[Internal orientations.]{Possible (internal) orientations on 0, 1, 2 and
  3 dimensional manifold.}
\label{fig_ori_Man}
\end{figure}
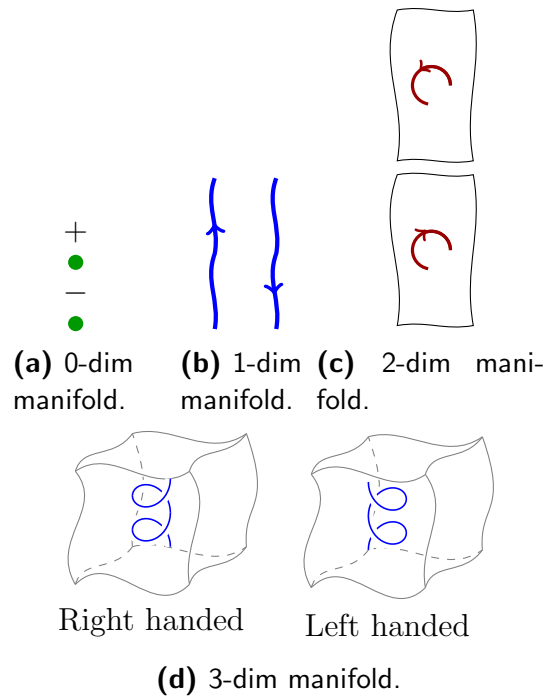

We start by showing the possible orientations on 
manifolds. 
Not all manifolds can have an orientation. If a
manifold can have an orientation it is called
orientable.  Examples of orientable manifolds include the sphere and
torus. Non-orientable manifolds include the Möbius strip and the Kline
bottle. Even if a manifold is orientable, one may
choose not to put an orientation on the it.
For manifolds of 0, 1, 2 or 3 dimensions the two
possible (internal) orientations are shown in figure
\ref{fig_ori_Man}. These are as follows:

\begin{jgitemize}
\item
An orientation for a 0-dimensional submanifold consists of either plus or
minus, figure \ref{fig_ori_Man0}.
However convention usually dictates plus.
\item
On a 1-dimensional manifold it is an arrow pointing one way or
another, figure \ref{fig_ori_Man1}.
\item
On a 2-dimensional manifold it is an arc pointing clockwise or
anticlockwise, figure \ref{fig_ori_Man2}.
\item
On a three dimensional manifold one considers a helix pointing one way
or another, figure \ref{fig_ori_Man3}. 
Interestingly, observe that unlike in the 1-dim and
2-dim cases, there is no need for an arrow on the helix.
\end{jgitemize}

\subsection{Orientation on submanifolds: Internal and external}


\begin{table}[tb]
\centering
\begin{tabular}{|l|l|@{}c@{}|@{}c@{}|@{\,}c@{\,}|@{}c@{}|}
\hline
\multicolumn {2}{|c|}{} & \multicolumn{4}{|c|}{\rule{0em}{1.1em}Embedding dimension}
\\\cline{3-6}
\multicolumn {2}{|c|}{} & \rule{0em}{1.1em} 0 & 1 & 2 & 3 
\\[.2em]\hline 
& \raisebox{2em}{0} & 
\begin{tikzpicture}
\fill [green] (0,0) circle (0.1) ; 
\draw (0,0.4) node[red] {+} ;
\fill [green] (0.5,0) circle (0.1) ; 
\draw (0.5,0.4) node[red] {$-$} ;
\end{tikzpicture}
& 
\begin{tikzpicture}
\draw [thick,black] (0,0) -- (0,1) ; 
\fill [ultra thick, green] (0,0.5) circle (0.1) ; 
\draw [ultra thick, red,->] (0,0.6) -- +(0,.3) ; 
\draw [thick, black] (0.5,0) -- +(0,1) ; 
\fill [ultra thick, green] (0.5,0.5) circle (0.1) ; 
\draw [ultra thick, red,->] (0.5,0.4) -- +(0,-.3) ; 
\end{tikzpicture}
& 
\begin{tikzpicture}[scale=1.3]
\draw [thick,black] (0,0) rectangle (.7,.7) ; 
\fill [ultra thick, green] (0.35,0.35) circle (0.1) ; 
\draw [very thick,red,->-=.7] (0.35,0.35) +(-10:.2) arc (-10:170:.2) ; 
\begin{scope}[shift={(1,0)}]
\draw [thick,black] (0,0) rectangle (.7,.7) ; 
\fill [ultra thick, green] (0.35,0.35) circle (0.1) ; 
\draw [very thick,red,-<-=.7] (0.35,0.35) +(-10:.2) arc (-10:170:.2) ; 
\end{scope}
\end{tikzpicture}
\ 
&
\begin{tikzpicture}[scale=.7]
\draw[ultra thick,draw=white,double=red,xshift=2cm]
(0.5,1) \foreach \x in {0,10,...,720}
{ to ( {0.5*cos(\x)}, {0.0025*\x+0.5*sin(\x)+1} ) } ;
\draw[ultra thick,draw=white,double=red,xshift=2cm]
(-.5,0.0025*180+1) \foreach \x in {180,190,...,360}
{ to ( {0.5*cos(\x)}, {0.0025*\x+0.5*sin(\x)+1} ) } ;
\fill [green] (2,1.9) circle (0.3) ; 
\draw[ultra thick,draw=white,double=red,xshift=2cm]
(-.5,0.0025*540+1) \foreach \x in {540,550,...,720}
{ to ( {0.5*cos(\x)}, {0.0025*\x+0.5*sin(\x)+1} ) } ;
\end{tikzpicture}
%
\begin{tikzpicture}[yscale=.7,xscale=-.7]
\draw[ultra thick,draw=white,double=red,xshift=2cm]
(0.5,1) \foreach \x in {0,10,...,720}
{ to ( {0.5*cos(\x)}, {0.0025*\x+0.5*sin(\x)+1} ) } ;
\draw[ultra thick,draw=white,double=red,xshift=2cm]
(-.5,0.0025*180+1) \foreach \x in {180,190,...,360}
{ to ( {0.5*cos(\x)}, {0.0025*\x+0.5*sin(\x)+1} ) } ;
\fill [green] (2,1.9) circle (0.3) ; 
\draw[ultra thick,draw=white,double=red,xshift=2cm]
(-.5,0.0025*540+1) \foreach \x in {540,550,...,720}
{ to ( {0.5*cos(\x)}, {0.0025*\x+0.5*sin(\x)+1} ) } ;
\end{tikzpicture}
\\ \cline{2-6}
\raisebox{-3em}[0em][0em]{\rotatebox{90}{Submanifold dimension}}
 & \raisebox{2em}{1} & &
\begin{tikzpicture}[yscale=1.5]
\draw [thick,black] (0,0) -- (0,1) ; 
\draw [ultra thick, green] (0,0.3) -- +(0,.5) ; 
\draw [red] (-.2,0.5) node {$+$} ;
\draw [thick,black] (.5,0) -- +(0,1) ; 
\draw [ultra thick, green] (0.5,0.3) -- +(0,.5) ; 
\draw [red] (.7,0.5) node {$-$} ;
\end{tikzpicture}
&
\begin{tikzpicture}[scale=1.3]
\draw [thick,black] (0,0) rectangle (.7,.7) ; 
\draw [ultra thick, green] (.35,0) -- +(0,.7) ; 
\draw [very thick,red,->] (0.35,.35) -- +(-.3,0) ; 
\begin{scope}[shift={(1,0)}]
\draw [thick,black] (0,0) rectangle (.7,.7) ; 
\draw [ultra thick, green] (.35,0) -- +(0,.7) ; 
\draw [very thick,red,->] (0.35,.35) -- +(.3,0) ; 
\end{scope}
\end{tikzpicture}
&
\begin{tikzpicture}
\draw [ultra thick, green] (.25,0) -- +(0,1) ; 
\draw [very thick,white,yscale=.6,double=white] (0.25,1) +(70:.2) arc (70:-250:.2) ; 
\draw [very thick,red,yscale=.6,->-=.5] (0.25,1) +(70:.2) arc (70:-250:.2) ; 
\draw [ultra thick, green] (1.25,0) -- +(0,1) ; 
\draw [very thick,white,yscale=.6,double=white] (1.25,1) +(70:.2) arc (70:-250:.2) ; 
\draw [very thick,red,yscale=.6,-<-=.7] (1.25,1) +(70:.2) arc (70:-250:.2) ; 
\end{tikzpicture}
\\\cline{2-6}
& \raisebox{1.5em}{2}\rule{0em}{3.5em} & &
&
\begin{tikzpicture}[scale=1.3]
\draw [thick,black] (0,0) rectangle (.7,.7) ; 
\filldraw [thick,fill=green] (.35,0.35) circle (.2) ; 
\draw  (.35,0.35) node[red] {$+$} ;
\begin{scope}[shift={(1,0)}]
\draw [thick,black] (0,0) rectangle (.7,.7) ; 
\filldraw [thick,fill=green] (.35,0.35) circle (.2) ; 
\draw  (.35,0.35) node[red] {$-$} ;
\end{scope}
\end{tikzpicture}
&
\begin{tikzpicture}
\filldraw [thick,fill=green,xscale=.6] (0,0) circle (.4) ; 
\draw [very thick,white,xscale=.6,double=white] (0,0) -- +(1,0) ; 
\draw [very thick,red,xscale=.6,->] (0,0) -- +(.8,0) ; 
\draw [very thick,red,xscale=.6,->] (2,0) -- +(-.8,0) ; 
\filldraw [thick,fill=green,xscale=.6] (2,0) circle (.4) ; 
\end{tikzpicture}
\\\cline{2-6}
& \raisebox{1.em}{3}\rule{0em}{2.5em} & & & &
\begin{tikzpicture}
\shadedraw [thick,ball color=green] (.25,0.5) circle (.35) ; 
\draw  (.25,0.5) node[red] {$+$} ;
\shadedraw [thick,ball color=green] (1.25,0.5) circle (.35) ; 
\draw  (1.25,0.5) node[red] {$-$} ;
\end{tikzpicture}
\\\hline

\end{tabular}

\caption[External orientations.]{Possible external orientations on a
  0, 1, 2 and 3-dim submanifold of a 0, 1, 2 and 3-dim embedding
  manifold.}
\label{tab_ori_twist_SubMan}
\end{table}

Again  submanifolds may be orientable or non-orientable. However with
orientable manifold one now has a choice as to which type of
orientations to chose: internal or external. With a manifold there is
no concept of an external orientation because there is no embedding
manifold in which to put the external orientation. 

The internal orientation of submanifolds are exactly the same as those
for manifolds, as shown in figure \ref{fig_ori_Man}, with the type of
the orientation: sign, arrow, arc or helix, depending on the dimension
of the submanifold.

By contrast external orientations depend not just on the dimension of
the submanifold but also on the dimension of the embedding manifold.
They are similar, plus or minus sign, arrow, arc or helix, as the
internal orientations. However in this case it takes place in the
embedding external to the submanifold.  In table
\ref{tab_ori_twist_SubMan} we see examples of external orientations
for 0, 1, 2, 3 dimensional submanifolds embedded into 0, 1, 2, 3
dimensional manifolds. As we see, the type of the external orientation
depends on the difference in the dimension of the embedding manifold
and submanifold. Thus the external orientation is a plus or minus sign
if the two manifolds have the same dimension. It is an arrow if
they differ by 1 dimension, an arc if they differ by 2 dimensions
and a helix if they differ by 3 dimensions.


\begin{figure}[tb]
\centering
\begin{tikzpicture}
\draw[ultra thick] (0,0) to[out=80,in=-110] (0,1) ; 
\draw[ultra thick] (0,1) to[out=80,in=-110] (0,2) ; 

\draw[ultra thick,yscale=0.5,draw=white,double=green!50!black] 
(-0.09,2) arc (260:-80:0.4) ; 
\draw[very thick,draw=green!50!black] 
(0.27,1.17) -- +(0.1,0.1) -- +(0.2,0) ; 

\draw[very thick,yshift=-0.7cm,xshift=0.62cm,rotate={30},
     yscale=0.6,xscale=0.8,draw=white,double=green!50!black] 
(-0.09,2) arc (260:-80:0.4) ; 
\draw[very thick,yshift=-0.6cm,xshift=-0.2cm,draw=green!50!black] 
(0.27,1.17) -- +(0.1,0.1) -- +(0.2,0) ; 
\end{tikzpicture}
\qquad
\begin{tikzpicture}
\filldraw[ultra thick,draw=black,fill=blue!40!white] (-0.8,-1.6) 
to [out=110,in=-70] +(-.5,2) 
to [out=70,in=-100] +(1.5,1) 
to [out=-60,in=110] +(0.5,-2) 
to [out=-120,in=70] cycle ;

\draw[very thick,draw=white,double=white] (0.,0) -- +(1,0) ;
\draw[->,very thick,draw=green!50!black] (0.,0) -- +(1,0) ;
\draw[very thick,draw=white,double=white,rotate=30] (0.,0) -- +(1,0) ;
\draw[->,very thick,draw=green!50!black,rotate=30] (0.,0) -- +(1,0) ;
\draw[very thick,draw=white,double=white,rotate=60] (0.3,0.2) -- +(1.5,0) ;
\draw[->,very thick,draw=green!50!black,rotate=60] (0.3,0.2) -- +(1.5,0) ;
\end{tikzpicture}
\caption[Equivalent orientations.]{Equivalent external orientations of 1 and 2
  dimensional submanifold in 3-dim.  
 The actual length of the arc or arrow leaving is
  irrelevant. What matters is simply what the circulation around
  1-dim submanifold and  side of the 2-dim submanifold
  it points.}
\label{fig_equiv_ori}
\end{figure}
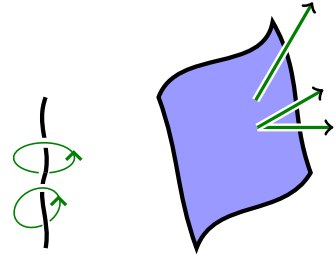

We emphasise that when specifying an orientation, we can point the
respective arrow or arcs in a variety of directions without changing the
orientation. For example, in figure \ref{fig_equiv_ori} we see how the
same orientation can be represented by several arcs or arrows.


\subsection{Internal and external, twisted and untwisted}
\label{ch_Orientation_Int_Ext}

\begin{table}[tb]
\footnotesize
\centering
\begin{tabular}{|l||c|l|l|}
\hline
& \multi{Always\\definable \& \\Essential} & Internal & External
\\\hline\hline
Manifold & No & Untwisted & Undefined
\\\hline
Submanifold & No & Untwisted & Twisted
\\\hline
Vectors & Yes & Untwisted & Twisted
\\\hline
Forms & Yes & Twisted & Untwisted
\\\hline
\end{tabular}
\caption[Types of orientation for all objects.]
{Objects which can have orientations and the
  correspondence between 
  untwisted versus twisted and internal versus external.}
\label{tab_ori_Int_Ext}
\end{table}

For manifolds and submanifold it is useful to label internal
orientations as untwisted and external orientations as twisted.
As we have stated, as well as manifolds and submanifolds, vectors and
forms also have orientations. These may also be twisted or
untwisted and may also be internal or external. The pictures of the
orientations depend on whether they are internal or external, but the
assignment to which are twisted or untwisted depend on the object,
according to table \ref{tab_ori_Int_Ext}.


\begin{table}[tb]
\centering
\def\dotplus{
      \tikz{\fill[blue] (0,0) circle (0.1) node[yshift=0.3cm] {+} ;}}
\def\followedby{\raisebox{1em}[3em][0em]{\footnotesize\multi{followed\\by}}}
\def\becomes{\raisebox{1em}[0em][0em]{\footnotesize becomes}}
\def\whichis{\raisebox{1em}[0em][0em]{\footnotesize\multi{which\\is}}}
\def\vertarrow{
      \tikz{ \draw [white] (-.1,0) -- (.1,0) ; 
         \draw[blue,very thick,->] (0,0) -- +(0,.5) ;
            }}
\def\horizarrow{
      \tikz{\draw [white] (0,-.1) -- (0,.1) ;
    \draw[blue,very thick,->] (0,0) -- +(0.5,0) ;}}
\def\looparrow{
\begin{tikzpicture}
  \draw[black!50,thick] (0.1,-.7) -- +(0,1.3) ;
  \draw[ultra thick,draw=white,double=blue,yscale=0.5] (0,0.5) arc (-260:80:0.5)  ;
  \draw[blue,very thick,->,yscale=0.5] (0,0.5) arc (-260:-90:0.5)  ;
  \draw[very thick,draw=blue,yscale=0.5] (0,0.5) arc (-260:80:0.5)  ;
  \end{tikzpicture}
}

\centering
\begin{tabular}{|@{}c@{}c@{}c@{}c@{}c@{}c@{}c@{\,}|}
\hline
\dotplus
&
\followedby
&
\vertarrow 
&
\becomes
&
\tikz{\draw[blue,very thick,->] (4.5,-.3) -- +(0,.5) ;
      \fill[blue] (4.5,-.3) circle (0.1) ;}
&
\whichis
&
\vertarrow
\\\hline
\vertarrow
&
\followedby
& 
\dotplus
&
\becomes
&
\tikz{\draw[blue,very thick,->] (4.5,-.3) -- +(0,.5) ; 
\fill[blue] (4.5,.3) circle (0.1) ;}
&
\whichis
&
\vertarrow
\\\hline
\horizarrow
&
\followedby
&
\vertarrow
&
\becomes
&
\begin{tikzpicture}
\draw[white] (-.1,-.1) -- (.85,.51) ;
\draw[blue,very thick,->] (0,0) -- +(0.7,0) ;
\draw[blue,very thick,->] (0.75,0) -- +(0,.5) ; 
\end{tikzpicture}
&
\whichis
&
\tikz{\draw[blue,very thick,<-] (6.5,-.3) arc (200:-20:0.5) ;} 
\\\hline
\vertarrow
&
\followedby
&
\horizarrow
&
\becomes
&
\begin{tikzpicture}
\draw[white] (-.1,0) -- (.8,.5) ;
\draw[blue,very thick,->] (0.03,.24) -- +(0.7,0) ;
\draw[blue,very thick,->] (0.,-.3) -- +(0,.5) ; 
\end{tikzpicture}
&
\whichis
&
\tikz{\draw[blue,very thick,->] (6.5,-.3) arc (200:-20:0.5) ;}
\\\hline
\looparrow
&
\followedby
&
\vertarrow
&
\becomes
&
\begin{tikzpicture}
\draw[black!50,thick] (4.6,-.7) -- +(0,1.3) ;
\draw[ultra thick,draw=white,double=blue,yscale=0.5] (4.5,0.) arc (-260:60:0.5)  ;
\draw[blue,very thick,->,yscale=0.5] (4.5,0.) arc (-260:-90:0.5)  ;
\draw[very thick,draw=blue,yscale=0.5] (4.5,0.) arc (-260:60:0.5)  ;
\draw[blue,very thick,->] (4.8,0) -- +(0,.4) ; 
\end{tikzpicture}
&
\whichis
&
\begin{tikzpicture}
\draw[ultra thick,draw=white,double=blue,xscale=.7,yscale=0.7,yshift=-2cm,xshift=9cm]
(0.5,1) \foreach \x in {0,10,...,720}
{ to ( {0.5*cos(\x)}, {0.0025*\x+0.5*sin(\x)+1} ) } ;
\draw[ultra thick,draw=white,double=blue,xscale=.7,yscale=0.7,yshift=-2cm,xshift=9cm]
(-.5,0.0025*180+1) \foreach \x in {180,190,...,360}
{ to ( {0.5*cos(\x)}, {0.0025*\x+0.5*sin(\x)+1} ) } ;
\draw[ultra thick,draw=white,double=blue,xscale=.7,yscale=0.7,yshift=-2cm,xshift=9cm]
(-.5,0.0025*540+1) \foreach \x in {540,550,...,720}
{ to ( {0.5*cos(\x)}, {0.0025*\x+0.5*sin(\x)+1} ) } ;
\end{tikzpicture}
\\\hline
\vertarrow
&
\followedby
&
\looparrow
&
\becomes
&
\begin{tikzpicture}
\draw[black!50,thick] (5,-.7) -- +(0,1.3) ;
\draw[ultra thick,draw=white,double=blue,yscale=0.5] (4.5,0.3) arc (-170:60:0.5)
  ;
\draw[blue,very thick,->,yscale=0.5] (4.5,0.3) arc (-170:-90:0.5)  ;
\draw[very thick,draw=blue,yscale=0.5] (4.5,0.3) arc (-170:80:0.5)  ;
\draw[blue,very thick,->-=.5] (4.5,-.45) -- +(0,.6) ; 
\end{tikzpicture}
&
\whichis
&
\begin{tikzpicture}
\draw[ultra thick,draw=white,double=blue,xscale=.7,yscale=0.7]
(0.5,1) \foreach \x in {0,10,...,720}
{ to ( {0.5*cos(\x)}, {0.0025*\x+0.5*sin(\x)+1} ) } ;
\draw[ultra thick,draw=white,double=blue,xscale=.7,yscale=0.7]
(-.5,0.0025*180+1) \foreach \x in {180,190,...,360}
{ to ( {0.5*cos(\x)}, {0.0025*\x+0.5*sin(\x)+1} ) } ;
\draw[ultra thick,draw=white,double=blue,xscale=.7,yscale=0.7]
(-.5,0.0025*540+1) \foreach \x in {540,550,...,720}
{ to ( {0.5*cos(\x)}, {0.0025*\x+0.5*sin(\x)+1} ) } ;
\end{tikzpicture}
\\\hline
\end{tabular}
\caption{Concatenation of orientations.}
\label{tab_concat_ori}
\end{table}


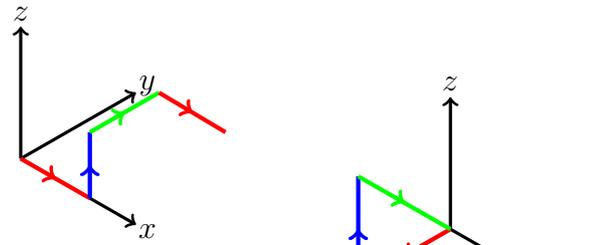
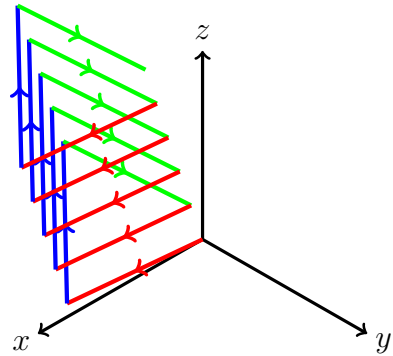
\begin{figure}[tb]
\begin{subfigure}{0.22\textwidth}
\pgfmathsetmacro{\len}{2.5}
\pgfmathsetmacro{\leny}{1}
\centering
\begin{tikzpicture}[scale=0.7]
\draw [very thick,->] (0,0) -- (0,\len) ;
\draw (0,1.1*\len) node {$z$} ;
\draw [very thick,->] (0,0) -- ({\len*cos(30)},{-\len*sin(30)}) ;
\draw ({1.1*\len*cos(30)},{-1.1*\len*sin(30)}) node {$x$} ;
\draw [very thick,->] (0,0) -- ({\len*cos(30)},{\len*sin(30)}) ;
\draw ({1.1*\len*cos(30)},{1.1*\len*sin(30)}) node {$y$} ;
\pgfmathsetmacro{\len}{1.5}
\pgfmathsetmacro{\dx}{0.0}
\pgfmathsetmacro{\dy}{0.0}

\draw [ultra thick,red,->-=.5] (0,0) -- ({\len*cos(30)},{-\len*sin(30)});
\draw [ultra thick,blue,->-=.5] 
({\len*cos(30)},{-\len*sin(30)}) -- ({\len*cos(30)},{\leny*sin(30)}) ;
\draw [ultra thick,green,->-=.5] 
({\len*cos(30)},{\leny*sin(30)}) --({2*\len*cos(30)},{(\len+\leny)*sin(30)}) ;
\draw [ultra thick,red,->-=.5] 
({2*\len*cos(30)},{(\len+\leny)*sin(30)}) -- 
({(\len+\leny)*sin(30)+2*\len*cos(30)},{\leny*sin(30)});
\end{tikzpicture}
\caption{
A curve created by going a short distance in $x$ (red), 
followed by a short distance in $z$ (blue) followed by going a short
distance in $y$ (green) and then repeating. 
}
\label{fig_ori_pipclearner_a}
\end{subfigure}
\begin{subfigure}{0.22\textwidth}
\pgfmathsetmacro{\len}{2.5}
\centering
\begin{tikzpicture}[scale=0.7]
\draw [very thick,->] (0,0) -- (0,\len) ;
\draw (0,1.1*\len) node {$z$} ;
\draw [very thick,->] (0,0) -- ({\len*cos(30)},{-\len*sin(30)}) ;
\draw ({1.1*\len*cos(30)},{-1.1*\len*sin(30)}) node {$y$} ;
\draw [very thick,->] (0,0) -- (-{\len*cos(30)},{-\len*sin(30)}) ;
\draw ({-1.1*\len*cos(30)},{-1.1*\len*sin(30)}) node {$x$} ;
\pgfmathsetmacro{\len}{2}
\pgfmathsetmacro{\dx}{0.0}
\pgfmathsetmacro{\dy}{0.0}
\foreach \t in {0} {
\draw [ultra thick,red,->-=.5] ({3*\t*\dx},{3*\t*\dy}) 
-- ({\len*cos(150)+(3*\t+1)*\dx},{-\len*sin(150)+(3*\t+1)*\dy});
\draw [ultra thick,blue,->-=.5] ({\len*cos(150)+(3*\t+1)*\dx},{-\len*sin(150)+(3*\t+1)*\dy}) 
-- ({\len*cos(150)+(3*\t+2)*\dx},{\len*sin(150)+(3*\t+2)*\dy})  ;
\draw [ultra thick,green,->-=.5]  ({\len*cos(150)+(3*\t+2)*\dx},{\len*sin(150)+(3*\t+2)*\dy}) 
-- ({(3*\t+3)*\dx},{(3*\t+3)*\dy}) ;
}
\end{tikzpicture}
\caption{The same diagram as in figure
  \ref{fig_ori_pipclearner_a} but now
  with the axes all pointing towards the reader.}
\label{fig_ori_pipclearner_b}
\end{subfigure}
\begin{subfigure}{0.4\textwidth}
\centering
\pgfmathsetmacro{\len}{2.5}
\begin{tikzpicture}
\draw [very thick,->] (0,0) -- (0,\len) ;
\draw (0,1.1*\len) node {$z$} ;
\draw [very thick,->] (0,0) -- ({\len*cos(30)},{-\len*sin(30)}) ;
\draw ({1.1*\len*cos(30)},{-1.1*\len*sin(30)}) node {$y$} ;
\draw [very thick,->] (0,0) -- (-{\len*cos(30)},{-\len*sin(30)}) ;
\draw ({-1.1*\len*cos(30)},{-1.1*\len*sin(30)}) node {$x$} ;
\pgfmathsetmacro{\len}{2}
\pgfmathsetmacro{\dx}{-0.05}
\pgfmathsetmacro{\dy}{0.15}
\foreach \t in {0,1,2,3,4} {
\draw [ultra thick,red,->-=.5] ({3*\t*\dx},{3*\t*\dy}) 
-- ({\len*cos(150)+(3*\t+1)*\dx},{-\len*sin(150)+(3*\t+1)*\dy});
\draw [ultra thick,blue,->-=.5] ({\len*cos(150)+(3*\t+1)*\dx},{-\len*sin(150)+(3*\t+1)*\dy}) 
-- ({\len*cos(150)+(3*\t+2)*\dx},{\len*sin(150)+(3*\t+2)*\dy})  ;
\draw [ultra thick,green,->-=.5]  ({\len*cos(150)+(3*\t+2)*\dx},{\len*sin(150)+(3*\t+2)*\dy}) 
-- ({(3*\t+3)*\dx},{(3*\t+3)*\dy}) ;
}
\end{tikzpicture}
\caption{The same diagram as in figure
  \ref{fig_ori_pipclearner_b} but with
  a offset each loop to show the helix.}
\label{fig_ori_pipclearner_c}
\end{subfigure}

\caption[Creating a helix.]{Demonstration that going a short distance
  in $x$ then $z$ then $y$ and repeating creates a helix. The reader
  is encouraged to experiment with
  some pipe cleaners.}
\label{fig_ori_pipclearner}
\end{figure}


\subsection{Concatenation of orientations} 
\label{ch_Orientation_Concat}

We can combine two orientations
together. Simply take the two orientations and join them together, the second
onto the end of the first. See table \ref{tab_concat_ori}. 

You can see that the untwisted orientations for two dimensions,
figure \ref{fig_ori_Man2}, i.e. clockwise or anticlockwise arise from
the concatenation of two orientations, as in rows three and four of
table \ref{tab_concat_ori}. Likewise the helices, figure
\ref{fig_ori_Man3} arise by
concatenating an arc with a perpendicular line as in rows five and six of
table \ref{tab_concat_ori}.

Using pipe cleaners may help you to see that that an arc followed by a
perpendicular line is equivalent to a helix. Likewise a helix arises
in 3 dimensions by tracing out a path, going in one direction, then a
second, then the third and repeating as in figure
\ref{fig_ori_pipclearner}.


\subsection{Inheritance of orientation by the boundary of a submanifold} 
\label{ch_Orientation_boundary}

\begin{figure}[tb]
\centering
\begin{subfigure}[t]{.25\figsize}
\centering
\begin{tikzpicture}
\draw[->,thick] (0,0) to[out=80,in=-110] (0,1) ; 
\draw[thick] (0,1) to[out=80,in=-110] (0,2) ; 
\fill[red] (0,2) circle (0.1) node[xshift=0.2cm,yshift=0.3] {+} ;
\fill[blue] (0,0) circle (0.1) node[xshift=0.2cm,yshift=-0.3] {-} ;
\end{tikzpicture}
\caption{Boundary of a 1-dim submanifold}
\label{fig_Ori_Bdd-Int13}
\end{subfigure}
\
\begin{subfigure}[t]{.25\figsize}
\centering
\begin{tikzpicture}
\filldraw[ultra thick,draw=green!60!black,fill=green,xscale=0.6] 
(0,0) circle (1) ;
\draw[thick,->] (0.25,0) arc(0:130:0.25) ; 
\draw[thick] (0.25,0) arc(0:260:0.25) ; 
\draw[ultra thick,xscale=0.6] (0.8,0.) -- +(.2,.2) -- +(.4,0) ;
\end{tikzpicture}
\caption{Boundary of a 2-dim submanifold}
\label{fig_Ori_Bdd-Int23}
\end{subfigure}
\
\begin{subfigure}[t]{.45\figsize}
\centering
\begin{tikzpicture}[scale=0.8]
\shadedraw[ball color=green!80] (1.8,1.6) circle (2) ;

\draw[ultra thick,draw=white,double=black,xshift=2cm,yshift=-0.7cm]
(0.5,1) \foreach \x in {0,10,...,720}
{ to ( {0.5*cos(\x)}, {0.0025*\x+0.5*sin(\x)+1} ) } ;
\draw[ultra thick,draw=white,double=black,xshift=2cm,yshift=-0.7cm]
(-.5,0.0025*180+1) \foreach \x in {180,190,...,360}
{ to ( {0.5*cos(\x)}, {0.0025*\x+0.5*sin(\x)+1} ) } ;
\draw[ultra thick,draw=white,double=black,xshift=2cm,yshift=-0.7cm]
(-.5,0.0025*540+1) \foreach \x in {540,550,...,720}
{ to ( {0.5*cos(\x)}, {0.0025*\x+0.5*sin(\x)+1} ) } ;

\draw[ultra thick,-<,blue,xscale=1.3] (1.7,2.8) arc(0:130:0.3) ; 
\draw[ultra thick,blue,xscale=1.3] (1.7,2.8) arc(0:260:0.3) ; 

\end{tikzpicture}
\caption{Boundary of a 3-dim submanifold. The helix is the orientation
of the 3-dim ball.}
\label{fig_Ori_Bdd-Int33}
\end{subfigure}

\caption[Inheritance of internal orientation on a boundary.]
{Inheritance of an internal orientation of a submanifold by its boundary for 
1,2,3-dimensional submanifold embedded into a 3-dimensional manifold.}
\label{fig_Ori_Bdd-Int}
\end{figure}
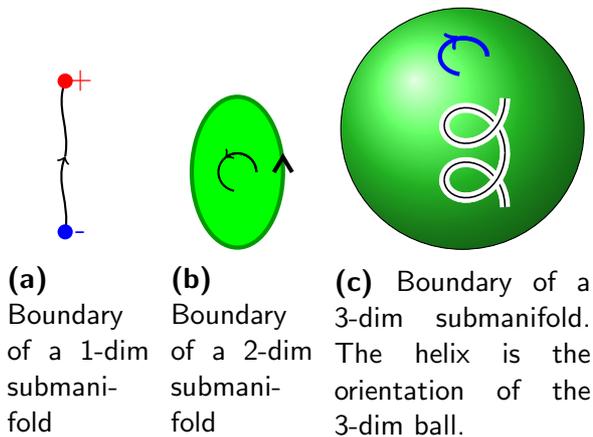


\begin{figure}[tb]
\centering
\begin{subfigure}[t]{.25\figsize}
\centering
\begin{tikzpicture}
\draw[thick] (0,0.3) to[out=80,in=-110] (0,1) ; 
\draw[thick] (0,1) to[out=80,in=-110] (0,1.7) ; 
\draw[thick,yscale=-0.5,draw=white,double=green!50!black] 
(-0.09,-2.8) arc (260:-80:0.4) ; 
\draw[thick,draw=green!50!black] 
(0.27,1.17) -- +(0.1,0.1) -- +(0.2,0) ; 
\draw[ultra thick,draw=white,double=red,xscale=-.4,yscale=0.4,yshift=3cm]
(0.5,1) \foreach \x in {0,10,...,720}
{ to ( {0.5*cos(\x)}, {0.0025*\x+0.5*sin(\x)+1} ) } ;
\draw[ultra thick,draw=white,double=red,xscale=-.4,yscale=0.4,yshift=3cm]
(-.5,0.0025*180+1) \foreach \x in {180,190,...,360}
{ to ( {0.5*cos(\x)}, {0.0025*\x+0.5*sin(\x)+1} ) } ;
\draw[ultra thick,draw=white,double=red,xscale=-.4,yscale=0.4,yshift=3cm]
(-.5,0.0025*540+1) \foreach \x in {540,550,...,720}
{ to ( {0.5*cos(\x)}, {0.0025*\x+0.5*sin(\x)+1} ) } ;
\draw[ultra thick,draw=white,double=blue,xscale=.4,yscale=0.4,yshift=-2cm]
(0.5,1) \foreach \x in {0,10,...,720}
{ to ( {0.5*cos(\x)}, {0.0025*\x+0.5*sin(\x)+1} ) } ;
\draw[ultra thick,draw=white,double=blue,xscale=.4,yscale=0.4,yshift=-2cm]
(-.5,0.0025*180+1) \foreach \x in {180,190,...,360}
{ to ( {0.5*cos(\x)}, {0.0025*\x+0.5*sin(\x)+1} ) } ;
\draw[ultra thick,draw=white,double=blue,xscale=.4,yscale=0.4,yshift=-2cm]
(-.5,0.0025*540+1) \foreach \x in {540,550,...,720}
{ to ( {0.5*cos(\x)}, {0.0025*\x+0.5*sin(\x)+1} ) } ;
\end{tikzpicture}
\caption{Boundary of a 1-dim submanifold}
\label{fig_Ori_Bdd-Ext13}
\end{subfigure}
\ 
\begin{subfigure}[t]{.25\figsize}
\centering
\begin{tikzpicture}
\filldraw[ultra thick,draw=green!60!black,fill=green,xscale=0.6] 
(0,0) circle (1) ;
\draw[ultra thick,draw=white,double=white] (0.,0.) -- +(1,0) ;
\draw[ultra thick,->] (0.,0.) -- +(1,0) ;
\draw[ultra thick,draw=white,double=black] (0.1,0.7) arc(-70:180:0.25) ; 
\draw[thick,->] (0.1,0.7) arc(-70:130:0.25) ; 
\end{tikzpicture}
\caption{Boundary of a 2-dim submanifold}
\label{fig_Ori_Bdd-Ext23}
\end{subfigure}
\ 
\begin{subfigure}[t]{.45\figsize}
\centering
\begin{tikzpicture}
\shadedraw[ball color=green!80] (0,0) circle (1.5) ;
\draw(0,0) node {\Large +} ;
\draw[ultra thick,draw=white,double=white] (0.8,0.5) -- +(0.8,0.5) ;
\draw[ultra thick,-<] (0.8,0.5) -- +(0.8,0.5) ;
\draw[ultra thick] (0.8,0.5) -- +(1.,0.5*5/4) ;
\end{tikzpicture}
\caption{Boundary of a 3-dim submanifold, The $+$ is the external orientation of the 3-dim ball and the arrow
  pointing inwards is the external orientation of the 2-dim sphere.}
\label{fig_Ori_Bdd-Ext33}
\end{subfigure}
\caption[Inheritance of external orientation on a boundary.]
        {Inheritance of an external orientation of a submanifold by
          its boundary for 1,2,3-dimensional submanifold embedded into
          a 3-dimensional manifold.}
\label{fig_Ori_Bdd-Ext}
\end{figure}
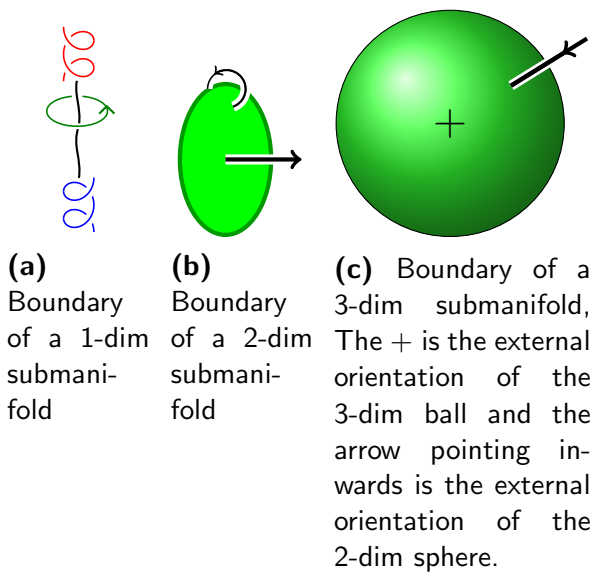

If a submanifold has an internal orientation and a boundary then the boundary
inherits an internal orientation, see figure \ref{fig_Ori_Bdd-Int}. 
Likewise if a submanifold has an
external orientation then the boundary inherits the external
orientation, figure \ref{fig_Ori_Bdd-Ext}.


\subsection{Orientations of vectors} 
\label{ch_Orientation_Vectors}

In section \ref{ch_DG_Vectors} and figure \ref{fig_Vector_field} we
saw the usual type of vectors. These are untwisted vectors, since they
have an untwisted orientations which is the internal
orientations. There exists another type of vector called a twisted
vector. These are also called axial or pseudo-vectors. An example of
such a vector is the magnetic field $\VB$. Twisted vectors have
external orientations. They look like the external orientations of
1-dimensional submanifolds as depicted in table
\ref{tab_ori_twist_SubMan}. In the case of vectors in 3-dimension one
can see why they are call axial.

\subsection{Twisting and untwisting}
\label{ch_Orientation_Twisting}

\begin{table}[tb]
\footnotesize
\centering
\begin{tabular}{|@{\,}c@{\,}|@{\,}c@{\,}|@{\,}c@{\,}|@{\,}c@{\,}|}
\hline
\raisebox{-.5em}{\multi{Embedding\\dimension}} 
& \raisebox{-.5em}{\multi{Embedding\\orientation}}
& \multi{Submanifold\\external\\orientation}
& \multi{Submanifold\\internal\\orientation}
\\\hline
\rule{0em}{3.4em}
\raisebox{1em}{2} &
\begin{tikzpicture}[scale=1.3]
\draw [thick,black] (0,0) rectangle (.7,.7) ; 
\draw [very thick,red,->-=.5] (0.35,.35) +(70:.2) arc (70:-250:.2) ; 
\end{tikzpicture}
&
\begin{tikzpicture}[scale=1.3]
\draw [thick,black] (0,0) rectangle (.7,.7) ; 
\draw [ultra thick, green] (.35,0) -- +(0,.7) ; 
\draw [very thick,red,->] (0.35,.35) -- +(-.3,0) ; 
\end{tikzpicture}
&
\begin{tikzpicture}[scale=1.3]
\draw [thick,black] (0,0) rectangle (.7,.7) ; 
\draw [ultra thick, green] (.35,0) -- +(0,.7) ; 
\draw [very thick,red,->] (0.35,.35) -- +(0,.3) ; 
\end{tikzpicture}
\\\hline
\rule{0em}{4em}
\raisebox{1em}{3} &
\begin{tikzpicture}[scale=1.3]
\draw[ultra thick,draw=white,double=red,xscale=-.4,yscale=0.4,yshift=3cm]
(0.5,1) \foreach \x in {0,10,...,720}
{ to ( {0.5*cos(\x)}, {0.0025*\x+0.5*sin(\x)+1} ) } ;
\draw[ultra thick,draw=white,double=red,xscale=-.4,yscale=0.4,yshift=3cm]
(-.5,0.0025*180+1) \foreach \x in {180,190,...,360}
{ to ( {0.5*cos(\x)}, {0.0025*\x+0.5*sin(\x)+1} ) } ;
\draw[ultra thick,draw=white,double=red,xscale=-.4,yscale=0.4,yshift=3cm]
(-.5,0.0025*540+1) \foreach \x in {540,550,...,720}
{ to ( {0.5*cos(\x)}, {0.0025*\x+0.5*sin(\x)+1} ) } ;
\end{tikzpicture}
&
\begin{tikzpicture}[scale=1.3]
\draw [ultra thick, green] (.25,0) -- +(0,1) ; 
\draw [very thick,white,yscale=.6,double=white] (0.25,1) +(70:.2) arc (70:-250:.2) ; 
\draw [very thick,red,yscale=.6,->-=.5] (0.25,1) +(70:.2) arc (70:-250:.2) ; 
\end{tikzpicture}
&
\begin{tikzpicture}[scale=1.3]
\draw (-.1,0) (.1,0) ;
\draw [ultra thick, green] (0,0) -- +(0,.7) ; 
\draw [very thick,red,->] (0,.35) -- +(0,.3) ; 
\end{tikzpicture}
\\\hline
\end{tabular}
\caption{Examples of twisting and untwisting.}
\label{tab_Ori_Twisting}
\end{table}

If a manifold has an orientation then we can use it to change the type
of orientation of a submanifold. Thus we can convert an external
orientation into an internal orientation. We can also convert an
internal orientation into an external orientation. The convention we
use here is to say that if we concatenate the external orientation
(first) followed by the internal orientation (second) one must arrive
at the orientation of the manifold. See examples in table
\ref{tab_Ori_Twisting}.

For example if the orientation of a 2-dim manifold is clockwise and the
untwisted orientation is up then the twisted orientation is right. See
fourth line of table \ref{tab_concat_ori}.

Starting with an
untwisted form, then twisting and then untwisting it does not change
the orientation.



\begin{figure}[tb]
\begin{subfigure}[t]{0.48\figsize}
\begin{tikzpicture}[scale=0.7]
\shadedraw [inner color=white, outer color = lightgray,
scale=1.05,shift={(-.2,-.2)}] 
(0,0) to [out=0,in=-135] (4,0) to [in=-135] (4,4) to [in=45]
(0,4) to [out=-135,in=180] cycle ;
\draw[color=red] (1,1.3) node  {\footnotesize $\boldsymbol\oplus$}
  \foreach \x in {0,60,...,300} {
 +(\x:0.5) node {\footnotesize $\boldsymbol\oplus$}
 +(\x:1.3) node {\footnotesize $\boldsymbol\oplus$}
  }   
\foreach \x in {-30,0,30} {
 +(\x:2.5) node {\footnotesize $\boldsymbol\oplus$}
  } 
\foreach \x in {120,150} {
 +(\x:1.9) node {\footnotesize $\boldsymbol\oplus$}
  } 
;
\fill[color=blue] (2,3.5)   node {\footnotesize $\boldsymbol\ominus$}
  \foreach \x in {0,60,...,300} {
 +(\x:0.5)  node {\footnotesize $\boldsymbol\ominus$}
  } 
 +(150:1.3) node {\footnotesize $\boldsymbol\ominus$}
 +(180:1.6) node {\footnotesize $\boldsymbol\ominus$}
;

\end{tikzpicture}
\caption{Low density.}
\label{fig_low-density}
\end{subfigure}
\quad
\begin{subfigure}[t]{0.48\figsize}
\centering
\begin{tikzpicture}[scale=0.7]
\shadedraw [inner color=white, outer color = lightgray,
scale=1.05,shift={(-.2,-.2)}] 
(0,0) to [out=0,in=-135] (4,0) to [in=-135] (4,4) to [in=45]
(0,4) to [out=-135,in=180] cycle ;
\fill[color=red] (1,1.3) 
  node {\tiny $\boldsymbol\oplus$}
  \foreach \x in {0,30,...,330} {
 +(\x:0.2) node {\tiny $\boldsymbol\oplus$}
 +(\x:0.5) node {\tiny $\boldsymbol\oplus$}
 +(\x:0.9) node {\tiny $\boldsymbol\oplus$}
 +(\x:1.3) node {\tiny $\boldsymbol\oplus$} 
  }   
\foreach \x in {-30,0,30} {
 +(\x:1.9) node {\tiny $\boldsymbol\oplus$} 
 +(\x:2.5) node {\tiny $\boldsymbol\oplus$} 
  } 
 +(130:2.4) node {\tiny $\boldsymbol\oplus$}
 +(110:2.4) node {\tiny $\boldsymbol\oplus$}
 +(110:1.8) node {\tiny $\boldsymbol\oplus$}
\foreach \x in {120,150,180,210,-60} {
 +(\x:1.85) node {\tiny $\boldsymbol\oplus$} 
  } 
 +(-15:2.4) node {\tiny $\boldsymbol\oplus$}
 +(-15:2.8) node {\tiny $\boldsymbol\oplus$}
;
\fill[color=blue] (2,3.5) 
  node {\tiny $\boldsymbol\ominus$} 
  \foreach \x in {0,30,...,330} {
 +(\x:0.2) node {\tiny $\boldsymbol\ominus$} 
 +(\x:0.5) node {\tiny $\boldsymbol\ominus$} 
 +(\x:0.9) node {\tiny $\boldsymbol\ominus$} 
  } 
 +(150:1.3) node {\tiny $\boldsymbol\ominus$} 
 +(180:1.6) node {\tiny $\boldsymbol\ominus$} 
;
\end{tikzpicture}
\caption{Increased density, while reducing the value of each dot.}
\label{fig_high-density}
\end{subfigure}
\caption[Representation of a top-form.]  {Representation of a
  twisted top-form, i.e. an $n$-form on an $n$-dimensional manifold. }
\label{fig_n-form_field}
\end{figure}
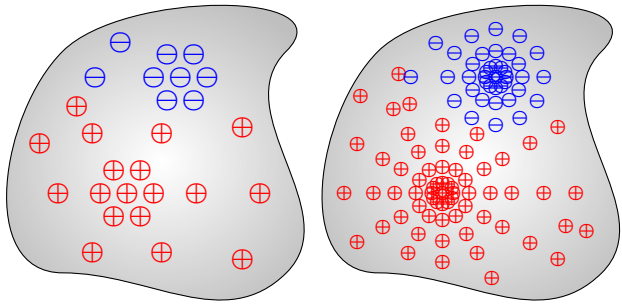

\begin{figure}[tb]
\begin{tikzpicture}[scale=1.3]
\shadedraw [inner color=white, outer color = lightgray] 
(0,0) to [out=0,in=-135] (4,0) to [in=-135] (4,4) to [in=45]
(0,4) to [out=-135,in=180] cycle ;
\draw[ultra thick,blue] (0,0.5) to [out=110,in=-70] +(0,3) ;

\draw[very thick,double=black,draw=black!10!white]  (0,1.8) arc (70:-250:.2 and 0.05) ;
\draw[very thick,->-=.7]  (0,1.8) arc (70:-250:.2 and 0.05) ;

\draw[ultra thick,blue] (0.3,0.8) to [out=110,in=-70] +(0,3) ;
\draw[very thick,double=black,draw=black!10!white] (.25,1.2)   arc (70:-250:.2 and 0.05) ;
\draw[very thick,->-=.7]  (.25,1.2) arc (70:-250:.2 and 0.05) ;

\draw[ultra thick,blue] (1,1.) to [out=110,in=-80] +(0,3) ;
\draw[very thick,double=black,draw=black!5!white] (.95,1.3)   arc (70:-250:.2 and 0.05) ;
\draw[very thick,->-=.7]  (.95,1.3) arc (70:-250:.2 and 0.05) ;

\draw[ultra thick,blue,rotate=-10] (1.1,2.3) ellipse (0.2 and 0.5) ;
\draw[very thick,double=black,draw=black!5!white]  (1.4,2.3)  arc (70:-250:.2 and 0.05) ;
\draw[very thick,->-=.7]  (1.4,2.3) arc (70:-250:.2 and 0.05) ;

\draw[ultra thick,blue] (2,0.5) to [out=110,in=-70] +(0,3) ;
\draw[very thick,double=black,draw=black!5!white]  (2.12,2.3)  arc (70:-250:.2 and 0.05) ;
\draw[very thick,-<-=.4]  (2.12,2.3) arc (70:-250:.2 and 0.05) ;

\draw[ultra thick,blue] (2.8,0.8) to [out=130,in=-80] +(0,3) ;
\draw[very thick,double=black,draw=black!5!white]  (2.65,2.3)  arc (70:-250:.2 and 0.05) ;
\draw[very thick,-<-=.4]  (2.65,2.3) arc (70:-250:.2 and 0.05) ;

\draw[ultra thick,blue] (3,1.0) to [out=120,in=-70] +(0,3) ;
\draw[very thick,double=black,draw=black!5!white]  (2.9,2.1)  arc (70:-250:.2 and 0.05) ;
\draw[very thick,-<-=.4]  (2.9,2.1) arc (70:-250:.2 and 0.05) ;

\draw[ultra thick,blue] (3.3,1.2) to [out=120,in=-70] +(0,3) ;
\draw[very thick,double=black,draw=black!5!white]  (3.2,1.9)  arc (70:-250:.2 and 0.05) ;
\draw[very thick,-<-=.4]  (3.2,1.9) arc (70:-250:.2 and 0.05) ;

\end{tikzpicture}
\caption[An untwisted 2-form on a 3-dim manifold.]{An untwisted 2-form
  on a 3-dim manifold.}
\label{fig_2-form_untw_field}
\end{figure}

\begin{figure}[tb]
\begin{tikzpicture}[scale=1.3]
\shadedraw [inner color=white, outer color = lightgray] 
(0,0) to [out=0,in=-135] (4,0) to [in=-135] (4,4) to [in=45]
(0,4) to [out=-135,in=180] cycle ;
\draw[ultra thick,blue,->-=.7] (0,0.5) to [out=110,in=-70] +(0,3) ;
\draw[ultra thick,blue,->-=.7] (0.3,0.8) to [out=110,in=-70] +(0,3) ;
\draw[ultra thick,blue,->-=.7] (1,1.) to [out=110,in=-80] +(0,3) ;
\draw[ultra thick,blue,rotate=-10] (1.1,2.3) ellipse (0.2 and 0.5) ;
\draw[ultra thick,blue,rotate=-10,<-] (.88,2.4)   +(0.005,.01) ;
\draw[ultra thick,blue,-<-=.7] (2,0.5) to [out=110,in=-70] +(0,3) ;
\draw[ultra thick,blue,-<-=.7] (2.8,0.8) to [out=130,in=-80] +(0,3) ;
\draw[ultra thick,blue,-<-=.7] (3,1.0) to [out=120,in=-70] +(0,3) ;
\draw[ultra thick,blue,-<-=.7] (3.3,1.2) to [out=120,in=-70] +(0,3) ;
\end{tikzpicture}
\caption[An untwisted 2-form on a 3-dim manifold.]{An untwisted 2-form on a 3-dim manifold.}
\label{fig_2-form_twist_field}
\end{figure}


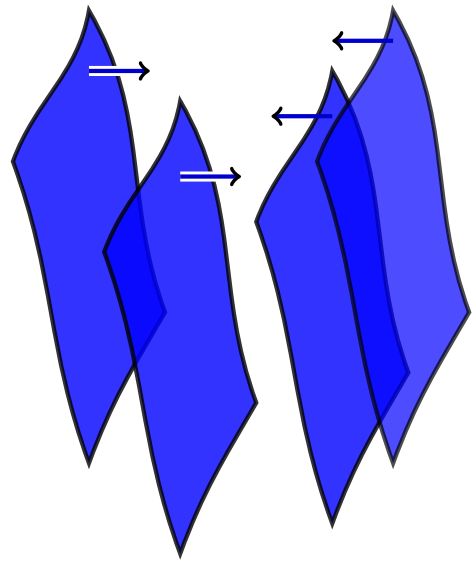
\begin{figure}[tb]
\begin{tikzpicture}[scale=2]
\filldraw[ultra thick,draw=black,fill=blue,opacity=0.8] (0.4,0.9) 
to [out=110,in=-70] +(-.5,2) 
to [out=70,in=-100] +(0.5,1) 
to [out=-60,in=110] +(0.5,-2) 
to [out=-120,in=70] cycle
;
\filldraw[ultra thick,draw=black,fill=blue,opacity=0.8] (1.0,0.3) 
to [out=110,in=-70] +(-.5,2) 
to [out=70,in=-100] +(0.5,1) 
to [out=-60,in=110] +(0.5,-2) 
to [out=-120,in=70] cycle
;
\filldraw[ultra thick,draw=black,fill=blue,opacity=0.8] (2,0.5) 
to [out=110,in=-70] +(-.5,2) 
to [out=70,in=-100] +(0.5,1) 
to [out=-60,in=110] +(0.5,-2) 
to [out=-120,in=70] cycle
;
\filldraw[ultra thick,draw=black,fill=blue,opacity=0.7] (2.4,0.9) 
to [out=110,in=-70] +(-.5,2) 
to [out=70,in=-100] +(0.5,1) 
to [out=-60,in=110] +(0.5,-2) 
to [out=-120,in=70] cycle
;
\draw[ultra thick,draw=white,double=white] (0.4,3.5) -- +(.4,0) ;
\draw[ultra thick,draw=blue!80!black,->] (0.4,3.5) -- +(.4,0) ;
\draw[ultra thick,draw=white,double=white] (1.,2.8) -- +(.4,0) ;
\draw[ultra thick,draw=blue!80!black,->] (1.,2.8) -- +(.4,0) ;

\draw[ultra thick,draw=blue!80!black,->] (2.0,3.2) -- +(-.4,0) ;
\draw[ultra thick,draw=blue!80!black,->] (2.4,3.7) -- +(-.4,0) ;
\end{tikzpicture}
\caption[An untwisted 1-form on a 3-dim manifold]{An untwisted 1-form
  on a 3-dim manifold. Note that the 1-form must vanish in the region
  between the two opposing orientations.  A diagram of a twisted
  1-form in 3 dimensions is given in figure
  \ref{fig_int_1-form_field}.}
\label{fig_1-form_field}
\end{figure}

\section{Closed forms.}
\label{ch_CForms}

In this section we introduce exterior differential forms, which we
call simply forms.  On an $n$-dimensional manifold a form has a
degree, which is an integer between 0 and $n$. We refer for a form
with degree $p$ as a $p$-form. A closed form of degree $p$ can be
thought of as a collection of submanifolds\footnote{Some authors
  \cite{schouten1954tensor,burke1985applied,hehl2012foundations} use a double sheet to indicate 1-forms and
  likewise they replace 1-dim form-submanifolds with cylinders to
  indicate $(n-1)$-forms. It is true that for a single $1$-form at a
  point it is necessary to have two $(n-1)$-dim
  form-elements-submanifold in order indicate the magnitude of the
  form. Consider in figure \ref{fig_vector_1form} one needs at least
  2-lines to indicate the combination of a 1-form with a vector where
  both are defined only at one point. However in this document we only
  consider smooth form-fields. I.e. a form defined at each point. Thus
  there is a form-submanifold passing through each point. Furthermore
  since the form-field is smooth form submanifolds for nearby points
  are nearly parallel. Therefore there is sufficient information to
  indicate the magnitude of the form via the density of the
  form-submanifolds. Looking at figure \ref{fig_vector_1form} again we
  see there is sufficient information to calculate the action of the
  vector on the form.} of dimension $n-p$.  We call these
\defn{form-submanifolds}. Similar to a vector field, there is a
submanifold of the correct dimension passing though each point in the
manifold.  They must not have boundaries, therefore they either
continue off to ``infinity'' or they close in on themselves, like
spheres. In section \ref{ch_NonClosed_Forms} we see that the
form-submanifolds with boundaries correspond to non-closed forms.

Although the points, lines, surfaces corresponding to $p$-forms, for
$p\ge1$, do not have values, they do have an orientation and as before
these can be untwisted or twisted. As stated in table
\ref{tab_ori_Int_Ext} these are the opposite way round than for
submanifolds. Thus an untwisted form has an external orientation and a
twisted form has an internal orientation. The untwisted forms are
mathematically simpler and much easier to take wedge products of. See
section \ref{ch_CForms_Wedge} below. Indeed the standard method of
prescribing the internal orientation of an $n$-dim manifold is in
terms of an untwisted $n$-form.

We have the following: 
\begin{jgitemize}
\item
An untwisted 0-form (internal orientation) is simply a scalar field
and is best visualised as above in figures \ref{fig_Scalar_field},
\ref{fig_Scalar_mountain} and \ref{fig_Scalar_contours}.  Note that a
closed 0-form is a constant and not very interesting. From the 
prescription above, a 0-form on an $n$-dim manifold is an $n$-dim
form-submanifold at each point. 

A twisted 0-form (external orientation) requires an external
orientation for a 0-dim submanifold (point) as in table
\ref{tab_ori_twist_SubMan}. 

\item
On an $n$-dimensional manifold, the form-submanifolds for an $n$-form
consists of a collection of dots.  
For the rest of this report, we shall refer to $n$-form as \defn{top-forms}. 

For twisted top forms (internal orientation) these dots
carry a plus or minus sign. See figure
\ref{fig_n-form_field}.  The density of the dots corresponds loosely to
the magnitude of the top-form. The smoothness of the field requires
that between a high density of dots and a low density dots, there
should be a region of intermediate density. 

A non-vanishing twisted top-form with only plus dots is called a
\defn{measure}. We see below that these are very useful for
integrating scalar fields (section \ref{ch_CForms_Intergration})
and for converting between vector fields and
$(n-1)$-forms (section \ref{ch_Internal_contraction}).
We see in section \ref{ch_Metric_measure} that a
metric automatically gives rise to a measure.

Although, as with vector fields, there is really a
``dot at each point'', therefore the density of dots is actually
infinite. However we can think of a limiting process. Start with a set
of discrete dots, then (approximately) double the density of dots,
while at the same time halving the value associated with each
points. As such the number of dots multiplied by their sign inside
a particular $n$-volume is approximately constant during the limiting
process. Compare this to the scalar field, where as we stated, the
value of each point remains constant during the limiting process. 

For untwisted top-forms 
each point carries an external orientation, as
in 0 dimensional submanifolds in table \ref{tab_ori_twist_SubMan}. 

\item
On an $n$-dimensional manifold the closed $(n-1)$-forms correspond to
$1$-dimensional form-submanifolds, i.e. curves as depicted in figures
\ref{fig_2-form_untw_field} and \ref{fig_2-form_twist_field}.  The
curves do not have start points or end points. Therefore they must
either close on themselves to form circles or they must go from and to
``infinity''. They also do not intersect.  Again there is really a
curve going though each point. The smoothness again requires that
between a high density of curves and a low density there should be an
intermediate density, and that between curves of opposite sign there
should be a gap.

Figure \ref{fig_2-form_untw_field} shows an untwisted 2-form (external
orientation) in 3 dimensions. 

Twisted $(n-1)$-forms are lines with direction. (Internal
orientation). One may consider these as flow lines of a
fluid or similar. An example is shown in figure
\ref{fig_2-form_twist_field}.

You can see from figures \ref{fig_Vector_field} and
\ref{fig_2-form_twist_field} that untwisted vectors and twisted
$(n-1)$-forms look similar, and indeed they are. With just a
measure we can convert a
vector into a 1-form via the internal contraction. See section
\ref{ch_Internal_contraction}. 

\item
On an $n$-dimensional manifold the closed 1-forms correspond to
$(n-1)$-dimensional form-submanifolds. 

Untwisted 1-forms
are the submanifold
corresponding to the contours or level surfaces of a scalar field, as
in figure \ref{fig_Scalar_contours}. The
external orientation point in the direction of increasing value of the
scalar field. The
faster the rate of increase of the scalar field, the closer the level
surfaces. 

One usually takes the gradient of a scalar field which is the
vector corresponding  to the direction of steepest assent. This
gradient is perpendicular to the level surfaces. However to
do this one needs a concept of orthogonality which comes with the
metric defined below in section \ref{ch_Metric}.

An example of an untwisted 1-form in 3-dimensions is given in figure
\ref{fig_1-form_field}.

\end{jgitemize}


\subsection{Multiplication of vectors and forms by scalar fields.}  
\label{ch_Mult_scalar}

To multiply a $p$-form field by a scalar field
simply increase the density of dots, curves, surfaces of the $p$-form
by the value of the scalar field. 

The action of multiplying a $p$-form by a negative scalar field
changes the orientation of the $p$-form.

It is worth noting that if you multiply a closed $p$-form with $p<n$
by a non constant scalar the result will not be a closed form.

\subsection{Addition of forms.}  
\label{ch_Add_Forms}

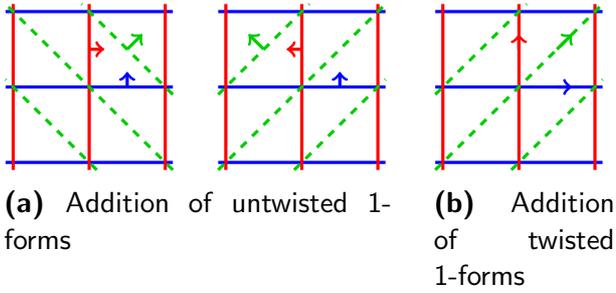
\begin{figure}[tb]
\begin{subfigure}[t]{0.3\textwidth}
\begin{tikzpicture}
\draw [very thick,blue] (-.1,0) -- (2.1,0) ;
\draw [very thick,blue] (-.1,1) -- (2.1,1) ;
\draw [very thick,blue] (-.1,2) -- (2.1,2) ;
\draw [very thick,blue,->] (1.5,1) -- +(0,.2) ;

\draw [very thick,red] (0,-.1) -- (0,2.1) ;
\draw [very thick,red] (1,-.1) -- (1,2.1) ;
\draw [very thick,red] (2,-.1) -- (2,2.1) ;
\draw [very thick,red,->] (1,1.5) -- +(.2,0) ;

\draw [very thick,green!80!black,dashed] (1.1,-.1) -- (-.1,1.1) ;
\draw [very thick,green!80!black,dashed] (2.1,-.1) -- (-.1,2.1) ;
\draw [very thick,green!80!black,dashed] (2.1,.9) -- (.9,2.1) ;
\draw [very thick,green!80!black,->] (1.5,1.5) -- +(.2,.2) ;
\end{tikzpicture}
\quad
\begin{tikzpicture}[xscale=-1]
\draw [very thick,blue] (-.1,0) -- (2.1,0) ;
\draw [very thick,blue] (-.1,1) -- (2.1,1) ;
\draw [very thick,blue] (-.1,2) -- (2.1,2) ;
\draw [very thick,blue,->] (0.5,1) -- +(0,.2) ;

\draw [very thick,red] (0,-.1) -- (0,2.1) ;
\draw [very thick,red] (1,-.1) -- (1,2.1) ;
\draw [very thick,red] (2,-.1) -- (2,2.1) ;
\draw [very thick,red,->] (1,1.5) -- +(.2,0) ;

\draw [very thick,green!80!black,dashed] (1.1,-.1) -- (-.1,1.1) ;
\draw [very thick,green!80!black,dashed] (2.1,-.1) -- (-.1,2.1) ;
\draw [very thick,green!80!black,dashed] (2.1,.9) -- (.9,2.1) ;
\draw [very thick,green!80!black,->] (1.5,1.5) -- +(.2,.2) ;
\end{tikzpicture}
\caption{Addition of untwisted 1-forms}
\label{fig_add_untw_1forms}
\end{subfigure}
\quad
\begin{subfigure}[t]{0.14\textwidth}
\begin{tikzpicture}
\draw [very thick,blue] (-.1,0) -- (2.1,0) ;
\draw [very thick,blue] (-.1,1) -- (2.1,1) ;
\draw [very thick,blue] (-.1,2) -- (2.1,2) ;
\draw [very thick,blue,->] (1.5,1) -- +(.2,0) ;

\draw [very thick,red] (0,-.1) -- (0,2.1) ;
\draw [very thick,red] (1,-.1) -- (1,2.1) ;
\draw [very thick,red] (2,-.1) -- (2,2.1) ;
\draw [very thick,red,->] (1,1.5) -- +(0,.2) ;

\draw [very thick,green!80!black,dashed] (0.9,-.1) -- (2.1,1.1) ;
\draw [very thick,green!80!black,dashed] (-.1,-.1) -- (2.1,2.1) ;
\draw [very thick,green!80!black,dashed] (-.1,.9) -- (1.1,2.1) ;
\draw [very thick,green!80!black,->] (1.5,1.5) -- +(.2,.2) ;
\end{tikzpicture}
\caption{Addition of twisted 1-forms}
\label{fig_add_twist_1forms}
\end{subfigure}
\caption[Addition of 1-forms in 2-dim]{Addition
  of untwisted (left two) and twisted (right) 1-forms in 2-dim. The
  solid 1-form-submanifold (blue and red) are added to produce the dashed
  1-form-submanifold (green).}
\label{fig_add_oriented_1forms}
\end{figure}

 It is only possible to add forms of the same degree and
twistedness. Thus two untwisted $p$-forms sum to an untwisted $p$-form
and two twisted $p$-forms sum to a twisted $p$-form.  
In contrast to multiplying
by a scalar, adding two closed forms gives rise to a closed form.

Although algebraically adding forms is trivial, the visualisation of
the addition of two $p$-forms is surprisingly complicated. Figure
\ref{fig_add_oriented_1forms} show how to add 1-forms in 2-dimension.


\subsection{The wedge product also known as the exterior product.}
\label{ch_CForms_Wedge}


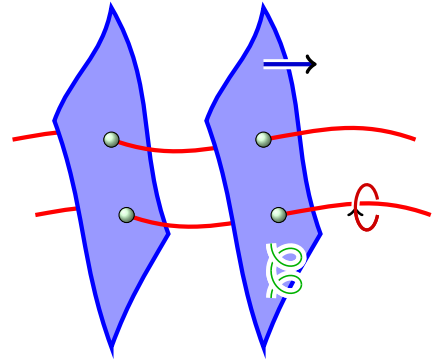
\begin{figure}[tb]
\centering
\begin{tikzpicture}
\draw[ultra thick,draw=red] (-1,2.5) to[out=10,in=160] (0.2,2.5) ;
\draw[ultra thick,draw=red] (-1.3,3.5) to[out=10,in=160] (0.,3.5) ;
\filldraw[ultra thick,draw=blue,fill=blue!40!white,scale=1.5] (0,0.5) 
to [out=110,in=-70] +(-.5,2) 
to [out=70,in=-100] +(0.5,1) 
to [out=-60,in=110] +(0.5,-2) 
to [out=-120,in=70] cycle
;
\draw[ultra thick,draw=red] (0.2,2.5) to[out=-20,in=190] (2.2,2.5) ;
\draw[ultra thick,draw=red] (0.,3.5) to[out=-20,in=190]  (2.,3.5) ;
\filldraw[ultra thick,draw=blue,fill=blue!40!white,xshift=2cm,scale=1.5] (0,0.5) 
to [out=110,in=-70] +(-.5,2) 
to [out=70,in=-100] +(0.5,1) 
to [out=-60,in=110] +(0.5,-2) 
to [out=-120,in=70] cycle
;
\draw[ultra thick,draw=red] (2.2,2.5) to[out=10,in=160] +(2,0) ;
\draw[ultra thick,draw=red] (2.,3.5) to[out=10,in=160] +(2,0) ;
\shadedraw[ball color=green!20] (0.2,2.5) circle (0.1) ;
\shadedraw[ball color=green!20] (0.,3.5) circle (0.1) ;
\shadedraw[ball color=green!20] (2.2,2.5) circle (0.1) ;
\shadedraw[ball color=green!20] (2.,3.5) circle (0.1) ;


\draw[ultra thick,draw=white,double=white]  (2,4.5) -- +(.5,0) ;
\draw[ultra thick,draw=blue!80!black,->] (2,4.5) -- +(.7,0) ;

\draw[very thick,draw=white,double=white,xscale=0.5] (7,2.5) arc
(-20:-340:0.3) ;
\draw[very thick,draw=red!80!black,xscale=0.5,->] (7,2.5) arc
(-20:-180:0.3) ;
\draw[very thick,draw=red!80!black,xscale=0.5] (7,2.5) arc
(-20:-340:0.3) ;

\draw[ultra thick,draw=white,double=green!70!black,shift={(2.3,1)},xscale=-.4,yscale=0.4]
(0.5,1) \foreach \x in {0,10,...,720}
{ to ( {0.5*cos(\x)}, {0.0025*\x+0.5*sin(\x)+1} ) } ;
\draw[ultra thick,draw=white,double=green!70!black,shift={(2.3,1)},xscale=-.4,yscale=0.4]
(-.5,0.0025*180+1) \foreach \x in {180,190,...,360}
{ to ( {0.5*cos(\x)}, {0.0025*\x+0.5*sin(\x)+1} ) } ;
\draw[ultra thick,draw=white,double=green!70!black,shift={(2.3,1)},xscale=-.4,yscale=0.4]
(-.5,0.0025*540+1) \foreach \x in {540,550,...,720}
{ to ( {0.5*cos(\x)}, {0.0025*\x+0.5*sin(\x)+1} ) } ;

\end{tikzpicture}
\caption[Wedge product of untwisted 1-form and 2-form.]{In
  3-dimensions, the wedge product of an untwisted 1-form (blue) with a
  untwisted 2-form (red) to give a untwisted 3-form (green).}
\label{fig_wedge_untw1_untw2_3dim}
\end{figure}

\begin{figure}[tb]
\centering
\begin{tikzpicture}[xscale=1.5,yscale=-1.5]

\filldraw[ultra thick,draw=red,fill=red!40!white]
(-1.2,3) -- (-.7,1.2)  -- (0.3,1.2)  to [out=120,in=-80] (-.2,3) -- cycle;

\filldraw[ultra thick,draw=blue,fill=blue!40!white] (0,0.5) 
to [out=110,in=-70] +(-.5,2) 
to [out=70,in=-100] +(0.5,1) 
to [out=-60,in=110] +(0.5,-2) 
to [out=-120,in=70] cycle
;
\filldraw[ultra thick,draw=red,fill=red!40!white]
(0.3,1.2)  to [out=110,in=-70] (-.2,3) -- (1.1,3)  -- (1.6,1.1)--  cycle;

\draw[ultra thick,draw=green]
(0.3,1.2)  to [out=110,in=-70] (-.2,3) ;

\filldraw[ultra thick,draw=blue,fill=blue!40!white,xshift=1.33cm] (0,0.5) 
to [out=110,in=-70] +(-.5,2) 
to [out=70,in=-100] +(0.5,1) 
to [out=-60,in=110] +(0.5,-2) 
to [out=-120,in=70] cycle
;
\filldraw[ultra thick,draw=red,fill=red!40!white]
(1.1,3) to [out=-70,in=110]  (1.6,1.1) -- (2.6,1.1) to [out=110,in=-70] (2.1,3) -- cycle ;

\draw[ultra thick,draw=green]
(1.1,3) to [out=-70,in=110]  (1.6,1.1) ;

\draw[ultra thick,draw=white,double=white]  (1.3,0.8) -- +(.4,0) ;
\draw[ultra thick,draw=blue!80!black,->] (1.3,0.8) -- +(.4,0) ;

\draw[ultra thick,draw=white,double=white]  (2.3,1.3) -- +(0,-.4) ;
\draw[ultra thick,draw=red!80!black,->] (2.3,1.3) -- +(0,-.4) ;

\draw[ultra thick,draw=white,double=white,yscale=0.7] (1.5,2/0.7) arc (0:-95:.3) ;
\draw[ultra thick,draw=green!60!black,->,yscale=0.7] (1.5,2/0.7) arc (0:-85:.3) ;
\draw[ultra thick,draw=green!60!black,yscale=0.7] (1.5,2/0.7) arc (0:-95:.3) ;

\end{tikzpicture}
\caption[Wedge product of two untwisted 2-form.]{In 3-dimensions, the
  wedge product of an untwisted 1-form (vertical, blue) with an
  untwisted 1-form (horizontal, red) to give an untwisted 2-form
  (green).}
\label{fig_wedge_3dim_untw1_untw1}
\end{figure}
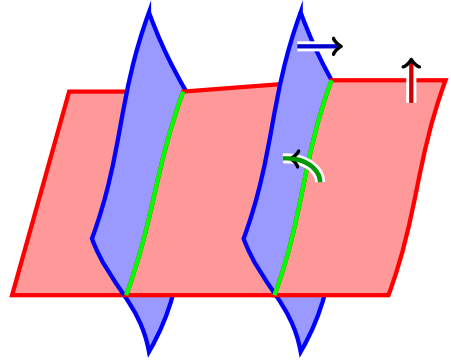

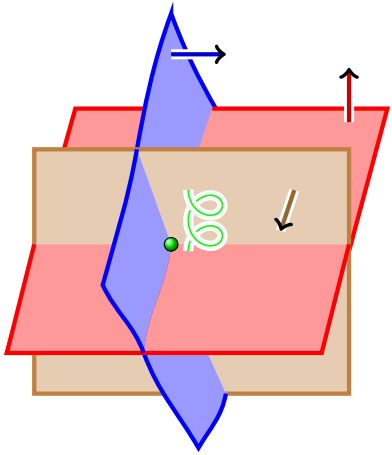
\begin{figure}[tb]
\centering
\begin{tikzpicture}[xscale=1.8,yscale=-1.8]

\filldraw[ultra thick,draw=red,fill=red!40!white]
(-1.2,3) -- (-.7,1.2)  -- (0.3,1.2)  to [out=120,in=-80]  (-.2,3) -- cycle;

\filldraw[ultra thick,draw=blue,fill=blue!40!white] (0,0.5) 
to [out=110,in=-70] +(-.5,2) 
to [out=70,in=-100] +(0.5,.5) 
to [out=-60,in=110] +(0.5,-1.5) 
to [out=-120,in=70] cycle
;

\fill[fill=red!40!white]
(0.3,1.2)  to [out=110,in=-70] (-.2,3) -- (1.1,3)  -- (1.58,1.2)--  cycle;
\draw[ultra thick,draw=red] (1.1,3)  -- (1.58,1.2) -- (0.3,1.2);

\draw[ultra thick,draw=brown,fill=brown!40!white] 
  (-1,1.5) rectangle +(2.3,1.8) ;

\fill [fill=blue!40!white] 
  (-.2,3) to[out=70,in=-120] (0.2,3.7) 
  to[out=-60,in=100]  (0.4,3.3) -- (-0.05,2.3) ;

\draw [ultra thick,draw=blue] 
  (-.2,3) to[out=70,in=-120] (0.2,3.7) 
  to[out=-60,in=100]  (0.4,3.3)  ;

\fill [fill=red!40!white] 
  (-1,2.2) -- (0,2.2) to[out=95,in=-85] (-.2,3) -- (-1.2,3) -- cycle  ;
\draw [ultra thick,draw=red] 
   (0,2.2) to[out=95,in=-85] (-.2,3) -- (-1.2,3) -- (-1,2.2) ;

\fill [fill=red!40!white] 
  (0,2.2) -- (1.3,2.2) -- (1.1,3) -- (-.2,3) to[out=-85,in=95]  cycle  ;

\draw [ultra thick,draw=red] 
  (0,2.2)  (1.3,2.2) -- (1.1,3) -- (-.2,3)   ;

\fill[fill=blue!40!white] 
  (-0.25,1.5) -- (0,2.2) to[out=95,in=-85] (-.2,3) 
to[out=-110,in=65] (-.5,2.5) 
to[out=-75,in=100] cycle  ;

\draw [ultra thick,draw=blue,fill=blue!40!white] 
  (-0.25,1.5)  (0,2.2)  (-.2,3) 
to[out=-110,in=65] (-.5,2.5) 
to[out=-75,in=100] (-0.25,1.5)  ;

 \shadedraw[ball color=green] (0,2.2) circle (.05) ;

\draw[ultra thick,draw=white,double=white]  (0,0.8) -- +(.4,0) ;
\draw[ultra thick,draw=blue!80!black,->] (0,0.8) -- +(.4,0) ;

\draw[ultra thick,draw=white,double=white]  (1.3,1.3) -- +(0,-.4) ;
\draw[ultra thick,draw=red!80!black,->] (1.3,1.3) -- +(0,-.4) ;

\draw[ultra thick,draw=white,double=white]  (.9,1.8) -- +(-.1,.3) ;
\draw[ultra thick,draw=brown!80!black,->] (.9,1.8) -- +(-.1,.3) ;


\begin{scope}[yscale=-.25,xscale=.25,shift={(-1,-10)}]
\draw[ultra thick,draw=white,double=green,xscale=-1,xshift=-2cm]
(0.5,1) \foreach \x in {0,10,...,720}
{ to ( {0.5*cos(\x)}, {0.0025*\x+0.5*sin(\x)+1} ) } ;
\draw[ultra thick,draw=white,double=green,xscale=-1,xshift=-2cm]
(-.5,0.0025*180+1) \foreach \x in {180,190,...,360}
{ to ( {0.5*cos(\x)}, {0.0025*\x+0.5*sin(\x)+1} ) } ;
\draw[ultra thick,draw=white,double=green,xscale=-1,xshift=-2cm]
(-.5,0.0025*540+1) \foreach \x in {540,550,...,720}
{ to ( {0.5*cos(\x)}, {0.0025*\x+0.5*sin(\x)+1} ) } ;
\end{scope}

\end{tikzpicture}
\caption[Wedge product of three untwisted 2-form.]{In 3-dimensions, the
  wedge product of three untwisted 1-forms. The first  vertical (blue)
  followed by the second facing reader
  (brown) then the third horizontal (red). This  gives an untwisted 3-form
  (green).}
\label{fig_wedge_3dim_untw1x3}
\end{figure}


\begin{figure}[tb]
\centering

\begin{tikzpicture}[xscale=1.5,yscale=-1.5]

\filldraw[ultra thick,draw=red,fill=red!40!white]
(-1.2,3) -- (-.7,1.2)  -- (0.3,1.2)  to [out=120,in=-80] (-.2,3) -- cycle;

\filldraw[ultra thick,draw=blue,fill=blue!40!white] (0,0.5) 
to [out=110,in=-70] +(-.5,2) 
to [out=70,in=-100] +(0.5,1) 
to [out=-60,in=110] +(0.5,-2) 
to [out=-120,in=70] cycle
;
\filldraw[ultra thick,draw=red,fill=red!40!white]
(0.3,1.2)  to [out=110,in=-70] (-.2,3) -- (1.1,3)  -- (1.6,1.1)--  cycle;

\draw[ultra thick,draw=green]
(0.3,1.2)  to [out=110,in=-70] (-.2,3) ;

\filldraw[ultra thick,draw=blue,fill=blue!40!white,xshift=1.33cm] (0,0.5) 
to [out=110,in=-70] +(-.5,2) 
to [out=70,in=-100] +(0.5,1) 
to [out=-60,in=110] +(0.5,-2) 
to [out=-120,in=70] cycle
;
\filldraw[ultra thick,draw=red,fill=red!40!white]
(1.1,3) to [out=-70,in=110]  (1.6,1.1) -- (2.6,1.1) to [out=110,in=-70] (2.1,3) -- cycle ;

\draw[ultra thick,draw=green]
(1.1,3) to [out=-70,in=110]  (1.6,1.1) ;

\draw[ultra thick,draw=blue!80!black,->] (1.2,0.9) arc(-120:90:.2) ;
\draw[ultra thick,draw=red!80!black,->] (1.9,1.3) arc(-100:120:.2) ;

\draw[ultra thick,draw=white,double=white,yscale=0.7] (1.5,2/0.7) arc (0:-95:.3) ;
\draw[ultra thick,draw=green!60!black,->,yscale=0.7] (1.5,2/0.7) arc (0:-85:.3) ;
\draw[ultra thick,draw=green!60!black,yscale=0.7] (1.5,2/0.7) arc (0:-95:.3) ;

\end{tikzpicture}
\caption[Wedge product of twisted forms]
{In 3-dimensions, the wedge product of the twisted 1-form
  (vertical, blue) with the twisted 1-form (horizontal, red) 
  to give the green untwisted 2-form.}
\label{fig_wedge_1tw_1tw}
\end{figure}
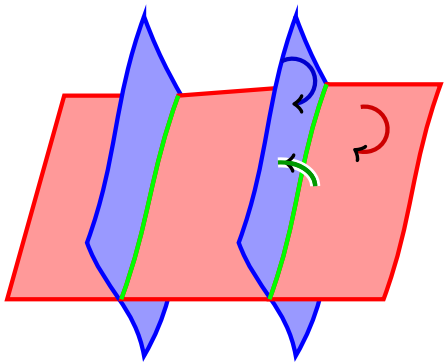


\begin{figure}[tb]
\centering

\begin{tikzpicture}[xscale=1.5,yscale=-1.5]

\filldraw[ultra thick,draw=red,fill=red!40!white]
(-1.2,3) -- (-.7,1.2)  -- (0.3,1.2)  to [out=120,in=-80] (-.2,3) -- cycle;

\filldraw[ultra thick,draw=blue,fill=blue!40!white] (0,0.5) 
to [out=110,in=-70] +(-.5,2) 
to [out=70,in=-100] +(0.5,1) 
to [out=-60,in=110] +(0.5,-2) 
to [out=-120,in=70] cycle
;
\filldraw[ultra thick,draw=red,fill=red!40!white]
(0.3,1.2)  to [out=110,in=-70] (-.2,3) -- (1.1,3)  -- (1.6,1.1)--  cycle;

\draw[line width=3pt,draw=green!70!white]
(0.3,1.2)  to [out=110,in=-70] (-.2,3) ;

\filldraw[ultra thick,draw=blue,fill=blue!40!white,xshift=1.33cm] (0,0.5) 
to [out=110,in=-70] +(-.5,2) 
to [out=70,in=-100] +(0.5,1) 
to [out=-60,in=110] +(0.5,-2) 
to [out=-120,in=70] cycle
;
\filldraw[ultra thick,draw=red,fill=red!40!white]
(1.1,3) to [out=-70,in=110]  (1.6,1.1) -- (2.6,1.1) to [out=110,in=-70] (2.1,3) -- cycle ;

\draw[line width=3pt,draw=green!70!white,->-=.7]
(1.1,3) to [out=-70,in=110]  (1.6,1.1) ;

\draw[ultra thick,draw=white,double=white]  (1.3,0.8) -- +(.4,0) ;
\draw[ultra thick,draw=blue!80!black,->] (1.3,0.8) -- +(.4,0) ;
\draw[ultra thick,draw=red!80!black,->] (1.9,1.3) arc(-100:120:.2) ;
\end{tikzpicture}
\caption[Wedge product of a twisted and an untwisted form.]
{In 3-dimensions, the wedge product of the untwisted 1-form (blue)
  with the twisted 1-form (red) to give the green twisted 2-form.}
\label{fig_wedge_1tw_1untw}
\end{figure}
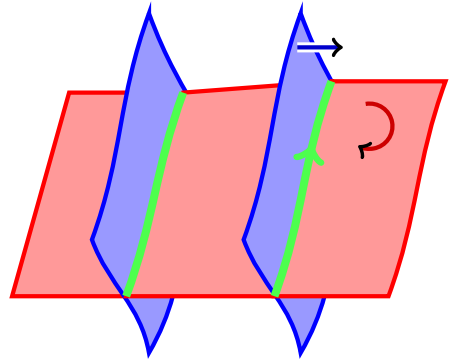


In order to take the wedge product of two forms simply requires taking
their intersection, figures
\ref{fig_wedge_untw1_untw2_3dim}-\ref{fig_wedge_1tw_1untw}. In places where
the two form-submanifolds are tangential then the resulting wedge
product vanishes. 

The wedge product of two untwisted forms is an untwisted form. Simply
concatenate the two orientations. See figures
\ref{fig_wedge_untw1_untw2_3dim} and
\ref{fig_wedge_3dim_untw1_untw1}. The wedge product for three untwisted
1-forms is given by the intersection of three $(n-1)$-dim
form-submanifolds, figure \ref{fig_wedge_3dim_untw1x3}.

The wedge product of two twisted forms is an untwisted form (figure
\ref{fig_wedge_1tw_1tw}) and the wedge product of
a twisted form and an untwisted form is a twisted form (figure
\ref{fig_wedge_1tw_1untw}). However the rules of
this is more complicated than simply concatenation. The only
guaranteed method is to choose an orientation for the manifold, then
convert the twisted forms into untwisted forms, take the wedge product
and then apply a twisting if necessary. One can check that the result
is invariant under the choice of orientation of the manifold. However
it may depend on the twisting convention described in section
\ref{ch_Orientation_Twisting}.


\subsection{Integration.}
\label{ch_CForms_Intergration}

\begin{figure}[tb]
\centering
\begin{tikzpicture}[scale=1.6]
\draw[->-=.8,ultra thick] (0,2.5) to [out=110,in=-70] (0,3.5) ; 
\draw[->-=.8,ultra thick] (0.3,2.6) to [out=110,in=-80] (0.3,3.5) ; 
\draw[->-=.8,ultra thick] (1,3) to [out=110,in=-80] (1,3.8) ; 
\draw[ultra thick] (2,2.5) to [out=110,in=-90] (2,3.8) ; 

\shadedraw[ball color=green!20] (1.1,2) circle [radius=1.3] ;

\draw[ultra thick] (0,0.5) to [out=110,in=-80] (0,1.5) ; 
\draw[ultra thick,opacity=0.3] (0,1.5) to [out=100,in=-70] (0,2.5) ; 
\fill (0,1.5) circle [radius=0.05] ;

\draw[ultra thick] (0.3,0.5) to [out=110,in=-70] (0.3,1.6) ;
\fill (0.3,1.6) circle [radius=0.05] ;
\draw[ultra thick,opacity=0.3] (0.3,1.6) to [out=100,in=-70] (0.3,2.6) ; 

\draw[ultra thick] (1,0.5) to [out=110,in=-80] (1,1.7) ;
\fill (1,1.7) circle [radius=0.05] ;
\draw[ultra thick,opacity=0.3] (1,1.7) to [out=100,in=-70] (1,3) ; 

\draw[-<-=.3,ultra thick] (2,0.5) to [out=110,in=-70] (2,1.6) ;
\fill (2,1.6) circle [radius=0.05] ;
\draw[ultra thick,opacity=0.3] (2,1.6) to [out=100,in=-70] (2,2.5) ; 

\draw[ultra thick,opacity=0.3,rotate=-10] (1.1,2) ellipse (0.2 and 0.5) ;
\draw[ultra thick,opacity=0.3,rotate=-10] (0.85,2.2) -- +(.1,.1) -- +(.1,-.1) ;

\draw[-<-=.3,ultra thick] (2.8,0.8) to [out=130,in=-80] +(0,3) ;
\draw[-<-=.3,ultra thick] (3,1.0) to [out=120,in=-70] +(0,3) ;
\draw[-<-=.3,ultra thick] (3.2,1.2) to [out=120,in=-70] +(0,3) ;
\draw[->-=.5,ultra thick] (-0.5,0.5) to [out=110,in=-80] +(0,3.4) ;

\end{tikzpicture}
\caption[Integration of a 2-form.]{Integration of a twisted 2-form on
  a 2-sphere embedded in a 3-dimensional manifold. Observe that every
  curve that enters the sphere also leaves.}
\label{fig_int_2-form_field}
\end{figure}
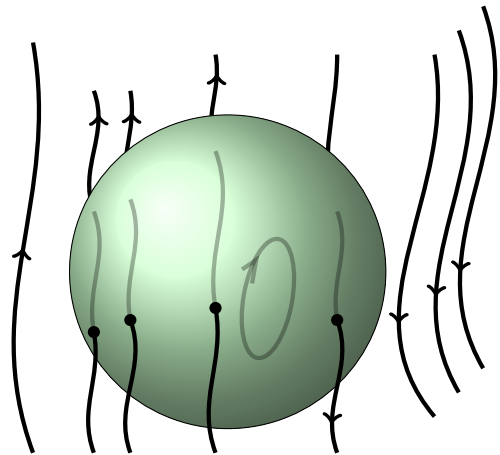


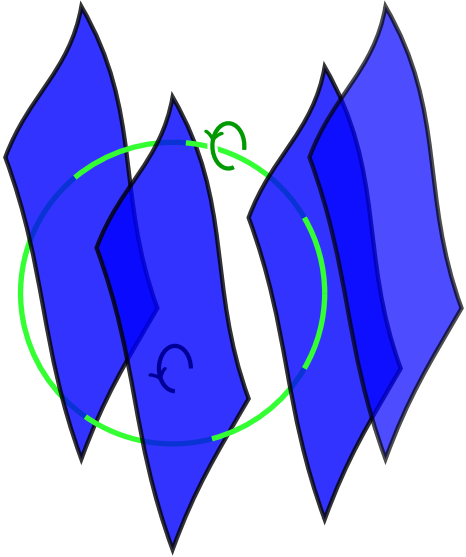
\begin{figure}[tb]
\begin{tikzpicture}[scale=2]
\draw[line width=2pt,green!80] (1,2) +(130:1) arc(130:235:1) ;
\filldraw[ultra thick,draw=black,fill=blue,opacity=0.8] (0.4,0.9) 
to [out=110,in=-70] +(-.5,2) 
to [out=70,in=-100] +(0.5,1) 
to [out=-60,in=110] +(0.5,-2) 
to [out=-120,in=70] cycle
;
\draw[line width=2pt,green!80] (1,2) +(130:1) arc(130:85:1) ;
\draw[line width=2pt,green!80] (1,2) +(235:1) arc(235:285:1) ;
\filldraw[ultra thick,draw=black,fill=blue,opacity=0.8] (1.0,0.3) 
to [out=110,in=-70] +(-.5,2) 
to [out=70,in=-100] +(0.5,1) 
to [out=-60,in=110] +(0.5,-2) 
to [out=-120,in=70] cycle
;
\draw[line width=2pt,green!80] (1,2) +(85:1) arc(85:30:1) ;
\draw[line width=2pt,green!80] (1,2) +(285:1) arc(285:330:1) ;
\filldraw[ultra thick,draw=black,fill=blue,opacity=0.8] (2,0.5) 
to [out=110,in=-70] +(-.5,2) 
to [out=70,in=-100] +(0.5,1) 
to [out=-60,in=110] +(0.5,-2) 
to [out=-120,in=70] cycle
;
\draw[line width=2pt,green!80] (1,2) +(30:1) arc(30:0:1) ;
\draw[line width=2pt,green!80] (1,2) +(330:1) arc(330:360:1) ;
\filldraw[ultra thick,draw=black,fill=blue,opacity=0.7] (2.4,0.9) 
to [out=110,in=-70] +(-.5,2) 
to [out=70,in=-100] +(0.5,1) 
to [out=-60,in=110] +(0.5,-2) 
to [out=-120,in=70] cycle
;
\draw[ultra thick,blue!60!black,->-=.8,xscale=.7] (1.6,1.5) arc (0:270:0.15) ;
\draw[ultra thick,double=green!60!black,white,xscale=.7] (2.1,2.95) arc (-10:290:0.15) ;
\draw[ultra thick,green!60!black,xscale=.7,->-=.6] (2.1,2.95) arc (-10:290:0.15) ;
\end{tikzpicture}

\caption[Integration of a 1-form.]{Integration of a twisted 1-form on
  a twisted circle embedded in a 3-dimensional manifold. Observe that
  every surface that crosses the circle must cross again in the
  opposite direction.}
\label{fig_int_1-form_field}
\end{figure}


We can only integrate twisted top-forms, i.e. $n$-forms on an
$n$-dimensional manifold, figure \ref{fig_n-form_field}. To do this we
simply add the number of positive dots and subtract the number of
negative dots. Then take the limit the the density tends to infinity
and the value tends to zero.  Let figure \ref{fig_n-form_field} refer
to a 2-dimensional manifold.  Thus in figure
\ref{fig_low-density}, assuming each dot
value has 1, then its integral equals 9. In figure
\ref{fig_high-density} each dot must have
value $\tfrac14$ and the integral is $6\tfrac12$. As the density of dots
increases and their value tends to zero, this sum will converge. Contrast
this to the case of the scalar field, figure \ref{fig_Scalar_field},
which does not converge as the density increases. In order to
integrate a scalar field it is necessary to multiply it by a measure.

To integrate an untwisted top-form it is necessary to use the
orientation of the manifold to convert the untwisted top-form into a
twisted top-form. If the manifold is non orientable, such
as the Möbius strip then one can only integrate twisted top-form.

What about integrating forms of lower dimension? To integrate a
$p$-form we must integrate it over a $p$-dimensional submanifold. In
addition untwisted forms are integrated over untwisted submanifolds and
likewise for twisted forms and twisted submanifolds; with plus if the
orientations agree and minus otherwise. We will see a generalisation
of this in section \ref{ch_PullBacks} on pullbacks.


\subsection{Conservation Laws}
\label{sch_conslaws}

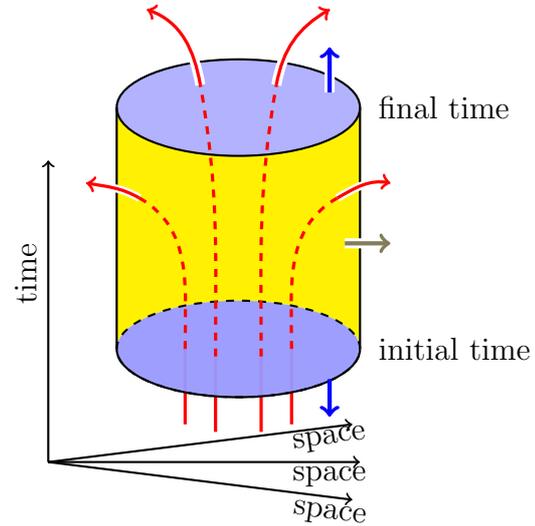
\begin{figure}[tb]
\centering
\begin{tikzpicture}
\draw [thick,->] (-2.5,-1.5) -- + (0,4) ;
\draw (-2.8,1) node {\rotatebox{90}{time}} ; 
\draw [thick,->] (-2.5,-1.5) -- + (4,.5) ;
\draw [thick,->] (-2.5,-1.5) -- + (4,-.5) ;
\draw [thick,->] (-2.5,-1.5) -- + (4.1,0) ;
\draw (1.2,-1.7) node {{space}} ; 
\draw (1.2,-1.2) node[rotate=8] {{space}} ; 
\draw (1.2,-2.2) node[rotate=-8] {{space}} ; 

\draw (1.7,0) node[right] {initial time} ;
\draw (1.7,3.2) node[right] {final time} ;

\draw[very thick, red] (-.7,0) -- +(0,-1);
\draw[very thick, red] (-.3,-.1) -- +(0,-1);
\draw[very thick, red] (.3,-.1) -- +(0,-1);
\draw[very thick, red] (.7,0) -- +(0,-1);

\fill[fill=yellow,yscale=0.4] (1.6,0) -- ++(0,3.2/.4)  -- ++(-3.2,0)  -- ++(0,-3.2/.4)
arc (180:0:1.6) ;
\filldraw[thick,yscale=0.4,dashed,fill=blue!40,opacity=.95] (0,0) circle (1.6);
\draw[thick,yscale=0.4] (1.6,0) arc (0:-180:1.6);
\filldraw[thick,yscale=0.4,fill=blue!30] (0,8) circle (1.6);
\draw[thick] (-1.6,0) -- +(0,3.2) ;
\draw[thick] (1.6,0) -- +(0,3.2) ;

\draw[very thick, red,dashed] (-.7,0) to[out=90,in=-30] (-1.3,2)
to[out=150,in=-10] +(-.7,.2);
\draw[very thick, double=red,white] (-1.3,2) to[out=150,in=-10] +(-.7,.2);
\draw[very thick, red,->] (-1.3,2) to[out=150,in=-10] +(-.7,.2);

\draw[very thick, red,dashed] (.7,0) to[out=90,in=-150] (1.3,2)
to[out=30,in=170] +(.7,.2);
\draw[very thick, double=white,white] (1.3,2) to[out=30,in=170] +(.7,.2);
\draw[very thick, red,->] (1.3,2) to[out=30,in=170] +(.7,.2);

\draw[very thick, red,dashed] (-.3,-0.1) to[out=90,in=-80] (-.5,3.5)
to[out=100,in=-20] +(-.7,1);
\draw[ultra thick, double=red,white] (-.5,3.5) to[out=100,in=-20] +(-.7,1);
\draw[very thick, red,->] (-.5,3.5) to[out=100,in=-20] +(-.7,1);

\draw[very thick, red,dashed] (.3,-0.1) to[out=90,in=-100] 
(.5,3.5) to[out=80,in=-150] +(.7,1);
\draw[ultra thick, double=red,white] (.5,3.5) to[out=80,in=-150] +(.7,1);
\draw[very thick, red,->] (.5,3.5) to[out=80,in=-150] +(.7,1);;

\draw[ultra  thick, double=white,white] (1.2,3.4) -- +(0,.6);
\draw[ultra  thick, blue,->]  (1.2,3.4) -- +(0,.6);
\draw[ultra  thick, double=white,white] (1.4,1.4) -- +(.6,0);
\draw[ultra  thick, yellow!40!black,->]  (1.4,1.4)  -- +(.6,0);
\draw[ultra  thick, blue,->]  (1.2,-.4) -- +(0,-.5);

\end{tikzpicture}

\caption[A conservation law.]{The closed twisted 3-form (red) is
  integrated over the 3-dim surface of a cylinder with net result
  zero. Any current which enters at the initial time  (lower blue disc) must
  either leave at the final time  (upper blue disc) or through the curved
  surface of the cylinder (yellow). Note that the orientation of the
  submanifold is out (given by the yellow and blue arrows). This is
  because the sum of the integrals is zero. However the contribution
  from initial time is negative.}
\label{fig_conserved_cylinder}
\end{figure}


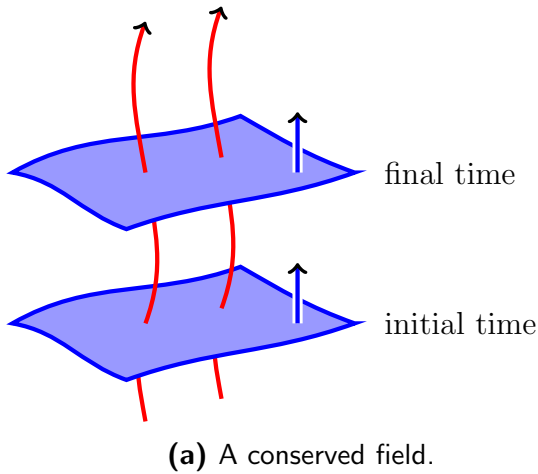
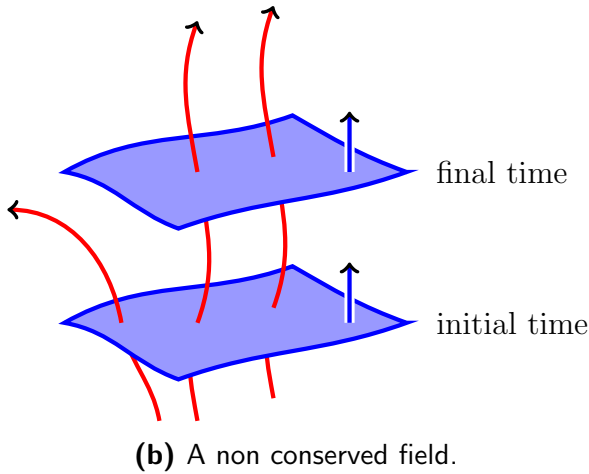
\begin{figure}[tb]
\centering
\begin{subfigure}{0.45\textwidth}
\begin{tikzpicture}[rotate=90]
\draw[ultra thick,draw=red] (-1,2.5) to[out=10,in=160] (0.2,2.5) ;
\draw[ultra thick,draw=red] (-1.3,3.5) to[out=10,in=160] (0.,3.5) ;
\filldraw[ultra thick,draw=blue,fill=blue!40!white,scale=1.5] (0,0.5) 
to [out=110,in=-70] +(-.5,2) 
to [out=70,in=-100] +(0.5,1) 
to [out=-60,in=110] +(0.5,-2) 
to [out=-120,in=70] cycle
;
\draw[ultra thick,draw=red] (0.2,2.5) to[out=-20,in=190] (2.2,2.5) ;
\draw[ultra thick,draw=red] (0.,3.5) to[out=-20,in=190]  (2.,3.5) ;
\filldraw[ultra thick,draw=blue,fill=blue!40!white,xshift=2cm,scale=1.5] (0,0.5) 
to [out=110,in=-70] +(-.5,2) 
to [out=70,in=-100] +(0.5,1) 
to [out=-60,in=110] +(0.5,-2) 
to [out=-120,in=70] cycle
;
\draw[ultra thick,draw=white,double=white] (0,+1.5) -- +(0.8,0) ; 
\draw[ultra thick,draw=blue,->] (0,+1.5) -- +(0.8,0) ; 
\draw[ultra thick,draw=white,double=white] (2,+1.5) -- +(0.8,0) ; 
\draw[ultra thick,draw=blue,->] (2,+1.5) -- +(.8,0) ; 

\draw[ultra thick,draw=red,->] (2.2,2.5) to[out=10,in=160] +(2,0) ;
\draw[ultra thick,draw=red,->] (2.,3.5) to[out=10,in=160] +(2,0) ;

\draw (0,.5) node[right] {initial time} ;
\draw (2.,.5) node[right] {final time} ;
\end{tikzpicture}
\caption{A conserved field.}
\label{fig_conserved_no_escape}
\end{subfigure}
\qquad
\begin{subfigure}{0.45\textwidth}
\begin{tikzpicture}[rotate=90]
\draw[ultra thick,draw=red] (-1,2.5) to[out=10,in=160] (0.2,2.5) ;
\draw[ultra thick,draw=red] (-1.3,3.5) to[out=10,in=160] (0.,3.5) ;
\draw[ultra thick,draw=red] (-1.3,4) to[out=10,in=-170] (0.,4.5) ;
\filldraw[ultra thick,draw=blue,fill=blue!40!white,scale=1.5] (0,0.5) 
to [out=110,in=-70] +(-.5,2) 
to [out=70,in=-100] +(0.5,1) 
to [out=-60,in=110] +(0.5,-2) 
to [out=-120,in=70] cycle
;
\draw[ultra thick,draw=red] (0.2,2.5) to[out=-20,in=190] (2.2,2.5) ;
\draw[ultra thick,draw=red] (0.,3.5) to[out=-20,in=190]  (2.,3.5) ;
\draw[ultra thick,draw=red,->] (0.,4.5)  to[out=10,in=-90] (1.5,6);
\filldraw[ultra thick,draw=blue,fill=blue!40!white,xshift=2cm,scale=1.5] (0,0.5) 
to [out=110,in=-70] +(-.5,2) 
to [out=70,in=-100] +(0.5,1) 
to [out=-60,in=110] +(0.5,-2) 
to [out=-120,in=70] cycle
;
\draw[ultra thick,draw=white,double=white] (0,+1.5) -- +(0.8,0) ; 
\draw[ultra thick,draw=blue,->] (0,+1.5) -- +(0.8,0) ; 
\draw[ultra thick,draw=white,double=white] (2,+1.5) -- +(0.8,0) ; 
\draw[ultra thick,draw=blue,->] (2,+1.5) -- +(.8,0) ; 

\draw[ultra thick,draw=red,->] (2.2,2.5) to[out=10,in=160] +(2,0) ;
\draw[ultra thick,draw=red,->] (2.,3.5) to[out=10,in=160] +(2,0) ;
\draw (0,.5) node[right] {initial time} ;
\draw (2.,.5) node[right] {final time} ;
\end{tikzpicture}
\caption{A non conserved field.}
\label{fig_conserved_escape}
\end{subfigure}
\caption[Conservation law for all space.]
{In 4-dimensions, integration of the twisted 3-form (red),
  showing it is conserved as long as nothing ``escapes to
  infinity''. The 3-dimensional submanifolds (blue) correspond to
  ``equal time slice''. 
Note that, comparing with figure
  \ref{fig_conserved_cylinder} we have let the curved surface (yellow)
pass to infinity. We have also swapped the orientation of the initial
time-slice so that it has the same orientation as the finial time slice.}
\label{fig_conserved}
\end{figure}

We can see from figures \ref{fig_int_2-form_field} and
\ref{fig_int_1-form_field} that every curve or surface that intersects
the submanifold also leaves. Therefore the integral over the sphere in
figure \ref{fig_int_2-form_field} and over the circle in figure
\ref{fig_int_1-form_field} are both zero. This gives rise to
conservation laws which are very important in physics. For example the
conservation of charge or the conservation of energy. In figure
\ref{fig_conserved_cylinder} we see that the closed twisted 3-form is
integrated over the surface of a 3 dimensional cylinder. In this
context we refer to the closed 3-form as a current.

However there are three conditions on the submanifold 
which together imply that the
integral of a closed form is zero. 
\begin{jgitemize}
\item
The submanifold (of dimension $p$) 
must be the boundary of another submanifold (of dimension $p+1$).
\item
The $(p+1)$-dimensional manifold must not go off to infinity.
\item
The larger manifold in which all the submanifolds are embedded 
must not have any ``holes''.
\end{jgitemize}
See section \ref{sch_derham} below for a comment about holes. 

In order to convert the integral of a closed current into a
conservation over all time it is necessary to expand the curved
surface of the cylinder out to infinity. Thus the initial and final
discs become all of space. One must therefore guarantee that no
current ``escapes to infinity''. See figure \ref{fig_conserved}.


\subsection{DeRham Cohomology.}
\label{sch_derham}

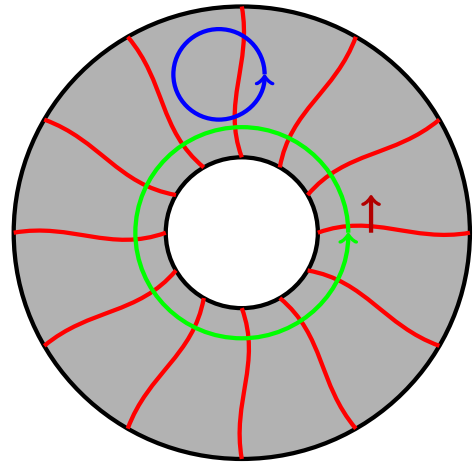
\begin{figure}[tb]
\centering
\begin{tikzpicture}
\filldraw [ultra thick,fill=black!30] (0,0) circle (3) ;
\filldraw [ultra thick,fill=white] (0,0) circle (1) ;

\foreach \x in {0,30,...,330} {
\draw [ultra thick,red,rotate={\x}] (1,0) to [out=20,in=-170] (3,0) ;
}
\draw (1.7,0) [ultra thick,red!70!black,->] -- +(0,0.5) ;
\draw [ultra thick,green,->] (1.4,0) arc (0:360:1.4) ;
\draw [ultra thick,blue,->] (0.3,2.1) arc (0:360:.6) ;
\end{tikzpicture}
\caption[DeRham Cohomology on an annulus.]{A annulus which is a 2-dim
  manifold with a hole. When the untwisted 1-form (red) is integrated
  in a circle (blue) which does not enclose the hole the result is
  zero. By contrast when one integrates over circle (green) which does
  enclose the hole the result is non zero.}
\label{fig_int_circle_around_hole}
\end{figure}

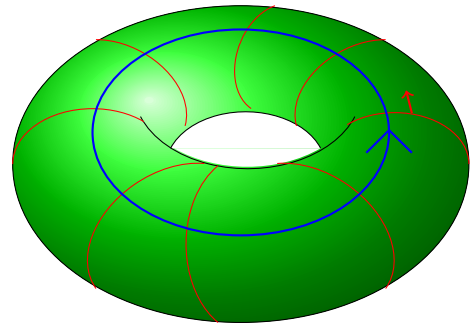
\begin{figure}[tb]
\centering
\begin{tikzpicture}[scale=1.5]
\shadedraw[ball color=green,yscale=0.7] (0,0) circle (2) ;
\filldraw[fill=white,yscale=0.7] (0.7,0.2) arc (20:160:0.7) ;
\fill[white,yscale=0.7] (0.7,0.2) arc (-50:-130:1) ;
\draw[yscale=0.7] (1,0.6) arc (-20:-160:1) ;
\draw (2,0) [yscale=0.7,red] arc (0:130:0.65) ;
\draw (1.27,1.1) [rotate=30,yscale=0.7,red]  arc (50:170:0.65) ;
\draw (.3,1.4) [rotate=70,yscale=0.7,red]  arc (30:170:0.5) ;
\draw (1.27,-1.1) [rotate=-35,yscale=0.7,red]  arc (0:140:0.8) ;
\draw[xscale=-1] (2,0) [yscale=0.7,red] arc (0:130:0.7) ;
\draw[xscale=-1] (1.27,1.1) [rotate=30,yscale=0.7,red]  arc (50:170:0.65) ;
\draw[xscale=-1] (1.27,-1.1) [rotate=-35,yscale=0.7,red]  arc (0:140:0.8) ;
\draw (-0.2,-1.4) [rotate=-90,yscale=-0.7,red]  arc (30:150:0.8) ;
\draw[thick,red,->] (1.5,0.45) -- +(-0.05,0.2) ;
\draw[thick,blue,yscale=0.7,->] (0,0.4) circle (1.3) ;
\draw[thick,blue] (1.1,0.1) -- +(0.2,0.2) -- +(0.4,0) ;
\end{tikzpicture}
\caption[DeRham Cohomology on a Torus.]  {On a 2-dim torus,
  integrating the 1-form (red) over the circle (blue) is not zero,
  despite the fact that the 1-form is closed.}
\label{fig_derham}
\end{figure}

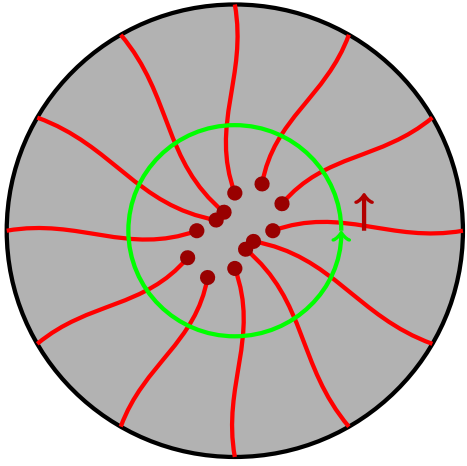
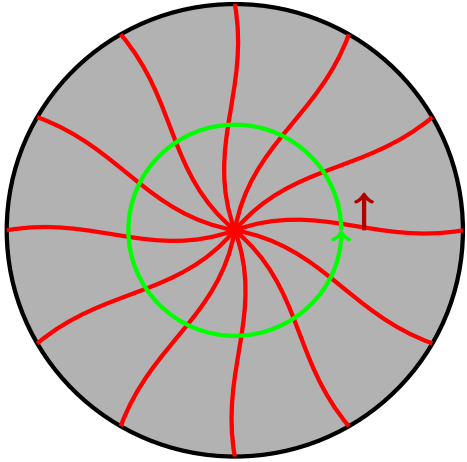
\begin{figure}[tb]
\begin{subfigure}{0.45\textwidth}
\centering
\begin{tikzpicture}
\filldraw [ultra thick,fill=black!30] (0,0) circle (3) ;

\foreach \x in {0,30,...,330} {
\draw [ultra thick,red,rotate={\x}] ({(2+sin(2*\x))/4},0) 
     to [out=20,in=-170] (3,0) ;
\fill [ultra thick,red!60!black,rotate={\x}] ({(2+sin(2*\x))/4},0)
     circle (0.1) ;
}
\draw (1.7,0) [ultra thick,red!70!black,->] -- +(0,0.5) ;
\draw [ultra thick,green,->] (1.4,0) arc (0:360:1.4) ;
\end{tikzpicture}
\caption{In this case the form-submanifolds terminate. This
  corresponds to the 1-form (red) not being closed.}
\label{fig_int_circle_non_closed}
\end{subfigure}
\quad
\begin{subfigure}{0.45\textwidth}
\centering
\begin{tikzpicture}
\filldraw [ultra thick,fill=black!30] (0,0) circle (3) ;

\foreach \x in {0,30,...,330} {
\draw [ultra thick,red,rotate={\x}] (0,0) 
     to [out=20,in=-170] (3,0) ;
}
\draw (1.7,0) [ultra thick,red!70!black,->] -- +(0,0.5) ;
\draw [ultra thick,green,->] (1.4,0) arc (0:360:1.4) ;
\end{tikzpicture}
\caption{In this case the form-submanifolds all intersect. This
  corresponds to the 1-form (red) not being continuous at this point.}
\label{fig_int_circle_not_cts}
\end{subfigure}
\caption[Failed non-zero integrals.]{Two attempts to create a non zero integral as in figure
  \ref{fig_int_circle_around_hole}, but on a manifold which does not
  have a hole. On the 2-dimensional disc, the 1-form (red) is
  integrated over circle (green).}
\label{fig_int_circle_nonzero}
\end{figure}

As mentioned in section \ref{sch_conslaws} in order to guarantee that
the integral of a closed $p$-form over a closed $p$-dimensional
submanifold is zero we need that the larger manifold in which all the
submanifolds are embedded has no ``holes''. As an example of a
manifold with holes consider the annulus, which is disc with a hole,
given in figure \ref{fig_int_circle_around_hole}. When integrating
around a circle which does not encircle the hole then the integral is
zero. However when integrating around the hole then the integral is
non-zero.

The use of integration to identify holes in a manifold is known as
DeRham Cohomology.

Note that the ``hole'' need not be removed from the manifold like a
hole is removed from a disc to create an annulus. In figure
\ref{fig_derham}, one integrates over a circle which encircles
the ``hole'' of a torus. 

If the manifold does not have any holes, then as stated the integral is
zero. It is interesting to see attempts to create a non zero integral
as in figure \ref{fig_int_circle_nonzero}. 
In figure
\ref{fig_int_circle_non_closed} the form-submanifolds terminate, 
whereas in figure
\ref{fig_int_circle_not_cts} the form-submanifolds all intersect at a
point. In both cases these correspond to non closed forms as we will
see below in section \ref{ch_NonClosed_Forms}. The difference is that
in \ref{fig_int_circle_non_closed} the 1-form is continuous, in
\ref{fig_int_circle_not_cts} it is discontinuous. We will see physical
examples of such discontinuous forms in electrostatics and
magnetostatics, section \ref{ch_Max_statics} which correspond to the
fields generated by point and line sources such as electrons and
wires.

\subsection{Four dimensions.}
\label{ch_four-dim}

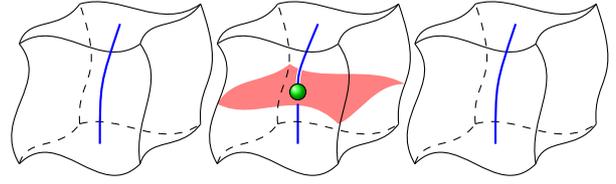
\begin{figure}[tb]
\centering
\begin{subfigure}[b]{0.31\figsize}
\begin{tikzpicture}[scale=0.53]
\draw (0,0) .. controls (1,1) and (2,-1) .. (3,0) ;
\draw (0,3) .. controls (1,3.7) and (2,2.3) .. (3,3) ;
\draw (0,0) .. controls (-0.7,1.2) and (.8,2) .. (0,3) ;
\draw (3,0) .. controls (2.8,1.2) and (3.7,2) .. (3,3) ;
\draw (3,3) .. controls (3.2,3.7) and (4.,3.3) .. (4.5,4) ;
\draw (0,3) .. controls (0.2,3.7) and (1.,3.3) .. (1.5,4) ;
\draw (3,0) .. controls (3.2,0.7) and (4.,0.3) .. (4.5,1) ;
\draw (1.5,4) .. controls (2.5,4.3) and (3.5,3.3) .. (4.5,4) ;
\draw (4.5,1) .. controls (4.4,2.2) and (5.2,2) .. (4.5,4) ;
\draw [dashed] (0,0) .. controls (0.2,0.7) and (1.,0.3) .. (1.5,1) ;
\draw [dashed] (1.5,1) .. controls (2.5,1.3) and (3.5,0.3) .. (4.5,1) ;
\draw [dashed] (1.5,1) .. controls (1.4,2.2) and (2.2,2) .. (1.5,4) ;
\draw [color=blue,thick] (2.0,0.5) .. controls (2.0,2) .. (2.5,3.5) ;
\end{tikzpicture}
\caption{$t<0$}
\end{subfigure}
\ 
\begin{subfigure}[b]{0.31\figsize}
\begin{tikzpicture}[scale=0.53]
\fill [color=red,opacity=0.5] (0.0,1.5) .. controls (1,1) and (2,2) .. (3,1) 
          .. controls (3.5,1.2) and (3.9,2) .. (4.65,2) 
          .. controls (3.8,2.5) and (2.5,2) .. (1.8,2.5) 
          .. controls (1.2,2) and (0.4,2.2) .. cycle ;
\draw [color=white,ultra thick] (2.0,1.8) .. controls (2.0,2.3) .. (2.5,3.5) ;
\draw [color=blue,thick] (2.0,1.8) .. controls (2.0,2.3) .. (2.5,3.5) ;
\draw [color=blue,thick] (2.0,0.5) -- (2.0,1.5) ;
\draw (0,0) .. controls (1,1) and (2,-1) .. (3,0) ;
\draw (0,3) .. controls (1,3.7) and (2,2.3) .. (3,3) ;
\draw (0,0) .. controls (-0.7,1.2) and (.8,2) .. (0,3) ;
\draw (3,0) .. controls (2.8,1.2) and (3.7,2) .. (3,3) ;
\draw (3,3) .. controls (3.2,3.7) and (4.,3.3) .. (4.5,4) ;
\draw (0,3) .. controls (0.2,3.7) and (1.,3.3) .. (1.5,4) ;
\draw (3,0) .. controls (3.2,0.7) and (4.,0.3) .. (4.5,1) ;
\draw (1.5,4) .. controls (2.5,4.3) and (3.5,3.3) .. (4.5,4) ;
\draw (4.5,1) .. controls (4.4,2.2) and (5.2,2) .. (4.5,4) ;
\draw [dashed] (0,0) .. controls (0.2,0.7) and (1.,0.3) .. (1.5,1) ;
\draw [dashed] (1.5,1) .. controls (2.5,1.3) and (3.5,0.3) .. (4.5,1) ;
\draw [dashed] (1.5,1) .. controls (1.4,2.2) and (2.2,2) .. (1.5,4) ;
\shadedraw [ball color=green] (2.0,1.8) circle (0.2) ;
\end{tikzpicture}
\caption{$t=0$}
\end{subfigure}
\ 
\begin{subfigure}[b]{0.31\figsize}
\begin{tikzpicture}[scale=0.53]
\draw (0,0) .. controls (1,1) and (2,-1) .. (3,0) ;
\draw (0,3) .. controls (1,3.7) and (2,2.3) .. (3,3) ;
\draw (0,0) .. controls (-0.7,1.2) and (.8,2) .. (0,3) ;
\draw (3,0) .. controls (2.8,1.2) and (3.7,2) .. (3,3) ;
\draw (3,3) .. controls (3.2,3.7) and (4.,3.3) .. (4.5,4) ;
\draw (0,3) .. controls (0.2,3.7) and (1.,3.3) .. (1.5,4) ;
\draw (3,0) .. controls (3.2,0.7) and (4.,0.3) .. (4.5,1) ;
\draw (1.5,4) .. controls (2.5,4.3) and (3.5,3.3) .. (4.5,4) ;
\draw (4.5,1) .. controls (4.4,2.2) and (5.2,2) .. (4.5,4) ;
\draw [dashed] (0,0) .. controls (0.2,0.7) and (1.,0.3) .. (1.5,1) ;
\draw [dashed] (1.5,1) .. controls (2.5,1.3) and (3.5,0.3) .. (4.5,1) ;
\draw [dashed] (1.5,1) .. controls (1.4,2.2) and (2.2,2) .. (1.5,4) ;
\draw [color=blue,thick] (2.0,0.5) .. controls (2.0,2) .. (2.5,3.5) ;
\end{tikzpicture}
\caption{$t>0$}
\end{subfigure}
\caption[2-form in 4-dimensions: plane picture.]{The plane
  representation of a 2-form-submanifold in 4-dimensions. The
  2-form-submanifold is represented by both the line (blue), for all
  time, and the momentary plane (red) at the time $t=0$. The self wedge
  is represented the dot (green).}
\label{fig_2form_4dim_plane}
\end{figure}

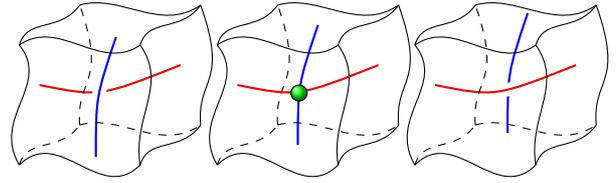
\begin{figure}[tb]
\centering
\begin{subfigure}[b]{0.31\figsize}
\begin{tikzpicture}[scale=0.53]
\draw [color=red,thick] (0.5,2.0) .. controls (2.0,1.7) .. (4,2.5) ;
\draw [line width=5pt,color=white,shift={(-.1,-.3)}] (2.0,0.5) .. controls (2.0,2) .. (2.5,3.5) ;
\draw [color=blue,thick,shift={(-.1,-.3)}] (2.0,0.5) .. controls (2.0,2) .. (2.5,3.5) ;
\draw (0,0) .. controls (1,1) and (2,-1) .. (3,0) ;
\draw (0,3) .. controls (1,3.7) and (2,2.3) .. (3,3) ;
\draw (0,0) .. controls (-0.7,1.2) and (.8,2) .. (0,3) ;
\draw (3,0) .. controls (2.8,1.2) and (3.7,2) .. (3,3) ;
\draw (3,3) .. controls (3.2,3.7) and (4.,3.3) .. (4.5,4) ;
\draw (0,3) .. controls (0.2,3.7) and (1.,3.3) .. (1.5,4) ;
\draw (3,0) .. controls (3.2,0.7) and (4.,0.3) .. (4.5,1) ;
\draw (1.5,4) .. controls (2.5,4.3) and (3.5,3.3) .. (4.5,4) ;
\draw (4.5,1) .. controls (4.4,2.2) and (5.2,2) .. (4.5,4) ;
\draw [dashed] (0,0) .. controls (0.2,0.7) and (1.,0.3) .. (1.5,1) ;
\draw [dashed] (1.5,1) .. controls (2.5,1.3) and (3.5,0.3) .. (4.5,1) ;
\draw [dashed] (1.5,1) .. controls (1.4,2.2) and (2.2,2) .. (1.5,4) ;
\end{tikzpicture}
\caption{$t<0$}
\end{subfigure}
\ 
\begin{subfigure}[b]{0.31\figsize}
\begin{tikzpicture}[scale=0.53]
\draw (0,0) .. controls (1,1) and (2,-1) .. (3,0) ;
\draw (0,3) .. controls (1,3.7) and (2,2.3) .. (3,3) ;
\draw (0,0) .. controls (-0.7,1.2) and (.8,2) .. (0,3) ;
\draw (3,0) .. controls (2.8,1.2) and (3.7,2) .. (3,3) ;
\draw (3,3) .. controls (3.2,3.7) and (4.,3.3) .. (4.5,4) ;
\draw (0,3) .. controls (0.2,3.7) and (1.,3.3) .. (1.5,4) ;
\draw (3,0) .. controls (3.2,0.7) and (4.,0.3) .. (4.5,1) ;
\draw (1.5,4) .. controls (2.5,4.3) and (3.5,3.3) .. (4.5,4) ;
\draw (4.5,1) .. controls (4.4,2.2) and (5.2,2) .. (4.5,4) ;
\draw [dashed] (0,0) .. controls (0.2,0.7) and (1.,0.3) .. (1.5,1) ;
\draw [dashed] (1.5,1) .. controls (2.5,1.3) and (3.5,0.3) .. (4.5,1) ;
\draw [dashed] (1.5,1) .. controls (1.4,2.2) and (2.2,2) .. (1.5,4) ;
\draw [color=blue,thick] (2.0,0.5) .. controls (2.0,2) .. (2.5,3.5) ;
\draw [color=red,thick] (0.5,2.0) .. controls (2.0,1.7) .. (4,2.5) ;
\shadedraw [ball color=green] (2.03,1.8) circle (0.2) ;
\end{tikzpicture}
\caption{$t=0$}
\end{subfigure}
\ 
\begin{subfigure}[b]{0.31\figsize}
\begin{tikzpicture}[scale=0.53]
\draw [color=blue,thick,shift={(.3,.3)}] (2.0,0.5) .. controls (2.0,2) .. (2.5,3.5) ;
\draw [line width=5pt,color=white] (0.5,2.0) .. controls (2.0,1.7) .. (4,2.5) ;
\draw [color=red,thick] (0.5,2.0) .. controls (2.0,1.7) .. (4,2.5) ;
\draw (0,0) .. controls (1,1) and (2,-1) .. (3,0) ;
\draw (0,3) .. controls (1,3.7) and (2,2.3) .. (3,3) ;
\draw (0,0) .. controls (-0.7,1.2) and (.8,2) .. (0,3) ;
\draw (3,0) .. controls (2.8,1.2) and (3.7,2) .. (3,3) ;
\draw (3,3) .. controls (3.2,3.7) and (4.,3.3) .. (4.5,4) ;
\draw (0,3) .. controls (0.2,3.7) and (1.,3.3) .. (1.5,4) ;
\draw (3,0) .. controls (3.2,0.7) and (4.,0.3) .. (4.5,1) ;
\draw (1.5,4) .. controls (2.5,4.3) and (3.5,3.3) .. (4.5,4) ;
\draw (4.5,1) .. controls (4.4,2.2) and (5.2,2) .. (4.5,4) ;
\draw [dashed] (0,0) .. controls (0.2,0.7) and (1.,0.3) .. (1.5,1) ;
\draw [dashed] (1.5,1) .. controls (2.5,1.3) and (3.5,0.3) .. (4.5,1) ;
\draw [dashed] (1.5,1) .. controls (1.4,2.2) and (2.2,2) .. (1.5,4) ;
\end{tikzpicture}
\caption{$t>0$}
\end{subfigure}
\caption[2-form in 4-dimensions: moving line picture.]{The
  moving line representation of a 2-form-submanifold in
  4-dimensions. The 2-form-submanifold is represented by two lines (red and
  blue). The self wedge is
  represented by the dot (green).}
\label{fig_2form_4dim_moving_line}
\end{figure}


There are significant yet obvious problems when trying to visualise
forms in 4-dimensions. However there are important reasons to try. As
we will see below, we of course live in a 4-dimensional universe and
the pictorial representation of electrodynamics would be particularly
enlightening. Secondly there are phenomena which only exist in four or
more dimensions. In lower dimensions any $p$-form, for $p\ge0$ wedged
with itself is automatically zero. This is because every
form-submanifold is tangential to itself. This does not apply to
scalar fields, viewed as 0-forms, which can be squared.  However there
are 2-forms in 4-dimensions which are none zero when wedged with
themselves. In particular, this includes the electromagnetic 2-form
which we write $F$. Thirdly, all closed form-submanifolds up to now
have a single smooth submanifold passing through each point. However
2-forms in 4-dimensions, which are non-zero when wedged with
themselves, correspond, at each point, to two  2-dimensional
form-submanifolds which intersect at that point.

One way to attempt to visualise a 2-form-submanifold in 4-dimensions,
is to use 3-dimensional space and 1-dimension of time. This is
particularly relevant for visualising spacetime. There are two
possibilities. Let's assume that the two 2-dimensional
form-submanifolds intersect at time $t=0$. Then one possibility, shown
in figure \ref{fig_2form_4dim_plane} is that for all $t\ne0$ the
form-submanifold is a line, whereas at $t=0$ it is that line and an
intersecting plane. The alternative representation is given by two
moving lines as seen in figure \ref{fig_2form_4dim_moving_line}.  The
result of wedging the 2-form with itself gives a 4-form. This top-form is
represented by a dot, which is the intersection of the
2-form-submanifold with itself.


\section{Non-closed forms.}
\label{ch_NonClosed_Forms}

\begin{figure}[tb]
\centering
\begin{tikzpicture}[scale=2.5]
\foreach \x in {1,...,9}
{ 
  \foreach \y in {0,1,...,\x}
  {
    \draw[thick] (\x*0.2,\y / \x * 1.5) -- +(.2,0) ;
  }
} ;
\draw[thick,draw=green!50!black] (0.5,0.1) 
  -- (1.7,0.1) -- (1.7,1.4) -- (0.5,1.4) -- cycle ;
\draw[->,thick,draw=green!50!black] (0.5,0.1) -- (1.0,0.1) ;
\draw[->,thick,draw=green!50!black] (1.7,0.1) -- +(0,1) ;
\draw[->,thick,draw=green!50!black] (0.5,0.1) -- +(0,1) ;
\draw[->,thick,draw=green!50!black] (0.5,1.4) -- +(0.7,0) ;
\end{tikzpicture}

\caption[Example of a non-closed 1-form (black) in 2-dim.]
{Example of a non-closed 1-form (black) in 2-dim.  
If we integrate from
bottom left to top right on the green square, we see we get a
different value depending on the path taken.}
\label{fig_integrate_nonclosed_1form}
\end{figure}
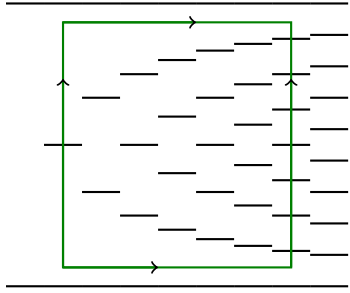

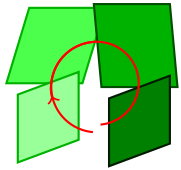
\begin{figure}
\centering
\begin{tikzpicture}
\filldraw[thick,draw=green!70!black,fill=green!70!white]
(0.15,0) -- +(1,0) -- +(1.3,1) -- +(0.3,1) -- cycle ;
\filldraw[thick,draw=green!30!black,fill=green!70!black]
(1.4,-0.05) -- +(1,0) -- +(0.9,1.1) -- +(-.1,1.1) -- cycle ;
\filldraw[thick,draw=green!10!black,fill=green!50!black]
(1.5,-.2) -- +(.8,.3) -- +(.8,-.6) -- +(0,-.9) -- cycle ;
\filldraw[thick,draw=green!70!black,fill=green!40!white]
(0.3,-.15)  -- +(.8,.3) -- +(.8,-.6) -- +(0,-.9) -- cycle ;

\draw[thick,red] (1.34,0.55) arc (90:-85:0.55) ;
\draw[thick,red] (1.34,0.55) arc (90:265:0.6) ;
\draw[thick,red,-<] (1.34,0.55) arc (90:200:0.6) ;

\end{tikzpicture}

\caption[A non-closed, non-integrable 1-form.]{An Attempt to show a
  non-closed, non-integrable 1-form in 3-dim (green). If we try and
  close the elements up in to surface we see that we fail. Starting
  the red arc on the lightest green and going round anticlockwise we
  end up on the darkest green, we have helixed down.  }
\label{fig_nonclosed_1-form}
\end{figure}


Non closed $p$-forms consist of a small submanifold of dimension
$(n-p)$ or surface at each point. We call these \defn{form-elements.}
The difference is that in non-closed forms, the form-elements do not
connect up to form a submanifold. Examples are given in figures
\ref{fig_integrate_nonclosed_1form} and \ref{fig_int_2-form_field}.
We can see from figure \ref{fig_integrate_nonclosed_1form} that the
integral now depends on path taken. Thus non-closed forms correspond
to non conservative fields.

\begin{figure}[tb]
\centering
\scalebox{0.70}{
\begin{tikzpicture}
\draw (0,0) node {
\tikz[ rotate=25] {
\filldraw[thick,fill=green,shift={(-2.5,0.5)},scale=1.4,rotate=-35] (0,0) to[out=10,in=190] (4,0)
to[out=40,in=-160] (6,3) to[out=-170,in=10] (2,3) 
to[out=-140,in=40] cycle;
\draw[red,thick] (0,0) 
to[out=90,in=180]  (1,1) to[out=0,in=180] (1.5,.5) 
to[out=0,in=180] (2.5,.7) to[out=0,in=90] (3,0)
to[out=-90,in=0] (2.6,-.6) to[out=180,in=0] (1.6,-.4) 
to[out=180,in=0] (1,-1.1) to[out=180,in=-90] (0,0) ;
\draw[red,thick,scale=1.4,shift={(-0.25,-0)}] (0,0) 
to[out=90,in=180]  (1,1) to[out=0,in=180] (1.5,.5) 
to[out=0,in=180] (2.5,.7) to[out=0,in=90] (3,0)
to[out=-90,in=0] (2.6,-.6) to[out=180,in=0] (1.6,-.4) 
to[out=180,in=0] (1,-1.1) to[out=180,in=-90] (0,0) ;
\draw[red,thick] (0.2,0.) to[out=90,in=180]  (.9,0.6) 
to[out=0,in=90] (1.5,0) to[out=-90,in=0](1,-0.6)  to[out=180,in=-90] cycle ;
\draw[red,thick] (0.7,0) circle (0.3) ;
\draw[red,thick] (2.3,0) circle (0.3) ;
}} ;
\draw (2,1) node {
\tikz[ rotate=25] {
\filldraw[thick,fill=green,shift={(-2.5,0.5)},scale=1.4,rotate=-35] (0,0) to[out=10,in=190] (4,0)
to[out=40,in=-160] (6,3) to[out=-170,in=10] (2,3) 
to[out=-140,in=40] cycle;
\draw[red,thick] (0,0) 
to[out=90,in=180]  (1,1) to[out=0,in=180] (1.5,.5) 
to[out=0,in=180] (2.5,.7) to[out=0,in=90] (3,0)
to[out=-90,in=0] (2.6,-.6) to[out=180,in=0] (1.6,-.4) 
to[out=180,in=0] (1,-1.1) to[out=180,in=-90] (0,0) ;
\draw[red,thick,scale=1.4,shift={(-0.25,-0.)}] (0,0) 
to[out=90,in=180]  (1,1) to[out=0,in=180] (1.5,.5) 
to[out=0,in=180] (2.5,.7) to[out=0,in=90] (3,0)
to[out=-90,in=0] (2.6,-.6) to[out=180,in=0] (1.6,-.4) 
to[out=180,in=0] (1,-1.1) to[out=180,in=-90] (0,0) ;
\draw[red,thick] (0.2,0.) to[out=90,in=180]  (.9,0.6) 
to[out=0,in=90] (1.5,0) to[out=-90,in=0](1,-0.6)  to[out=180,in=-90] cycle ;
\draw[red,thick] (0.7,0) circle (0.3) ;
\draw[red,thick] (2.3,0) circle (0.3) ;
}} ;
\draw[ultra thick,double=white,draw=white] (-3.5,-1) -- +(0,0.7) ;
\draw[ultra thick,draw=blue,->] (-3.5,-1) -- +(0,0.7) ;
\draw[very thick,double=white,draw=white,yscale=0.7] (-2.1,-1.2) arc (0:180:0.3)  ;
\draw[very thick,yscale=0.7,draw=blue,->] (-2.1,-1.2) arc (0:90:0.3)  ;
\draw[very thick,yscale=0.7,draw=blue] (-2.1,-1.2) arc (0:190:0.3)  ;

\draw[very thick,black!60] (5,2) arc (0:270:0.4) ;
\draw[very thick,black!60,->] (5,2) arc (0:90:0.4) ;
\draw[very thick,black!60,->,scale=1.2] (2.9,1.4) arc (0:30:0.4) ;

\end{tikzpicture} 
}
\caption[More non-closed forms.]
{An integrable non-closed 1-form in 3-dim: Two 1-form-submanifolds (green)
  each with a scalar field, represented by contour lines (red). The
  exterior differential of this 1-form are the closed red contour
  lines. We assume that the scalar field is increasing towards the
  inner contours. Both possible twistedness orientations are shown on the same
  diagram.  For untwisted orientations (lower sheet, blue) concatenate the
  inward arrow of the contours with the up arrow of the submanifold
  sheets. The twisted orientation (upper sheet, grey) simply inherits the
  orientation of the sheet.}
\label{fig_NonClosed_Contours}
\end{figure}
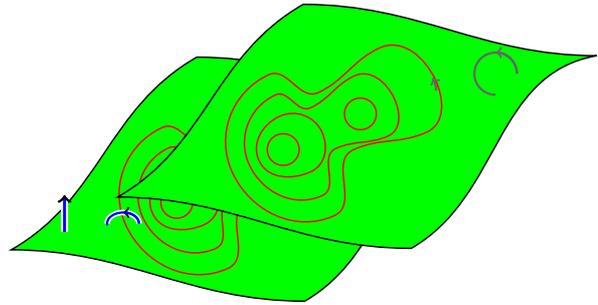

In some (but not all) cases a non-closed form can be thought of as a
collection of closed form-submanifolds, but each of these has a scalar
field on it, see figure \ref{fig_NonClosed_Contours}. Such forms are
known as \defn{integrable}. For example the
1-form (in 2-dim) in figure
\ref{fig_integrate_nonclosed_1form}
is integrable whereas the 1-form (in 3-dim) in figure
\ref{fig_nonclosed_1-form} is not.


\subsection{The exterior differential operator.} 
\label{ch_d}
The exterior
differential operator takes a $p$-form and gives a $(p+1)$-form. One may think of this as taking the
boundary of the form-submanifolds as given by figures \ref{fig_Ori_Bdd-Int} and \ref{fig_Ori_Bdd-Ext}. 

The exterior derivative of a scalar field is simply the
1-form-submanifolds corresponding to the contours, also known as level
$(n-1)$-surfaces, of the scalar field. For untwisted scalar fields the
orientation points in the direction the scalar field 
increases, as in figure \ref{fig_Scalar_contours}.

The exterior differential of an integrable non-closed form are the
contours of the scalar
field on each form-submanifold. See figure \ref{fig_NonClosed_Contours}.

Clearly if the $p$-form is closed then it has no boundary and the
result is zero. In addition, since a boundary has no
boundary, the action of the exterior derivatives operator twice
vanishes.


\begin{figure}[tb]
\def\larw#1{{\tikz[rotate=#1] {
\draw[thick,->] (0,0) to[out=10,in=170] +(1,0) ; 
\draw[thick] (1,0) to[out=-10,in=190] +(1,0) ; 
\fill[blue] (0,0) circle (0.1) node[xshift=-.3cm]{$-$} ;
\fill[red] (2,0) circle (0.1) node[xshift=.3cm]{$+$} ;
}}}
\centering
\begin{tikzpicture}
\filldraw[ultra thick,draw=green!80!black,fill=green!30] 
(0,0) circle(1.8) node {\Large +} ;
\draw[ultra thick,draw=green!20!black,->] (-0.8,1.6) -- (-0.4,0.8) ;
\draw (1,2) node {\larw{40}} ;
\draw (2,1) node {\larw{10}} ;
\draw (2,-1) node {\larw{-40}} ;
\draw (1,-2) node {\larw{-90}} ;
\draw (-2,1) node {\larw{170}} ;
\draw (-2,-1) node {\larw{-140}} ;
\draw (-1,-2) node {\larw{-90}} ;
\end{tikzpicture}
\caption[Demonstration of Stokes's theorem.]
{Demonstration of Stokes's theorem in 2-dimensions. The
  orientation of the twisted disk is $+$ and the twisted circle is
  inwards. The twisted 1-form is not closed and its exterior
  derivative is the positive and negative dots. The integration of
  these dots over the disk is $-7$. Likewise the integration of the
  1-form over the circle is also $-7$.}
\label{fig_integrate_stokes}
\end{figure}
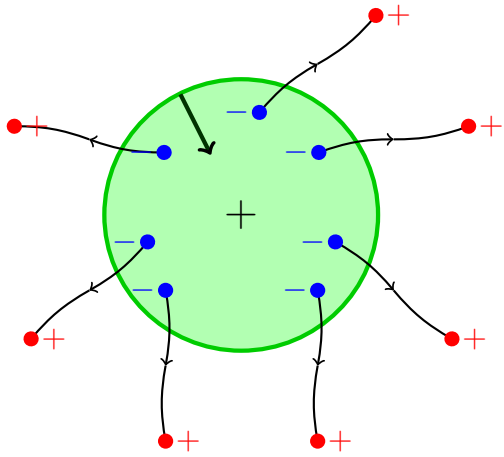


\subsection{Integration and Stokes's theorem for non-closed forms.}
\label{ch_Integration}

The integration of non-closed forms is no longer necessarily zero, even if there
are no holes in the manifold.
Recall to integrate a
$p$-form we must integrate it over a $p$-dimensional submanifold. In
addition untwisted forms are integrated over untwisted submanifolds an
likewise for twisted forms and twisted submanifolds; with plus if the
orientations agree and minus otherwise. 

In figure \ref{fig_integrate_stokes} we see a demonstration of Stokes's
theorem. That is, the integral of a $(p-1)$-form over a
$(p-1)$-dimensional boundary of a submanifold equals the integral of
the exterior derivative of the $(p-1)$-form over the $p$-dimensional
submanifold. Again we have to assume the embedding manifold has no holes.

With a bit of work defining the divergence and curl in terms of
exterior differential operator, one can see that our version of
Stokes's theorem is a generalisation of the divergence and Stokes's
theorem in 3-dimensional vector calculus.

\section{Other operations with scalars, vectors and forms}
\label{ch_other_ops}

\subsection{Combining a 1-form and a vector field to give a scalar field.}  
\label{ch_vec_1form}

\begin{figure}[tb]
\centering
\begin{tikzpicture}[scale=1.1,rotate=-90]
\foreach \x in {1,...,9}
{
  \draw[very thick] (0,\x*0.4)  to[out=0,in=180]  +(2,0.1+0.03*\x) ;
}
\draw[very thick,->] (0.2,2.) -- +(0,-0.2) ;
\draw[ultra thick,->,draw=green!50!black]  (0.4,1.4)  -- +(1.2,2.2) ;
\fill[green!50!black]  (0.4,1.4) circle (0.1)  ;
\end{tikzpicture}
\caption[Combining vectors and 1-forms]
{Combining an untwisted 1-form (black) and an untwisted vector
  field (green) to give a scalar
  field, in 2-dim. The scalar at the base point of the vector would have value $-5$.}
\label{fig_vector_1form}
\end{figure}
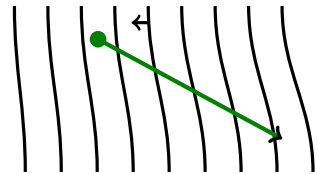


Given a vector field and a 1-form field, one can combine
them together to create a scalar field. The scalar is given by the
number form-submanifolds the vector crosses. See figure
\ref{fig_vector_1form}.
From this we see that if a 1-form contracted with a vector is zero,
then the vector must lie within the 
1-form-submanifold.

\subsection{Vectors acting on a scalar field.}
\label{ch_vec_scalars}
A vector acting on a scalar field gives a new scalar field. This is
given by contracting the exterior derivative of the scalar field with
the vector, as defined above. Thus it corresponds to the number of
times a vector crosses the contour submanifolds of the scalar field.

\subsection{Internal contractions.}
\label{ch_Internal_contraction}

\begin{figure}[tb]
\centering
\begin{tikzpicture}[scale=1.5]
\filldraw[fill=green!50] (0,0) -- (1,.5) -- (1,1.5) -- (0,1) ;
\draw[ultra thick,red] (0,0) -- (0,1) ;
\draw[ultra thick,blue,->] (0,.3) -- +(0.4,0.2) ;
\draw[very thick,double=white,draw=white] (0.2,0.8) arc (0:-180:0.2);
\draw[very thick,->] (0.2,0.8) arc (0:-180:0.2);
\draw[ultra thick,double=white,draw=white] (0.8,1.2) -- +(0.3,-0.2) ;
\draw[ultra thick,->,green!30!black] (0.8,1.2) -- +(0.3,-0.2) ;
\end{tikzpicture}
\caption[Demonstration of internal contraction.]
{In a 3-dim manifold. The internal contraction of the
  untwisted $2$-form (red) with the untwisted vector (blue) to produce 
  the untwisted $1$-form (green). The orientation of the vector
  followed by the orientation of the $2$-form gives the orientation of
the $1$-form.}
\label{fig_internal_contraction}
\end{figure}
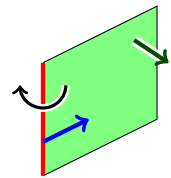

Internal contractions is an operation that takes a vector field and a
$p$-form and gives a $(p-1)$-form. Pictorially it corresponds to the
$(p-1)$-form-submanifold created by extending the $p$-form-submanifold
in the direction of the vector. Thus increases the dimension of the
form-submanifold by 1. See figure \ref{fig_internal_contraction}.
If the vector lies in the $p$-form-submanifold the 
internal contractions vanishes. 

This is a generalisation of the contracting of a vector and a
1-form given in subsection \ref{ch_vec_1form}.
Since the vector lies in the resulting $(p-1)$-form-submanifold,
internally contracting by the same vector twice is zero.


\section{The Metric}
\label{ch_Metric}

\begin{figure}[tb]
\centering
\begin{tikzpicture}[scale=1]
\draw [thick] (0,0) [rotate=40,xscale=2] circle (0.5) ; 
\draw (120:0.75)  node [rotate=30,right] {\footnotesize 0.5} ;
\draw [thick] (0,0) [rotate=40,xscale=2] circle (1) ;
\draw (120:1.25) node [rotate=30,right] {\footnotesize 1} ;
\draw [thick] (0,0) [rotate=40,xscale=2] circle (1.5) ;
\draw (120:1.75) node [rotate=30,right] {\footnotesize 1.5} ;
\draw [thick] (0,0) [rotate=40,xscale=2] circle (2) ;
\draw (120:2.25) node [rotate=30,right] {\footnotesize 2} ;
\draw [ultra thick,blue,->] (0,0) to +(3,0.5) ; 
\draw (0,-3) (0,3) ;
\end{tikzpicture}
\caption[The metric ellipsoid and length.]
{Using the metric ellipsoid to calculate the length of a vector. The
  unit ellipsoid is scaled to various sizes. The vector has length 2.}
\label{fig_metric_ellipse_vector_length}
\end{figure}
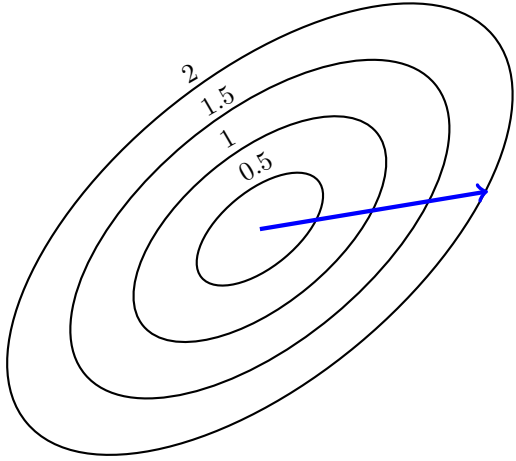

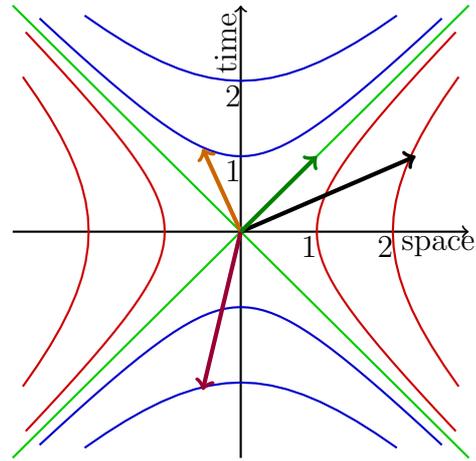
\begin{figure}[tb]
\centering
\begin{tikzpicture}
\draw [thick,->] (0,-3) -- (0,3) ; 
\draw [thick,->] (-3,0) -- (3,0) ; 
\draw (-0.1,0.8) node {1} (-0.1,1.8) node {2} 
(-.20,2.5) node [rotate=90] {time} ;
\draw (0.9,-0.2) node {1} (1.9,-0.2) node {2} (2.6,-.20) node {space} ;
\draw [thick,green!80!black] (-3,3) -- (3,-3) ;
\draw [thick,green!80!black] (3,3) -- (-3,-3) ;
\draw [thick,red!80!black] 
\foreach \a in {1,-1} {
(\a,0) \foreach \x in {0,.1,...,{1.76}} {
 -- ({\a*cosh(\x)},{\a*sinh(\x)}) }
(\a,0) \foreach \x in {0,.1,...,1.76} {
 -- ({\a*cosh(\x)},{-\a*sinh(\x)}) }
(2*\a,0) \foreach \x in {0,.1,...,{1.}} {
 -- ({2*\a*cosh(\x)},{2*\a*sinh(\x)}) }
(2*\a,0) \foreach \x in {0,.1,...,1.} {
 -- ({2*\a*cosh(\x)},{-2*\a*sinh(\x)}) }
} ;
\draw [thick,blue!80!black] 
\foreach \a in {1,-1} {
(0,\a) \foreach \x in {0,.1,...,{1.76}} {
 -- ({\a*sinh(\x)},{\a*cosh(\x)}) }
(0,\a) \foreach \x in {0,.1,...,1.76} {
 -- ({-\a*sinh(\x)},{\a*cosh(\x)}) }
(0,2*\a) \foreach \x in {0,.1,...,{1.}} {
 -- ({2*\a*sinh(\x)},{2*\a*cosh(\x)}) }
(0,2*\a) \foreach \x in {0,.1,...,1.} {
 -- ({-2*\a*sinh(\x)},{2*\a*cosh(\x)}) }
};
\draw [ultra thick,black,->] (0,0) -- (2.3,1) ;
\draw [ultra thick,orange!80!black,->] (0,0) -- (-0.5,1.1) ;
\draw [ultra thick,purple!80!black,->] (0,0) -- (-0.5,-2.1) ;
\draw [ultra thick,green!50!black,->] (0,0) -- (1,1) ;
\end{tikzpicture}
\caption[The metric hyperboloid and length.]{
{
Using the metric hyperboloid to measure the length of
  vectors. 
\\$\bullet$
The diagonal (green) lines are the path of light rays. Thus the
diagonal vector (dark green) is a lightlike vector.
\\$\bullet$
The
  left and right (red) curves are the hyperboloid for measuring the
  length of spatial vectors. (In two dimensions these
  two curves are disjoint. However in three or more dimensions these
  two lines connect to create the one-sheet hyperboloid.) The right-pointing (black) vector is
  spatial and has length 2. 
\\$\bullet$
 The upper and lower (blue) lines are for
  measuring the length of timelike vectors.  The upward
  (orange) line is timelike and has length 1 in the forward
  direction. The downward (purple) line is timelike and has length 2
  in the backward direction. }}
\label{fig_metric_hyper_vector_length}
\end{figure}

The metric gives you the length of vectors (figures
\ref{fig_metric_ellipse_vector_length} and
\ref{fig_metric_hyper_vector_length}) and the angle between two
vectors (figures \ref{fig_metic_elipse_orthogonal} and \ref{fig_metic_hyper_orthogonal}). It can also tell you how long a curve is between to points
and given two points it can tell you which is the shortest such curve,
which is called a geodesic. Thus it can tell you the distance between
two points. This is the length along the geodesic between them.

In higher dimensions the metric tells you the area, volume or
$n$-volume of your manifold and submanifolds within it. It gives you a
measure so that you can integrate scalar fields, figure
\ref{fig_metric_measure}. It enables you to
convert 1-forms into vectors and vectors into 1-forms (called the
metric dual). It also enables you to convert $p$-forms into
$(n-p)$-forms (called the Hodge dual). It doesn't, however, give you
an orientation.

Metrics possess a signature. A list of $n$ positive or negative
signs, where $n$ is the dimension of the manifold. 
There are a number of important dimensions and signatures:
\begin{jgitemize}
\item
One dimensional Riemannian manifolds. 
\item
Two dimensional Riemannian manifolds, i.e. signature $(+,+)$. These are
surfaces, which may be closed like spheres and torii, with a number of
``holes''. They may also be open like planes. They also included non
orientable surfaces such a the Möbius strip.
\item
Three dimensional Riemannian manifolds, i.e. signature $(+,+,+)$. These
include the three dimensional space that we live in.
\item
Four dimensional spacetime, i.e. signature $(-,+,+,+)$. Here the minus
sign refers to the time direction and the three pluses to space.
\end{jgitemize}
However for this document we will also be looking at three dimensional
spacetime, signature $(-,+,+)$, and two dimensional spacetime,
signature $(-,+)$. This is because they are easier to draw.


\begin{figure}[tb]
\centering
\begin{subfigure}{0.22\textwidth}
\begin{tikzpicture}
\shadedraw [ball color=green,shift={(0.0,0.0)},rotate=-20,xscale=0.7] circle (.5) ;
\shadedraw [ball color=green,shift={(1.5,0.0)},rotate=-20,xscale=0.7] circle (.5) ;
\shadedraw [ball color=green,shift={(3.0,0.0)},rotate=-20,xscale=0.7] circle (.5) ;
\shadedraw [ball color=green,shift={(0.0,1.5)},rotate=-20,xscale=0.7] circle (.5) ;
\shadedraw [ball color=green,shift={(1.5,1.5)},rotate=-20,xscale=0.7] circle (.5) ;
\shadedraw [ball color=green,shift={(3.0,1.5)},rotate=-20,xscale=0.7] circle (.5) ;
\shadedraw [ball color=green,shift={(0.0,3.0)},rotate=-20,xscale=0.7] circle (.5) ;
\shadedraw [ball color=green,shift={(1.5,3.0)},rotate=-20,xscale=0.7] circle (.5) ;
\shadedraw [ball color=green,shift={(3.0,3.0)},rotate=-20,xscale=0.7] circle (.5) ;
\end{tikzpicture}
\caption{The unit spheres in flat Riemannian space.}
\label{fig_metric_Remann_flat}
\end{subfigure}
\quad
\begin{subfigure}{0.22\textwidth}
\begin{tikzpicture}
\shadedraw [ball color=green,shift={(0.0,0.0)},rotate=-30,xscale=0.6] circle (.5) ;
\shadedraw [ball color=green,shift={(1.5,0.0)},rotate=-20,xscale=0.7] circle (.5) ;
\shadedraw [ball color=green,shift={(3.0,0.0)},rotate=-10,xscale=0.8] circle (.5) ;
\shadedraw [ball color=green,shift={(0.0,1.5)},xscale=0.7] circle (.5) ;
\shadedraw [ball color=green,shift={(1.5,1.5)},xscale=0.8] circle (.5) ;
\shadedraw [ball color=green,shift={(3.0,1.5)},xscale=1]   circle (.5) ;
\shadedraw [ball color=green,shift={(0.0,3.0)},rotate=30,xscale=0.6] circle (.5) ;
\shadedraw [ball color=green,shift={(1.5,3.0)},rotate=20,xscale=0.7] circle (.5) ;
\shadedraw [ball color=green,shift={(3.0,3.0)},rotate=10,xscale=0.8] circle (.5) ;
\end{tikzpicture}
\caption{The unit spheres in Riemannian space.}
\label{fig_metric_Remann_spheres}
\end{subfigure}
\caption[Representations of the metric in 3-dim space.]
        {Representations of the metric in 3-dim space.}
\label{fig_Metric}
\end{figure}
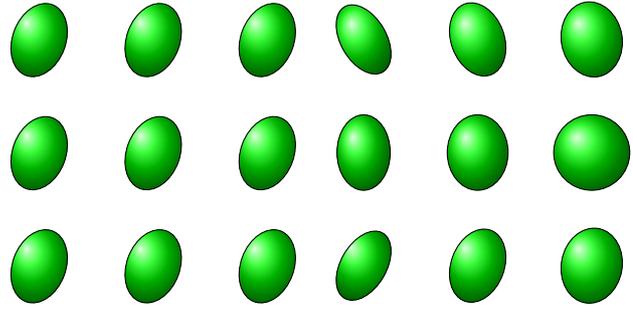


\begin{figure}[tb]
\begin{subfigure}[t]{0.13\textwidth}
\begin{tikzpicture}
[scale=0.4]
\shadedraw[ball color=blue,thick] (0,0) 
 .. controls (2.5,-2.5) and (3.5,-2.5) .. +(6,0) ; 
\shadedraw[ball color=blue,thick] (0,-6) 
 .. controls (2.5,-3.5) and (3.5,-3.5) .. (6,-6) arc (-60:-120:6) ; 
\shadedraw[left color=black, right color=red ,thick,yscale=0.3]
(3.0,0.9) circle (3.1) ; 
\draw [thick,->] (0,-3) -- +(0,2.5) ;
\draw (-0.5,-2) node[rotate=90] {time} ;
\draw [thick,->] (0,-3) -- +(2.6,.4) ;
\draw (1.5,-2.4) node[rotate=10] {space} ;
\draw [thick,->] (0,-3) -- +(2.6,-.4) ;
\draw (1.5,-3.6) node[rotate=-10] {space} ;
\end{tikzpicture}
\caption{
The two-sheet (upper and lower) hyperboloid.}
\label{fig_hyperbolid_twosheet}
\end{subfigure}
\quad
\begin{subfigure}[t]{0.13\textwidth}
\begin{tikzpicture}
[scale=0.4]
\draw [thick,->] (0,-3) -- +(2.6,.4) ;
\draw (1.5,-2.4) node[rotate=10] {space} ;
\shadedraw[ball color=yellow!50!black,thick] (0,0) 
 --(3,-3) -- (0,-6) 
arc (-120:-60:6) -- (3,-3) -- (6,0) ;
\shadedraw[left color=black, right color=red ,thick,yscale=0.3]
(3.0,0.9) circle (3.1) ; 
\draw [thick,->] (0,-3) -- +(0,2.5) ;
\draw (-0.5,-2) node[rotate=90] {time} ;
\draw [thick,->] (0,-3) -- +(2.6,-.4) ;
\draw (1.5,-3.6) node[rotate=-10] {space} ;
\end{tikzpicture}
\caption{The double cone. }
\label{fig_hyperbolid_doublecone}
\end{subfigure}
\quad
\begin{subfigure}[t]{0.13\textwidth}
\begin{tikzpicture}
[scale=0.4]
\draw [thick,->] (0,-3) -- +(2.6,.4) ;
\draw (1.5,-2.4) node[rotate=10] {space} ;
\shadedraw[ball color=green,thick] (0,0) 
 .. controls (2.5,-2.5) and (2.5,-3.5) .. +(0,-6) 
arc (-120:-60:6) .. controls (3.5,-3.5) and (3.5,-2.5) .. (6,0) ;
\shadedraw[left color=black, right color=red ,thick,yscale=0.3]
(3.0,0.9) circle (3.1) ; 
\draw [thick,->] (0,-3) -- +(0,2.5) ;
\draw (-0.5,-2) node[rotate=90] {time} ;
\draw [thick,->] (0,-3) -- +(2.6,-.4) ;
\draw (1.5,-3.6) node[rotate=-10] {space} ;
\end{tikzpicture}
\caption{The one-sheet hyperboloid.}
\label{fig_hyperbolid_onesheet}
\end{subfigure}
\caption[The spacetime metric.]{The spacetime metric (in 3
  dimensions). The two-sheet hyperboloid
  (\ref{sub@fig_hyperbolid_twosheet}) is used for measuring timelike
  vectors. The upper hyperboloid is for future pointing vectors, and
  the lower for past pointing vectors.  The one-sheet hyperboloid
  (\ref{sub@fig_hyperbolid_onesheet}) is used for measuring the length
  of spacelike vectors. Between the two hyperboloids is the double
  cone, which all lightlike vectors lie on.
} 
\label{fig_hyperbolids}
\end{figure}
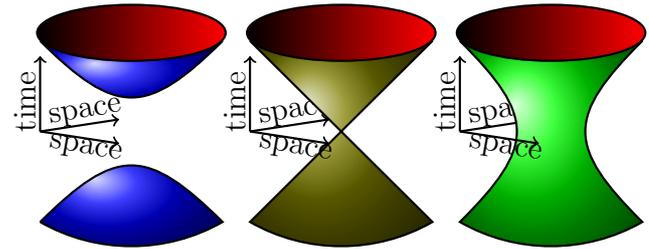


\subsection{Representing the metric: Unit ellipsoids and hyperboloids} 
\label{ch_Metric_unit_ellipse}

Here we represent the metric is by the unit ellipsoid or
hyperboloid. See also \cite{misner1957classical}. Around a point
consider all the vectors which have unit length. For a Riemannian
manifold, these vectors form an ellipse in 2 dimensions and an
ellipsoid in 3 dimensions. In figure \ref{fig_Metric} we show the
ellipsoid for a Riemannian manifold. If the ellipsoids are the same at
each point then the metric is flat as in figure
\ref{fig_metric_Remann_flat}. By contrast for a curved metric the
ellipsoids must be different at each point as in figure
\ref{fig_metric_Remann_spheres}.

For spacetime metrics the ellipsoids is replaced by a hyperboloid, as
in figure \ref{fig_hyperbolids}. Since there is a minus in the
signature, vectors are divided into three types:
\begin{jgitemize}
\item
Timelike vectors, these represent the 4-velocity of massive particles. These
vectors point to one of the components of the two-sheet hyperboloid,
as seen in figure \ref{fig_hyperbolid_twosheet}. As
the name suggest, the two-sheet hyperboloid has two components, the
upper and lower. This distinguishes between future point and past
pointing vectors. 
\item
Lightlike vectors, the 4-velocity of light. These vector lie on the double
cone. See figure \ref{fig_hyperbolid_doublecone}. Again one can
distinguish between past and future pointing lightlike vectors.
\item
Spacelike vectors, These are  not the 4-velocity of anything. These
vectors point to the one-sheet hyperboloid. One cannot distinguish
between past and future spacelike vectors.
\end{jgitemize}


\subsection{Length of a vector}

Using the metric ellipsoid or hyperboloid to calculate the length of a vector is
easy. For a Riemannian metric the ellipsoid represent all vectors of
unit length. Therefore simply compare length of the vector with 
the radius of ellipsoid in the direction of the vector as in
figure \ref{fig_metric_ellipse_vector_length}.

We see in figure \ref{fig_metric_hyper_vector_length}, we use the
two-sheet (upper and lower) hyperboloid to measure the length of
forward and backward pointing timelike vectors. We also use the
one-sheet hyperboloid (left and right curves) to measure the spacelike
vectors.

\subsection{Orthogonal subspace}
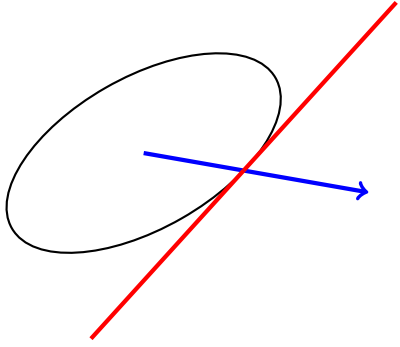
\begin{figure}[tb]
\centering
\begin{tikzpicture}[scale=1,rotate=-10]
\draw [thick] (0,0) [rotate=40,xscale=2] circle (1) ;
\draw [ultra thick,blue,->] (0,0) to +(3,0) ; 
\draw [ultra thick,red] (1.335,0) -- +(58:3) ;
\draw [ultra thick,red] (1.335,0) -- +(58:-3) ;
\end{tikzpicture}
\caption[The metric ellipsoid and orthogonality.]
{Using the metric ellipsoid to prescribe the orthogonal space (red) to a
  vector (blue). The orthogonal space is tangential to the ellipsoid.}
\label{fig_metic_elipse_orthogonal}
\end{figure}

Given a vector then the metric can be used to specify all the directions
which are orthogonal to that vector. Orthogonal vectors are also
known as perpendicular vectors or vectors at right angles to each
other. This is given by the line tangent to the ellipsoid (figure \ref{fig_metic_elipse_orthogonal}) or
hyperboloid (figure \ref{fig_metic_hyper_orthogonal}) at the point
where the vector crosses it. In 3 dimensions the set of orthogonal
directions form a plane and in 4 dimensions they would form a 3-volume.

\subsection{Angles and $\gamma$-factors}

\begin{figure}[tb]
\centering
\begin{tikzpicture}[scale=1,rotate=-5]
\draw [thick] (0,0) [rotate=40,xscale=2] circle (1) ;
\draw [ultra thick,black,->] (0,0) to +(3,0) ; 
\draw [ultra thick,blue,->] (0,0) to +(34:2.5) ;
\draw [very thick,red] (1.335,0) -- +(58:3) ;
\draw [very thick,red] (0.94,0) -- +(58:3) ;
\draw [thick,<->]  (0,-.1) +(0.5,-.2) node {$\ell$} (0,-.1)  -- +(0.94,0) ;
\end{tikzpicture}
\caption[The metric ellipsoid and angles.]{Using the metric
  ellipsoid to measure the angle between two vectors. The (red) lines
  are orthogonal to one of the vectors (black). The length $\ell$ is
  the cosine of the angle between the two vectors.}
\label{fig_metic_elipse_angle}
\end{figure}
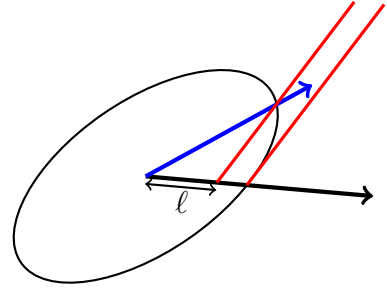

Given two vectors then a Riemannian metric can be used to find the
cosine of the angle between them as shown in figure
\ref{fig_metic_elipse_angle}. This cosine will always be a value
between $-1$ and $1$. 

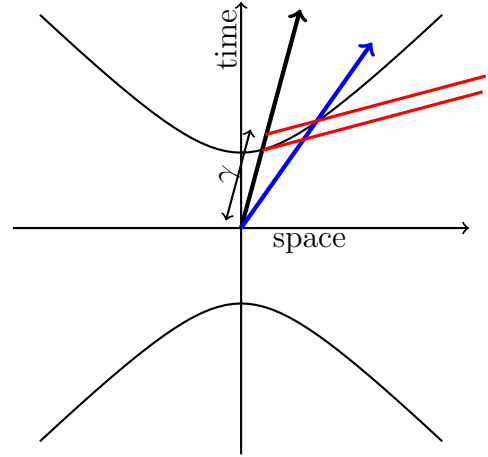
\begin{figure}[tb]
\centering
\begin{tikzpicture}
\draw [thick,->] (0,-3) -- (0,3) ; 
\draw [thick,->] (-3,0) -- (3,0) ; 
\draw (-.20,2.5) node [rotate=90] {time} ;
\draw (0.9,-0.2) node {space} ;
\draw [thick] 
\foreach \a in {1,-1} {
(0,\a) \foreach \x in {0,.1,...,{1.76}} {
 -- ({\a*sinh(\x)},{\a*cosh(\x)}) }
(0,\a) \foreach \x in {0,.1,...,1.76} {
 -- ({-\a*sinh(\x)},{\a*cosh(\x)}) }
};
\draw [ultra thick,->] (0,0) -- (75:3) ;
\draw [ultra thick,blue,->] (0,0) -- +(55:3) ;
\draw [very thick,red] (.32,1.24) -- +(15:3) ;
\draw [very thick,red] (.28,1.03) -- +(15:3) ;
\draw [thick,<->] (-.2,.1) -- +(75:1.25) ;
\draw (-.15,.7) node[rotate=75] {$\gamma$} ;
\end{tikzpicture}
\caption[The metric hyperboloid and $\gamma$-factor.]{Using the metric
  hyperboloid to measure the angle between two vectors. The (red)
  lines are orthogonal to one of the vectors (black). The length
  $\gamma$ is the $\gamma$-factor required to boost one of the vectors
  to the other.}
\label{fig_metic_hyper_angle}
\end{figure}

In the case of a spacetime metric, then in figure
\ref{fig_metic_hyper_angle} we get a value which is greater than 1. We
call this value the $\gamma$-factor. It is one of the key values in
special relativity. Assuming that the two vectors are the 4-velocity
of two particles, then the $\gamma$-factor is a  
function of the relative velocity between the two particles.


\begin{figure}[tb]
\centering
\begin{tikzpicture}
\draw [thick,->] (0,-3) -- (0,3) ; 
\draw [thick,->] (-3,0) -- (3,0) ; 
\draw (-.20,2.5) node [rotate=90] {time} ;
\draw (0.9,-0.2) node {space} ;
\draw [thick,black] 
\foreach \a in {1,-1} {
(0,\a) \foreach \x in {0,.1,...,{1.76}} {
 -- ({\a*sinh(\x)},{\a*cosh(\x)}) }
(0,\a) \foreach \x in {0,.1,...,1.76} {
 -- ({-\a*sinh(\x)},{\a*cosh(\x)}) }
};
\draw [ultra thick,->] (0,0) -- (75:3) ;
\draw [ultra thick,red] (-3,.16) -- +(15:6) ;
\end{tikzpicture}
\caption[The metric hyperboloid and orthogonality.]{Using the metric
  hyperboloid to prescribe the orthogonal space (red) to a vector
  (black). The orthogonal space is tangential to the hyperboloid.}
\label{fig_metic_hyper_orthogonal}
\end{figure}
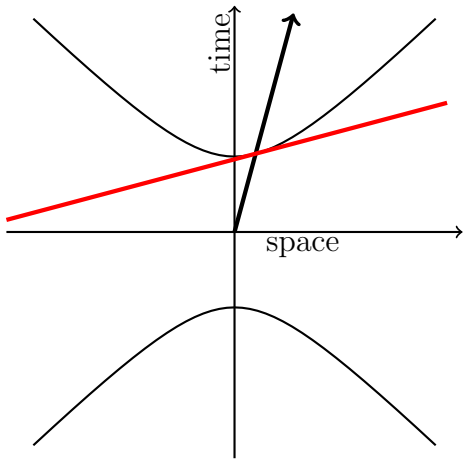


\subsection{Magnitude of a form.}
\label{ch_Metric_Mag_forms}

\begin{figure}[tb]
\centering
\begin{tikzpicture}[scale=1.1,rotate=0]
\foreach \x in {1,...,9}
{
  \draw[very thick] (0,\x*0.4)  to[out=0,in=180]  +(2,0.1+0.03*\x) ;
}
\draw[very thick,->] (0.2,2.) -- +(0,-0.2) ;
\draw[ultra thick,draw=green!50!black,rotate=10]  (1.4,1.8) ellipse (.7 and 1.4) ;
\end{tikzpicture}
\caption[The metric, magnitude of a 1-form]{Using the metric to measure the
  magnitude of a 1-form in 2 dimensions. Simply count the number of
  form-submanifold inside the ellipse and divide by 2. In this case the 1-form
  has magnitude $3.5$}
\label{fig_metric_magnitude_form}
\end{figure}
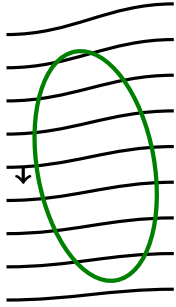


One can use the metric to measure the magnitude of forms. On an $n-$dimensional
Riemannian manifold the magnitude to a 1-form is given by counting the
$(n-1)$-dimensional form-submanifolds inside the ellipsoid and dividing by 2. In
figure \ref{fig_metric_magnitude_form} we see how to use the metric ellipse to
measure the magnitude to a 1-form in 2 dimensions. Measuring the magnitude of a
$2-$form, requires counting the number of $(n-2)$-dimensional form-submanifolds
inside the ellipsoid and dividing by $\pi$. Similar formula exist for higher
degrees and for spacetime manifolds.


\subsection{Metric dual.}
\label{ch_Metric_dual}

\begin{figure}[tb]
\centering
\begin{subfigure}{0.45\textwidth}
\centering
\begin{tikzpicture}[scale=0.85]
\fill[red!70,rotate=-10] (-2,0) -- (0,1) -- (2,0) -- (0,-1) -- cycle ;
\draw[ultra thick,->,color=blue] (0,0) -- (0,2) ;
\draw[ultra thick,draw=white,double=red!80!black] (1,0) -- (1,1) ; 
\draw[ultra thick,draw=red!80!black,->] (1,0) -- (1,1) ; 
\draw[scale=0.8] (0,1) to [rotate=-60] (0,1) to [ rotate=-120] (0,1) ;
\draw[yscale=0.6,xscale=-0.6] (0,1) to [rotate=-60] (0,1) to [ rotate=-120] (0,1) ;
\end{tikzpicture}
\caption{Spacelike vector (blue). The orientation
  (red) of the 1-form is unchanged.}
\end{subfigure}
\quad
\begin{subfigure}{0.45\textwidth}
\centering
\begin{tikzpicture}[scale=0.85]
\draw[ultra thick,draw=red!80!black,->] (1,0) -- (1,-1) ; 
\filldraw[draw=white,fill=red!70,rotate=-10] (-2,0) -- (0,1) -- (2,0) -- (0,-1) -- cycle ;
\draw[ultra thick,->,color=blue] (0,0) -- (0,2) ;
\draw[scale=0.8] (0,1) to [rotate=-60] (0,1) to [ rotate=-120] (0,1) ;
\draw[yscale=0.6,xscale=-0.6] (0,1) to [rotate=-60] (0,1) to [ rotate=-120] (0,1) ;
\end{tikzpicture}
\caption{Timelike vector (blue). The orientation
  (red) of the 1-form is reversed.}
\end{subfigure}
\quad
\begin{subfigure}{0.45\textwidth}
\centering
\begin{tikzpicture}[scale=0.85]
\filldraw[draw=white,fill=red!70,rotate=-10] (-2,0) -- (0,1) -- (2,0) -- (0,-1) -- cycle ;
\draw[ultra thick,draw=white,double=red!80!black] (-1,0) -- +(0,1) ; 
\draw[ultra thick,draw=red!80!black,->] (-1,0) -- +(0,1) ; 
\draw[thick,double=white,draw=white] (0,0) -- (1.9,0.3) ;
\draw[thick,->,color=green!80!black] (0,0) -- (1.9,0.3) ;
\draw[ultra thick,->,color=blue] (0,0) -- (2,0) ;
\end{tikzpicture}
\caption{Lightlike vector (blue). The orientation
  (red) of the 1-form points in the direction of the nearby spacelike
  vectors (green).}
\end{subfigure}

\caption[The metric dual.]{The metric dual converts
  a vector into a 1-form and visa versa. The untwisted vector (blue)
  is converted into the untwisted 1-form (red).}
\label{fig_Metric_dual}
\end{figure}
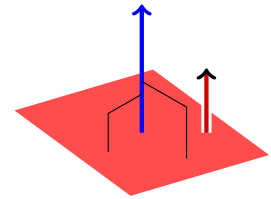
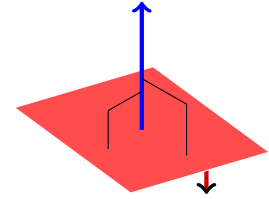
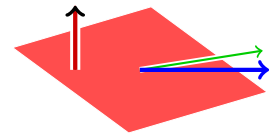


The metric can be used to convert a vector into a 1-form and a 1-form
into a vector. Given a vector then the 1-form is the
$(n-1)$-form-submanifold which is perpendicular to the vector, as
given in figures \ref{fig_metic_elipse_orthogonal} and
\ref{fig_metic_hyper_orthogonal}. See also figure \ref{fig_Metric_dual}.

The twistedness of the orientation is unchanged.
For spacelike vectors the orientation is preserved, whereas for
timelike vectors the orientation is reversed. For lightlike vectors,
the vector is orthogonal to itself. Therefore it lies in its own
orthogonal compliment. For these the orientation points in the
direction of the nearby spacelike vectors.


\subsection{Hodge dual.}

\begin{figure}[tb]
\centering
\begin{subfigure}{0.45\textwidth}
\centering
\begin{tikzpicture}[scale=1.3,rotate=10]
\draw[ultra thick,red] (0,-1) to[out=105,in=-75] (0,1) ;
\draw[thick,->] (0.03,0.7) -- +(.3,0) ;
\draw[ultra thick,blue] (-1,0) to[out=10,in=-170] (1,0) ;
\draw[ultra thick,->,blue] (0.5,-.03) to[out=10,in=-170] +(0.1,0) ;
\draw (-0.3,0) -- (-0.3,0.3) -- (0.05,0.3) ;
\end{tikzpicture}
\caption{In 2-dim. Hodge dual of the untwisted 1-form (red), becomes the twisted
  1-form (blue).}
\end{subfigure}
\begin{subfigure}{0.45\textwidth}
\centering
\begin{tikzpicture}[scale=0.85,rotate=-70]
\draw[ultra thick,color=blue] (0,-1.7) to[out=95,in=-85] (0,0) ;
\fill[red!70,rotate=-10] (-2,0) to[out=10,in=190] (0,1) -- (2,0)
   to[out=190,in=10](0,-1) -- cycle ;
\draw[ultra thick,->,color=blue] (0,0) to[out=95,in=-85] (0,1.5) ;
\draw[ultra thick,color=blue] (0,1.5) to[out=95,in=-90] (0,2) ;
\draw[ultra thick,draw=white,double=red!80!black] (1,0) -- (1,1) ; 
\draw[ultra thick,draw=red!80!black,->] (1,0) -- (1,1) ; 
\draw[scale=0.8] (0,1) to [rotate=-60] (0,1) to [ rotate=-120] (0,1) ;
\draw[yscale=0.6,xscale=-0.6] (0,1) to [rotate=-60] (0,1) to [ rotate=-120] (0,1) ;
\end{tikzpicture}
\caption{In 3-dim. Hodge dual of the untwisted 1-form (red), becomes the twisted
  2-form (blue).}
\end{subfigure}
\quad
\begin{subfigure}{0.45\textwidth}
\centering
\begin{tikzpicture}[scale=0.85,rotate=-70]
\draw[ultra thick,color=red] (0,-1.7) to[out=95,in=-85] (0,0) ;
\fill[blue!60!white,rotate=-10] (-2,0) to[out=10,in=190] (0,1) -- (2,0)
   to[out=190,in=10](0,-1) -- cycle ;
\draw[ultra thick,color=red!90!black] (0,0) to[out=95,in=-90] (0,2) ;
\draw[scale=0.8] (0,1) to [rotate=-60] (0,1) to [ rotate=-120] (0,1) ;
\draw[yscale=0.6,xscale=-0.6] (0,1) to [rotate=-60] (0,1) to [ rotate=-120] (0,1) ;
\draw[ultra thick,color=red!90!black] (0,0) to[out=95,in=-90] (0,2) ;
\draw[thick,yshift=1cm,yscale=-0.5,draw=white,double=green!50!black] 
(-0.09,-1.2) arc (260:-80:0.4) ; 
\draw[thick,yshift=1cm,draw=green!50!black] 
(0.27,0.4) -- +(0.1,0.1) -- +(0.2,0) ; 
\draw[thick,yshift=1cm,yscale=-0.5,draw=black,->] 
(-1.5, 1.5) arc (200:-40:0.4) ; 
\end{tikzpicture}
\caption{In 3-dim. Hodge dual of the untwisted 2-form (red), becomes the twisted
  1-form (blue).}
\end{subfigure}

\caption[The Hodge dual.]{Hodge dual of untwisted forms (red) to
  twisted forms (blue) in 2 and 3-dimensional Riemannian geometry. For
  spacetime geometry the resulting orientation may be reversed.}
\label{fig_Hodge_dual}
\end{figure}
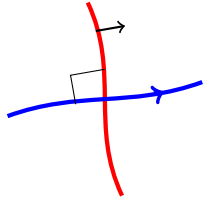
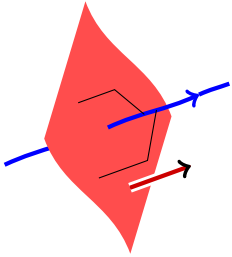
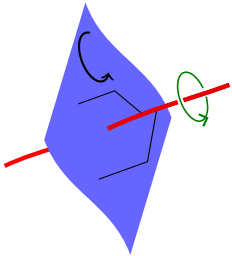

The Hodge dual takes a $p$-form and gives the unique $(n-p)$-form which is
orthogonal to it. 

It changes the twistedness of the form. In 
Riemannian geometry, if one
starts with an untwisted $p$-form it gives a twisted $(n-p)$-form with
the same orientation. If one starts with a twisted $p$-form it gives
the untwisted form with either the same or the opposite
orientation. This is given by double twisting as follows:
\begin{jgitemize}
\item
If $p$ is even or $(n-p)$ is even then the orientation is unchanged.
\item
If both $p$ and $(n-p)$ are odd then the orientation is reversed.
\end{jgitemize}
For non Riemannian geometry the formula for the resulting orientation
is more complicated.


\subsection{The measure from a metric}
\label{ch_Metric_measure}

\begin{figure}[tb]
\centering
\begin{tikzpicture}
\draw[ultra thick,green!50!black] (0,0) ellipse (0.75 and 0.5)  
node[blue] {$\boldsymbol\oplus$};
\draw[ultra thick,green!50!black,rotate=10] (0,1) ellipse (0.85 and 0.5) 
node[blue] {$\boldsymbol\oplus$};
\draw[ultra thick,green!50!black,rotate=20] (0,2) ellipse (0.95 and 0.5) 
node[blue] {$\boldsymbol\oplus$};
\draw[ultra thick,green!50!black,rotate=10] (1.4,0) ellipse (0.65 and 0.5) 
node[blue] {$\boldsymbol\oplus$};
\draw[ultra thick,green!50!black,rotate=20] (1.8,0.75) ellipse (0.75 and 0.5) 
node[blue] {$\boldsymbol\oplus$};
\draw[ultra thick,green!50!black,rotate=30] (2.0,1.45) ellipse (0.85 and 0.5) 
node[blue] {$\boldsymbol\oplus$};
\draw[ultra thick,green!50!black,rotate=20] (2.7,-0.5) ellipse (0.65 and 0.45) 
node[blue] {$\boldsymbol\oplus$};
\draw[ultra thick,green!50!black,rotate=30] (3.2,-0.15) ellipse (0.75 and 0.45) 
node[blue] {$\boldsymbol\oplus$};
\draw[ultra thick,green!50!black,rotate=40] (3.5,0.2) ellipse (0.85 and 0.45) 
node[blue] {$\boldsymbol\oplus$};
\end{tikzpicture}
\caption[The measure from the metric.]{The measure from the
  metric. In 2 dim, the twisted top form is given by packing together the metric
  ellipses.}
\label{fig_metric_measure} 
\end{figure}
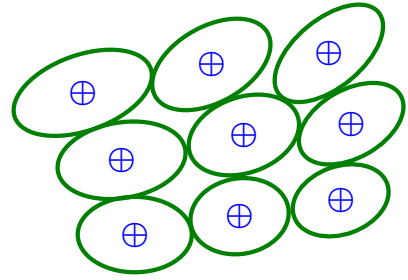

A metric naturally gives rise to a measure. That is a twisted
top-form. This is actually the Hodge dual of the untwisted constant scalar
field $1$. It is given by packing together the metric ellipsoids. See figure
\ref{fig_metric_measure}. We can
integrate this measure to give the size of any manifold even of 
it is not orientable.

\section{Smooth maps}
\label{ch_Maps}
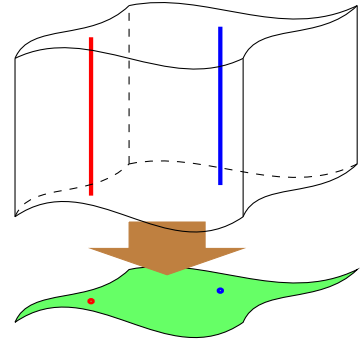
\begin{figure}[tb]
\centering
\begin{tikzpicture}[yscale=.7]
\filldraw[yshift=-2cm,fill=green!60!white] 
 (0,0) .. controls (1,1) and (2,-1) .. (3,0) 
       .. controls (3.2,0.7) and (4.,0.3) .. (4.5,1) 
       .. controls (3.5,0.3) and (2.5,1.3)  .. (1.5,1)  
       .. controls (1.,0.3) and (0.2,0.7)  .. (0,0) ;

\draw[blue,ultra thick] (2.7,0.6) -- +(0,3.) ;
\draw[blue,ultra thick,yshift=-2cm] (2.7,0.6) circle (0.02);
\draw[red,ultra thick] (1,0.4) -- +(0,3.) ;
\draw[red,ultra thick,yshift=-2cm] (1,0.4) circle (0.02) ;

\filldraw[xshift=1.5cm,color=brown,yshift=-0.1cm,yscale=-0.5] 
(0,0) -- (1,0) -- (1,1) -- (1.5,1) -- (0.5,2)  --
  (-0.5,1) -- (0,1) -- (0,0) ;

\draw (0,0) .. controls (1,1) and (2,-1) .. (3,0) ;
\draw (0,3) .. controls (1,3.7) and (2,2.3) .. (3,3) ;
\draw (0,0) -- (0,3) ;
\draw (3,0) -- (3,3) ;
\draw (3,3) .. controls (3.2,3.7) and (4.,3.3) .. (4.5,4) ;
\draw (0,3) .. controls (0.2,3.7) and (1.,3.3) .. (1.5,4) ;
\draw (3,0) .. controls (3.2,0.7) and (4.,0.3) .. (4.5,1) ;
\draw (1.5,4) .. controls (2.5,4.3) and (3.5,3.3) .. (4.5,4) ;
\draw (4.5,1) -- (4.5,4) +(.1,0) ;
\draw [dashed] (0,0) .. controls (0.2,0.7) and (1.,0.3) .. (1.5,1) ;
\draw [dashed] (1.5,1) .. controls (2.5,1.3) and (3.5,0.3) .. (4.5,1) ;
\draw [dashed] (1.5,1) -- (1.5,4) ;
\end{tikzpicture}
\caption[Projection 3-dim to 2-dim.]
{A projection mapping a 3-dim manifold onto a 2-dim manifold
  (green). All the points on each line (red and blue) 
  are projected onto the respective dots.}
\label{fig_proj32}
\end{figure}


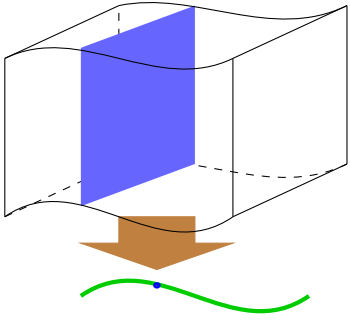
\begin{figure}[tb]
\begin{tikzpicture}[yscale=.7]
\centering
\draw[ultra thick,yshift=-1.5cm,xshift=1cm,draw=green!80!black] 
 (0,0) .. controls (1,1) and (2,-1) .. (3,0) ;

\filldraw[xshift=1.5cm,color=brown,yshift=0cm,yscale=-0.5] 
(0,0) -- (1,0) -- (1,1) -- (1.5,1) -- (0.5,2)  --
  (-0.5,1) -- (0,1) -- (0,0) ;

\draw [dashed] (0,0) -- (1.5,1) ;
\draw [dashed] (1.5,1) .. controls (2.5,1.3) and (3.5,0.3) .. (4.5,1) ;
\draw [dashed] (1.5,1) -- (1.5,4) ;

\fill[fill=blue!60!white,xshift=1cm] (0,0.2) -- (0,3.2) -- (1.5,4) -- (1.5,1) -- cycle ; 

\draw (0,0) .. controls (1,1) and (2,-1) .. (3,0) ;
\draw (0,3) .. controls (1,3.7) and (2,2.3) .. (3,3) ;
\draw (0,0) -- (0,3) ;
\draw (3,0) -- (3,3) ;
\draw (3,3) -- (4.5,4) ;
\draw (0,3) -- (1.5,4) ;
\draw (3,0) -- (4.5,1) ;
\draw (1.5,4) .. controls (2.5,4.3) and (3.5,3.3) .. (4.5,4) ;
\draw (4.5,1) -- (4.5,4) ;

\draw[blue,ultra thick] (2,-1.3) circle (0.02) ;

\end{tikzpicture}
\caption[Projection 3-dim to 1-dim.]
{A projection mapping a 3-dim manifold onto a 1-dim manifold
  (green). All the points on the plane (blue) are projected onto the dot.}
\label{fig_proj31}
\end{figure}


A smooth map is a function that maps one manifold into another. 

The simplest case of a smooth map is a diffeomorphism where there is a
one to one correspondence between points in the source manifold and
points in the target manifold. As a result the source and target
manifolds must have the same dimension and in some respects they look
the same, that is all the concepts in sections
\ref{ch_DG}-\ref{ch_other_ops}. By contrast when a metric is included
the two manifolds can look very different.  One example is in
transformation optics where one maps between one manifold where
light rays travel in straight lines and another manifold where
light rays travel around the hidden object.

One class of smooth maps that we have encountered is that of the embedding of
submanifolds as seen in figure \ref{fig_Man_submanifolds}. In this
case we map a lower dimensional (embedded) manifold into a higher
dimensional (embedding)
manifold. In this case every point in the embedded submanifold corresponds
to a unique point in the embedding manifold. However some points in the
embedding manifold do not correspond to a point in the embedded manifold.

Another class are regular projections. These map a higher dimensional manifold
into a lower dimensional manifold. See figures \ref{fig_proj32} and
\ref{fig_proj31}. 
In this case every point in the image manifold corresponds to a
submanifold of points in the domain.

Other maps however are none of the above and may include self
intersections (figure \ref{fig_Man_patho_submanifolds}), folding and a
variety of other  more complicated situations.


\subsection{Pullbacks of $p$-forms.}
\label{ch_PullBacks}

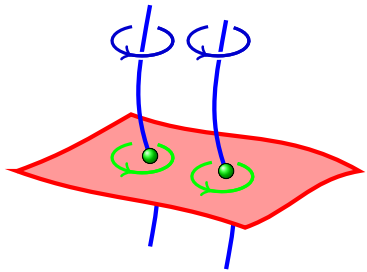
\begin{figure}[tb]
\centering
\begin{tikzpicture}[xscale=-1,rotate=90]
\draw[ultra thick,draw=blue] (-1,2.5) to[out=10,in=160] (0.2,2.5) ;
\draw[ultra thick,draw=blue] (-1.3,3.5) to[out=10,in=160] (0.,3.5) ;
\filldraw[ultra thick,draw=red,fill=red!40!white,scale=1.5] (0,0.5) 
to [out=110,in=-70] +(-.5,2) 
to [out=70,in=-100] +(0.5,1) 
to [out=-60,in=110] +(0.5,-2) 
to [out=-120,in=70] cycle
;
\draw[ultra thick,draw=blue] (0.2,2.5) to[out=-20,in=190] (2.2,2.5) ;
\draw[ultra thick,draw=blue] (0.,3.5) to[out=-20,in=190]  (2.,3.5) ;

\shadedraw[ball color=green] (0.2,2.5) circle (0.1) ;
\shadedraw[ball color=green] (0.,3.5) circle (0.1) ;

\draw[very thick, draw=white,double]  (1.9,2.6) 
   [xscale=0.5] arc (30:340:.4);
\draw[very thick,blue!80!black]  (1.9,2.6) 
   [xscale=0.5] arc (30:340:.4);
\draw[very thick,blue!80!black]  (1.9,2.6) 
   [xscale=0.5,-<] arc (30:220:.4);

\draw[very thick, draw=white,double]  (1.9,3.6) 
   [xscale=0.5] arc (30:340:.4);
\draw[very thick,blue!80!black]  (1.9,3.6) 
   [xscale=0.5] arc (30:340:.4);
\draw[very thick,blue!80!black]  (1.9,3.6) 
   [xscale=0.5,-<] arc (30:220:.4);

\draw[very thick, draw=green]  (.35,2.6) 
   [xscale=0.5] arc (30:340:.4);
\draw[very thick,green]  (0.35,2.6) 
   [xscale=0.5,-<] arc (30:220:.4);

\draw[very thick, draw=green]  (.1,3.65) 
   [xscale=0.5] arc (30:340:.4);
\draw[very thick,green]  (0.1,3.65) 
   [xscale=0.5,-<] arc (30:220:.4);

\end{tikzpicture}
\caption[Pullback, 2-dim to 3-dim, 2-form]
{The pullback with respect to a 2-dim submanifold (red plane)
  embedded in a 3-dim manifold. 
  The pullback of an untwisted 2-form (blue curves) gives an untwisted
  2-form (green dot).}
\label{fig_pull23_untw2}
\end{figure}


\begin{figure}[tb]
\centering
\begin{tikzpicture}[xscale=1.2,yscale=-1.2]

\filldraw[ultra thick,draw=red,fill=red!40!white]
(-1.2,3) -- (-.7,1.2)  -- (0.3,1.2)  to [out=120,in=-80] (-.2,3) -- cycle;

\filldraw[ultra thick,draw=blue,fill=blue!40!white] (0,0.5) 
to [out=110,in=-70] +(-.5,2) 
to [out=70,in=-100] +(0.5,1) 
to [out=-60,in=110] +(0.5,-2) 
to [out=-120,in=70] cycle
;
\filldraw[ultra thick,draw=red,fill=red!40!white]
(0.3,1.2)  to [out=110,in=-70] (-.2,3) -- (1.1,3)  -- (1.6,1.1)--  cycle;

\draw[ultra thick,draw=green]
(0.3,1.2)  to [out=110,in=-70] (-.2,3) ;

\filldraw[ultra thick,draw=blue,fill=blue!40!white,xshift=1.33cm] (0,0.5) 
to [out=110,in=-70] +(-.5,2) 
to [out=70,in=-100] +(0.5,1) 
to [out=-60,in=110] +(0.5,-2) 
to [out=-120,in=70] cycle
;
\filldraw[ultra thick,draw=red,fill=red!40!white]
(1.1,3) to [out=-70,in=110]  (1.6,1.1) -- (2.6,1.1) to [out=110,in=-70] (2.1,3) -- cycle ;

\draw[ultra thick,draw=green]
(1.1,3) to [out=-70,in=110]  (1.6,1.1) ;

\draw[ultra thick, draw=white,double=white] (1.3,0.9) -- +(0.4,0) ;
\draw[ultra thick, draw=blue!80!black,->] (1.3,0.9) -- +(0.4,0) ;
\draw[ultra thick, draw=white, double=white] (1.35,2) -- +(0.4,0) ;
\draw[ultra thick, draw=green!80!black,->] (1.35,2) -- +(0.4,0) ;
\end{tikzpicture}
\caption[Pullback, 2-dim to 3-dim, 1-form]{The pullback with respect to a 2-dim submanifold (red plane)
  embedded in a 3-dim manifold. The pullback of an untwisted 1-form
  (blue planes) gives an untwisted 1-form (green lines).}
\label{fig_pull23_untw1}
\end{figure}
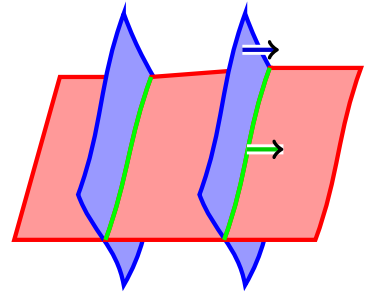


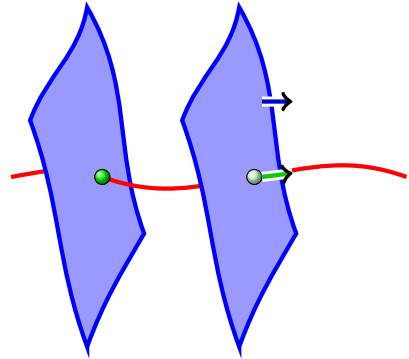
\begin{figure}[tb]
\centering
\begin{tikzpicture}[rotate=0]
\draw[ultra thick,draw=red] (-1,3) to[out=10,in=160] (0.2,3) ;
\filldraw[ultra thick,draw=blue,fill=blue!40!white,scale=1.5] (0,0.5) 
to [out=110,in=-70] +(-.5,2) 
to [out=70,in=-100] +(0.5,1) 
to [out=-60,in=110] +(0.5,-2) 
to [out=-120,in=70] cycle
;
\draw[ultra thick,draw=red] (0.2,3) to[out=-20,in=190] (2.2,3) ;
\filldraw[ultra thick,draw=blue,fill=blue!40!white,xshift=2cm,scale=1.5] (0,0.5) 
to [out=110,in=-70] +(-.5,2) 
to [out=70,in=-100] +(0.5,1) 
to [out=-60,in=110] +(0.5,-2) 
to [out=-120,in=70] cycle
;
\draw[ultra thick,draw=red] (2.2,3) to[out=10,in=160] +(2,0) ;
\shadedraw[ball color=green!800] (0.2,3) circle (0.1) ;
\shadedraw[ball color=green!20] (2.2,3) circle (0.1) ;

\draw[ultra thick, draw=white,double=white] (2.3,4) -- +(0.4,0) ;
\draw[ultra thick, draw=blue!80!black,->] (2.3,4) -- +(0.4,0) ;
\draw[ultra thick, draw=white,double=white] (2.3,3) -- +(0.4,0.05) ;
\draw[ultra thick, draw=green!80!black,->] (2.3,3) -- +(0.4,0.05) ;
\end{tikzpicture}
\caption[Pullback, 1-dim to 3-dim, 1-form]
{The pullback with respect to a 1-dim submanifold (red curve)
  embedded in a 3-dim manifold. 
  The pullback of an untwisted 1-form (blue plane) gives an untwisted 
  1-form (green dots).}
\label{fig_pull13_untw1}
\end{figure}


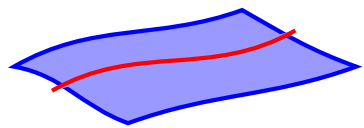
\begin{figure}[tb]
\begin{tikzpicture}
\filldraw[ultra thick,draw=blue,fill=blue!40!white,scale=1.5,rotate=90] (0,0) 
to [out=110,in=-70] +(-.5,2) 
to [out=70,in=-100] +(0.5,1) 
to [out=-60,in=110] +(0.5,-2) 
to [out=-120,in=70] cycle ;
\draw[ultra thick,draw=red,fill=blue!40!white,scale=1.6,rotate=90] (0.3,0.5) 
to [out=120,in=-60] +(-.5,2) ;
\end{tikzpicture}
\caption[Pullback, tangential]{The pullback with respect to a 1-dim submanifold (red curve)
  embedded in a 3-dim manifold.  The 1-form (blue plane) is tangential
  to 1-dim submanifold. The result is zero.}
\label{fig_pullback_tan}
\end{figure}

\begin{figure}[tb]
\centering
\begin{tikzpicture}[yscale=.7]
\filldraw[yshift=-2cm,fill=green!60!white] 
 (0,0) .. controls (1,1) and (2,-1) .. (3,0) 
       .. controls (3.2,0.7) and (4.,0.3) .. (4.5,1) 
       .. controls (3.5,0.3) and (2.5,1.3)  .. (1.5,1)  
       .. controls (1.,0.3) and (0.2,0.7)  .. (0,0) ;

\draw[red,ultra thick] (2.7,0.6) -- +(0,3.) ;
\draw[blue,ultra thick,yshift=-2cm] (2.7,0.6) circle (0.02);
\draw[red,ultra thick] (1,0.4) -- +(0,3.) ;
\draw[blue,ultra thick,yshift=-2cm] (1,0.4) circle (0.02) ;

\filldraw[xshift=1.5cm,color=brown,yshift=-0.1cm,yscale=-0.5] 
(0,0) -- (1,0) -- (1,1) -- (1.5,1) -- (0.5,2)  --
  (-0.5,1) -- (0,1) -- (0,0) ;

\draw (0,0) .. controls (1,1) and (2,-1) .. (3,0) ;
\draw (0,3) .. controls (1,3.7) and (2,2.3) .. (3,3) ;
\draw (0,0) -- (0,3) ;
\draw (3,0) -- (3,3) ;
\draw (3,3) .. controls (3.2,3.7) and (4.,3.3) .. (4.5,4) ;
\draw (0,3) .. controls (0.2,3.7) and (1.,3.3) .. (1.5,4) ;
\draw (3,0) .. controls (3.2,0.7) and (4.,0.3) .. (4.5,1) ;
\draw (1.5,4) .. controls (2.5,4.3) and (3.5,3.3) .. (4.5,4) ;
\draw (4.5,1) -- (4.5,4) +(.1,0) ;
\draw [dashed] (0,0) .. controls (0.2,0.7) and (1.,0.3) .. (1.5,1) ;
\draw [dashed] (1.5,1) .. controls (2.5,1.3) and (3.5,0.3) .. (4.5,1) ;
\draw [dashed] (1.5,1) -- (1.5,4) ;

\draw[very thick, draw=white,double]  (2.75,2.6) 
   [yscale=0.5] arc (80:-260:.4);
\draw[very thick,red!80!black]  (2.75,2.6) 
   [yscale=0.5] arc (80:-260:.4);
\draw[very thick,red!80!black]  (2.75,2.6) 
   [yscale=0.5,-<] arc (80:-80:.4);

\draw[very thick, draw=white,double]  (1.05,2.6) 
   [yscale=0.5] arc (80:-260:.4);
\draw[very thick,red!80!black]  (1.05,2.6) 
   [yscale=0.5] arc (80:-260:.4);
\draw[very thick,red!80!black]  (1.05,2.6) 
   [yscale=0.5,-<] arc (80:-80:.4);

\draw[very thick,blue]  (2.75,-1.25) 
   [yscale=0.5] arc (80:-260:.3);
\draw[very thick,blue]  (2.75,-1.25) 
   [yscale=0.5,-<] arc (80:-80:.3);

\draw[very thick,blue]  (1.05,-1.45) 
   [yscale=0.5] arc (80:-260:.3);
\draw[very thick,blue]  (1.05,-1.45) 
   [yscale=0.5,-<] arc (80:-80:.3);
\end{tikzpicture}
\caption[Pullback, 3-dim to 2-dim, 2-form]
{The pullback with respect to a projection from a 3-dim
  manifold to a 2-dim manifold (green plane).
  The pullback of the untwisted 2-form (blue dots) gives the 2-form
  (red lines).}
\label{fig_pull32_untw2}
\end{figure}
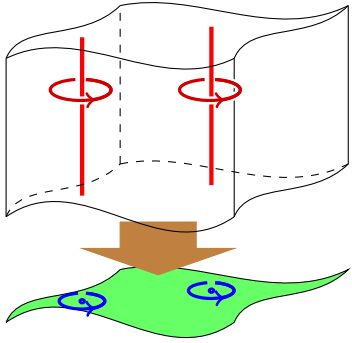


\begin{figure}[tb]
\centering
\begin{tikzpicture}[yscale=.7]
\filldraw[yshift=-2cm,fill=green!60!white] 
 (0,0) .. controls (1,1) and (2,-1) .. (3,0) 
       .. controls (3.2,0.7) and (4.,0.3) .. (4.5,1) 
       .. controls (3.5,0.3) and (2.5,1.3)  .. (1.5,1)  
       .. controls (1.,0.3) and (0.2,0.7)  .. (0,0) ;

\draw [dashed] (0,0) .. controls (0.2,0.7) and (1.,0.3) .. (1.5,1) ;
\draw [dashed] (1.5,1) .. controls (2.5,1.3) and (3.5,0.3) .. (4.5,1) ;
\draw [dashed] (1.5,1) -- (1.5,4) ;

\fill[fill=red,xshift=1cm] (0,0.2) -- (0,3.2) .. controls (0.2,3.7)
and (1.,3.3) .. (1.5,4) -- (1.5,1) .. controls (1.1,0.4) and (0.2,0.7)  ..  cycle ;

\filldraw[xshift=1.5cm,color=brown,yshift=-0.1cm,yscale=-0.5] 
(0,0) -- (1,0) -- (1,1) -- (1.5,1) -- (0.5,2)  --
  (-0.5,1) -- (0,1) -- (0,0) ;

\draw (0,0) .. controls (1,1) and (2,-1) .. (3,0) ;
\draw (0,3) .. controls (1,3.7) and (2,2.3) .. (3,3) ;
\draw (0,0) -- (0,3) ;
\draw (3,0) -- (3,3) ;
\draw (3,3) .. controls (3.2,3.7) and (4.,3.3) .. (4.5,4) ;
\draw (0,3) .. controls (0.2,3.7) and (1.,3.3) .. (1.5,4) ;
\draw (3,0) .. controls (3.2,0.7) and (4.,0.3) .. (4.5,1) ;
\draw (1.5,4) .. controls (2.5,4.3) and (3.5,3.3) .. (4.5,4) ;
\draw (4.5,1) -- (4.5,4) +(.1,0) ;

\draw[ultra thick,blue,xshift=1cm,yshift=-2cm] (1.5,1) .. controls (1.1,0.4) and
(0.2,0.7) .. (0,0.2)   ; 

\draw[ultra thick, draw=white,double=white] (2.3,3.2) -- +(0.6,-0.1) ;
\draw[ultra thick, draw=red!80!black,->] (2.3,3.2) -- +(0.6,-0.1) ;
\draw[ultra thick, draw=blue!80!black,->] (2.1,-1.3) -- +(0.4,-0.3) ;

\end{tikzpicture}
\caption[Pullback, 3-dim to 2-dim, 1-form]
{The pullback with respect to a projection from a 3-dim
  manifold to a 2-dim manifold (green plane).
  The pullback of the untwisted 1-form (blue curve) gives the 1-form
  (red plane).}
\label{fig_pull32_untw1}
\end{figure}
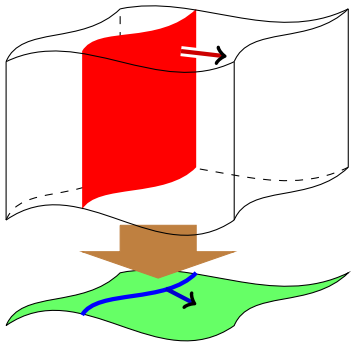


\begin{figure}[tb]
\centering
\begin{tikzpicture}[yscale=.7]
\draw[ultra thick,yshift=-1.5cm,xshift=1cm,draw=green!80!black] 
 (0,0) .. controls (1,1) and (2,-1) .. (3,0) ;

\filldraw[xshift=1.5cm,color=brown,yshift=0cm,yscale=-0.5] 
(0,0) -- (1,0) -- (1,1) -- (1.5,1) -- (0.5,2)  --
  (-0.5,1) -- (0,1) -- (0,0) ;

\draw [dashed] (0,0) -- (1.5,1) ;
\draw [dashed] (1.5,1) .. controls (2.5,1.3) and (3.5,0.3) .. (4.5,1) ;
\draw [dashed] (1.5,1) -- (1.5,4) ;

\fill[fill=red,xshift=1cm] (0,0.2) -- (0,3.2) -- (1.5,4) -- (1.5,1) -- cycle ; 

\draw (0,0) .. controls (1,1) and (2,-1) .. (3,0) ;
\draw (0,3) .. controls (1,3.7) and (2,2.3) .. (3,3) ;
\draw (0,0) -- (0,3) ;
\draw (3,0) -- (3,3) ;
\draw (3,3) -- (4.5,4) +(.1,0);
\draw (0,3) -- (1.5,4) ;
\draw (3,0) -- (4.5,1) ;
\draw (1.5,4) .. controls (2.5,4.3) and (3.5,3.3) .. (4.5,4) ;
\draw (4.5,1) -- (4.5,4) ;

\draw[blue,ultra thick] (2,-1.3) circle (0.02) ;

\draw[ultra thick, draw=white,double=white] (2.3,3.2) -- +(0.6,-0.10) ;
\draw[ultra thick, draw=red!80!black,->] (2.3,3.2) -- +(0.6,-0.10) ;
\draw[ultra thick, draw=blue!80!black,->] (2.1,-1.3) -- +(0.4,-0.2) ;

\end{tikzpicture}
\caption[Pullback, 3-dim to 1-dim, 1-form]
{The pullback with respect to a projection from a 3-dim
  manifold to a 1-dim manifold (green curve).
  The pullback of the untwisted 1-form (blue dot) gives the 1-form
  (red plane).}
\label{fig_pull31_untw1}
\end{figure}

No matter how complicated the map
is it is always possible to pullback the form. We will only consider
here the pullback with respect to diffeomorphisms, embeddings and
projections. Pullbacks preserve the degree of the form and also the 
twistedness of the forms. 

The pullback with respect to diffeomorphisms simply deforms the shape
of the form-submanifolds. However unless there is a metric the
question is: deforms with respect to what? This corresponds to our
statement that without a metric one cannot distinguish between
manifolds which are diffeomorphic and therefore as we said in the
introduction, one is dealing with rubber sheet geometry.

For the pullback with respect to an embedding, one simply takes the
intersection of submanifold with the form-submanifolds as seen in
figures \ref{fig_pull23_untw2}, \ref{fig_pull23_untw1} and
\ref{fig_pull13_untw1}.
This preserves the degree of the
form. If we pullback a $p$-form from an $m$-dimensional manifold to an
$n$-dimensional manifold then the $p$-form-submanifolds are of
dimension $(m-p)$, hence the intersection of the $(m-p)$-dim
form-submanifolds with the $n$-dimensional manifolds will have
dimension $(n-p)$ corresponding correctly to a $p$-form.
This assumes the form-submanifolds and the submanifold intersects
transversely. However if the submanifold and the
form-submanifold intersect tangentially then the result is zero. See
figure \ref{fig_pullback_tan}.

For untwisted forms there is a
natural pullback of the orientation. For twisted forms it is necessary
to pullback onto a twisted submanifold. To be guaranteed to get the
orientation correct it is easiest to choose an orientation for the
embedding manifold, untwist both the form and the submanifold and then
twist the pullbacked form. This gives the result that: The orientation
of the pullbacked form concatenated with the orientation of the
submanifold equals the orientation of the original form. In particular
the pullback of a twisted $p$-form onto a twisted submanifold is plus
if the orientations agree and minus otherwise, as was discussed in
section \ref{ch_Integration}.

For the pullback with respect to projections the resulting
form-submanifold is simply all the points which mapped onto the
original form-submanifold. See figures \ref{fig_pull32_untw2} and
\ref{fig_pull32_untw1}. For untwisted forms the orientation is
obvious. For twisted forms one can use orientation for both the source
and target manifolds, untwist the form and convert it
back. Alternatively one can decide an orientation of the projection
and use that. We will not describe this here.


\subsection{Pullbacks of metrics.}
\label{ch_PullBacks_Metric}

\begin{figure}[tb]
\centering
\begin{tikzpicture}
\filldraw[fill=red!50,draw=red] (0,0) to[out=10,in=190] (5,0) -- (7,2)
to[out=-170,in=10]  (2,2) -- cycle ;

\shadedraw [ball color=green,shift={(2.4,1.0)}]
(0,0) [yscale=0.5] arc(0:-180:.5) (0,0) [rotate=-30,yscale=2.5] arc
(15:190:.46) ;
\draw [ultra thick,blue,shift={(2.4,1.0)}]
(0,0) [yscale=0.5] arc(0:-180:.5)  ;
\shadedraw [ball color=green,shift={(3.9,1.0)}]
(0,0) [yscale=0.5] arc(0:-180:.5) (0,0) [rotate=-15,yscale=2.8] arc (7:184:.48) ;
\draw [ultra thick,blue,shift={(3.9,1.0)}]
(0,0) [yscale=0.5] arc(0:-180:.5)  ;
\shadedraw [ball color=green,shift={(5.4,1.0)}]
(0,0) [yscale=0.5] arc(0:-180:.5) (0,0) [yscale=3.0,rotate=-0] arc (0:180:.5) ;
\draw [ultra thick,blue,shift={(5.4,1.0)}]
(0,0) [yscale=0.5] arc(0:-180:.5)  ;

\end{tikzpicture}
\caption[The pullback of a metric onto an embedding.]{The
  pullback of a Riemannian metric onto an embedding as represented by
  unit spheres. The unit spheres in 3-dim (green) are pulled back into
  unit circles (blue) on the 2-dim submanifold (red).}
\label{fig_metric_pullback}
\end{figure}
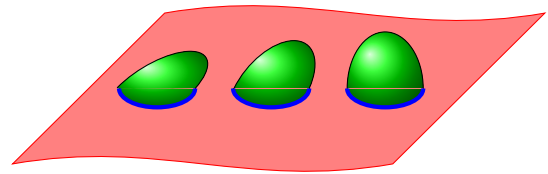

We can only pullback a metric with respect to diffeomorphisms and
embedding. With respect to diffeomorphisms it makes the two manifolds
have the same geometry.

The pullback of a metric with respect to an embedding is called the
induced metric. This is the metric we are familiar with on the 
sphere or torus which is embedded into flat 3-space.
Figure \ref{fig_metric_pullback} shows the result of the induced metric.

We cannot pullback a metric with respect to a projection as we don't
have enough information. In particular we do not know the length of
vectors in the direction of the projection.

\subsection{Pushforwards of vectors.}
\label{ch_PushForwards}

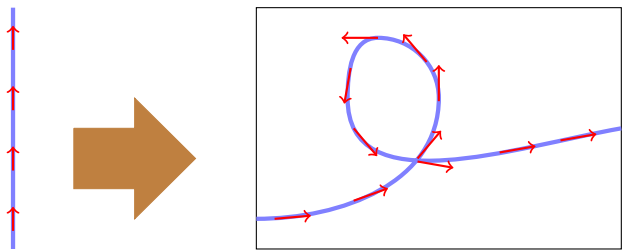
\begin{figure}[tb]
\centering
\begin{tikzpicture}[scale=0.8]
\draw[ultra thick, blue!50] (0,0) -- (0,4) ;
\draw[thick,red,->] (0,.3) -- +(0,.4) ; 
\draw[thick,red,->] (0,1.3) -- +(0,.4) ; 
\draw[thick,red,->] (0,2.3) -- +(0,.4) ; 
\draw[thick,red,->] (0,3.3) -- +(0,.4) ; 

\filldraw[xshift=1cm,yshift=2cm,color=brown,yscale=1,rotate=-90] 
(0,0) -- (1,0) -- (1,1) -- (1.5,1) -- (0.5,2)  --
  (-0.5,1) -- (0,1) -- (0,0) ;

\draw (4,0) -- (10,0) -- (10,4) -- (4,4) -- cycle ;

\draw[ultra thick, blue!50,shift={(4,0)}] 
(0,.5)  to[out=0,in=-90]  (3,2.5) 
to[out=90,in=0] (2,3.5) 
to[out=180,in=90] (1.5,2.5) 
to[out=-90,in=190] (6,2) ;


\draw[thick,red,->,shift={(4,0)},shift={(0.3,0.5)}] (0,0) to [rotate=5] +(.6,0) ;
\draw[thick,red,->,shift={(4,0)},shift={(1.6,0.8)}] (0,0) to [rotate=20] +(.6,0) ;
\draw[thick,red,->,shift={(4,0)},shift={(2.65,1.5)}] (0,0) to [rotate=50] +(.6,0) ;
\draw[thick,red,->,shift={(4,0)},shift={(3.0,2.45)}] (0,0) to [rotate=90] +(.6,0) ;
\draw[thick,red,->,shift={(4,0)},shift={(2.8,3.1)}] (0,0) to [rotate=130] +(.6,0) ;
                                                 
\draw[thick,red,->,shift={(4,0)},shift={(2.0,3.5)}] (0,0) to [rotate=180] +(.6,0) ;
\draw[thick,red,->,shift={(4,0)},shift={(1.55,3.0)}] (0,0) to [rotate=-100] +(.6,0) ;
\draw[thick,red,->,shift={(4,0)},shift={(1.6,2.0)}] (0,0) to [rotate=-50] +(.6,0) ;
\draw[thick,red,->,shift={(4,0)},shift={(2.65,1.45)}] (0,0) to [rotate=-10] +(.6,0) ;
\draw[thick,red,->,shift={(4,0)},shift={(4.0,1.6)}] (0,0) to [rotate=10] +(.6,0) ;
\draw[thick,red,->,shift={(4,0)},shift={(5.0,1.8)}] (0,0) to [rotate=10] +(.6,0) ;

\end{tikzpicture}
\caption[The Push forward of the vector field.]  {Push forward of the
  vector field (red) of the 1-dim manifold (blue) embedding into the
  2-dim manifold (white) produces vectors (red). However these vectors
  do not form a vector field as they are not defined in some places
  and in others have multiple values.}
\label{fig_pushforward_vectors}
\end{figure}

We can always pushforward a vector field with respect to a
diffeomorphism. However, in general the pushforward of a vector field
with respect to any other map will not result in a vector field. As
can be seen in figure \ref{fig_pushforward_vectors} some points in the
target manifold have no vector associated with them and others have
more than one.

However the pushforward of vectors are very useful. For example the
velocity of a single particle is the tangent to its worldline. This is
only defined on its worldline. 

I will not discuss these further here.

\section{Electromagnetism and Electrodynamics}
\label{ch_Maxwell}

\begin{table*}[tb]
\centering
\begin{tabular}{|c@{\quad}l|l|c@{\quad}l|l|}
\hline
\multicolumn{2}{|c|}{\multi{Vector\\notation}} & 
names &
\multicolumn{2}{|c|}{\multi{Form\\ notation}} & 
\multi{Relation between\\ vector  and form} 
\\\hline
$\VE$ & \multi{\rule{0em}{1em}untwisted\\vector} &
Electric field. &
$E$ &\multi{untwisted\\1-form} & 
Metric dual 
\\\hline
$\VD$ & \multi{\rule{0em}{1em}untwisted\\vector} &
Displacement field. &
$D$ & \multi{twisted\\2-form}  & 
\multi{Metric dual and\\Hodge dual} 
\\\hline
$\VB$ & \multi{twisted\\vector} & 
\multi{``B-field'', Magnetic field or \\Magnetic flux density.}
&
$B$ & \multi{untwisted\\2-form} & 
\multi{Metric dual and\\Hodge dual} 
\\\hline
$\VH$ & \multi{\rule{0em}{1em}twisted\\vector} &
\multi{\rule{0em}{1em}``H-field'', Magnetic field intensity, 
    \\  Magnetic field strength,
    \\Magnetic field  or
    Magnetising field.} &
$H$ & \multi{untwisted\\1-form} & 
Metric dual 
\\\hline
$\rho$ & \multi{untwisted\\scalar} &
charge density.&
$\varrho$ & \multi{twisted\\3-form} & 
Hodge dual 
\\\hline
$\VJ$ & \multi{untwisted\\vector} &
current density. &
$J$ & \multi{twisted\\2-form}& 
\multi{Metric dual and\\Hodge dual} 
\\\hline
\end{tabular}
\caption[The fields in electromagnetism]{The fields in
  electromagnetism, in terms of the familiar vector notation, the
  corresponding form notation and the operations needed to pass
  between the two. }
\label{tab_Max_Fields}
\end{table*}


\begin{table}[tb]
\begin{tabular}{|l@{\ \ }l|l@{\ \ }l|l@{\ \ }l|}
\hline
\multicolumn{2}{|c|}{\multi{Spacetime\\form}}
&
\multicolumn{2}{|c|}{\multi{Pull back\\ onto space}}
&
\multicolumn{2}{|c|}{\multi{Remaining\\component}}
\\\hline
$F$ & \multi{untwisted\\ 2-form} & 
$B$ &  \multi{untwisted\\ 2-form} & 
$E$ & \multi{untwisted\\ 1-form}
\\\hline
$\Hform$ & \multi{twisted\\2-form} &  
$D$ &\multi{twisted\\2-form}  & 
$H$ &\multi{twisted\\1-form}
\\\hline
$\Jcurr$ & \multi{twisted\\3-form} &   
$\varrho$ & \multi{twisted\\3-form} & 
$J$ &\multi{twisted\\2-form}
\\\hline
\end{tabular}
\caption[The relativistic field in electromagnetism]
{Combining the fields in electromagnetism into the
  relativistic quantities.}
\label{tab_fields-FHJ}
\end{table}

Maxwell's equations, when written in standard Gibbs-Heaviside notation
are four differential equations of four vector fields $\VE$, $\VB$, $\VD$
and $\VH$. The source for the electromagnetic fields is the charge
density $\rho$ and the current $\VJ$. In our notation it is clearer to
replace these with forms. In table \ref{tab_Max_Fields}
the relationship between the vector notation and the form notation is
given.

As can be seen in table \ref{tab_Max_Fields},
both the fields $\VB$ and $\VH$ may be
referred to as the magnetic field. Here we will simply refer to them as
the ``B-field'' and ``H-field''.  


As we said, Maxwell's equations are differential equations in $\VE$,
$\VB$, $\VD$ and $\VH$. However there is insufficient information to
be able to solve these equations. They need to be related by
additional information which relate the four quantities $\VE$,
$\VB$, $\VD$ and $\VH$. These are called the constitutive
relations. We have just stated the vacuum constitutive relations,
namely that $\VD=\epsilon_0\VE$ and $\VB=\mu_0\VH$. However in a
medium, such as glass, water, a plasma, a magnet, a human body or
whatever else you may wish to study, these constitutive relations are
more complicated. The simplest constitutive relations consider a
constant permittivity $\epsilon$ and permeability $\mu$. In this medium
there is no dispersion, i.e. radiation of different frequency travel
at the same speed. This would mean for example that rain drops would
not give rise to rainbows and may therefore be described as {\em
  ``antediluvian''}! They are however easy to work with.
A more physically reasonable constitutive relations is to allow the
permittivity or permeability to depend on frequency. This is a good
model for water, wax or the human body. If the electric or magnet
fields are very strong then the permittivity and permeability depend
on the electric or magnet fields. These are non linear
materials. Exotic and metamaterials can have more complicated
constitutive relations. The degree of complexity of these relations is
limited only by the imagination of the scientists studying the material
and their ability to build them.
Indeed physicists are currently looking at models of the vacuum where
the simple relations $\VD=\epsilon_0\VE$ and $\VB=\mu_0\VH$ break down
in the presents of intense fields, either due to quantum effects or
because Maxwell's equations themselves breakdown.

Just as in relativistic mechanics where we combine energy and momentum
into a single relativistic vector the 4-momentum, so we combine the
electric and magnetic B-field into a single spacetime untwisted 2-form
$F$. Likewise we combine the displacement current and H-field into a
single twisted 2-form $\Hform$, and we combine $\rho$ and $J$ in a
single twisted 3-form $\Jcurr$, as summarised in table
\ref{tab_fields-FHJ}.  As a result relativity mixes the $\VE$ and
$\VB$ fields and separately mixes the $\VD$ and $\VH$ fields. The
constitutive relations then relate $F$ and $\Hform$. The simplest case
is the vacuum where the relationship is simply given by the Hodge
dual.

As stated, Maxwell's equations give the electric and magnetic fields
for a given current. It is also necessary to know how the charges,
that produce the current, respond to the electromagnetic fields. This
is given by the Lorentz force equation.

In this section we start, section \ref{ch_Max_statics},
with electrostatics and magnetostatics as these are in three
dimensions and therefore easier to depict pictorially. We then, in
section \ref{ch_Max_Maxwell} to draw Maxwell's equations as four
spacetime pictures. In section \ref{ch_Max_Lorentz}, we depict the
Lorentz force equation and a spacetime picture. In section
\ref{ch_Max_3dim}
we look at Maxwell's equations in 3-dimensional
spacetime. This is the pullback of the 4-dimensional spacetime
equations onto a fixed plane. This contrasts the the static cases when
we pullback onto a fixed timeslice.

\subsection{Electrostatics and Magnetostatics.}
\label{ch_Max_statics}

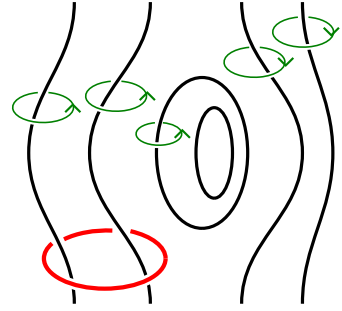
\begin{figure}[tb]
\centering
\begin{tikzpicture}[scale=2]

\draw [very thick] (0,0) to[out=90,in=-90] +(-0.3,1) to[out=90,in=-90] +(0.3,1) ;
\draw [very thick] (.5,0) to[out=90,in=-90] +(-0.4,1) to[out=90,in=-90] +(0.4,1) ;
\draw [very thick,xscale=0.6] (1.4,.5) arc(-90:270:0.5) ;
\draw [very thick,xscale=0.4] (2.3,.7) arc(-90:270:0.3) ;
\draw [very thick] (1.1,0) to[out=90,in=-90] +(0.4,1) to[out=90,in=-90] +(-0.4,1) ;
\draw [very thick] (1.5,0) to[out=90,in=-90] +(0.2,1) to[out=90,in=-90] +(-0.2,1) ;
\draw[thick,yscale=-0.5,draw=white,double=green!50!black] 
(-0.24,-2.8) arc (260:-80:0.2) ; 
\draw[thick,draw=green!50!black] 
(-0.06,1.27) -- +(0.05,0.05) -- +(0.1,0) ; 

\draw[thick,yscale=-0.5,draw=white,double=green!50!black] 
(0.24,-2.95) arc (260:-80:0.2) ; 
\draw[thick,draw=green!50!black] 
(0.43,1.35) -- +(0.05,0.05) -- +(0.1,0) ; 

\draw[thick,yscale=-0.5,draw=white,double=green!50!black] 
(0.53,-2.4) arc (260:-80:0.15) ; 
\draw[thick,draw=green!50!black] 
(0.66,1.1) -- +(0.05,0.05) -- +(0.1,0) ; 

\draw[thick,yscale=-0.5,draw=white,double=green!50!black] 
(1.15,-3.4) arc (260:-80:0.2) ; 
\draw[thick,draw=green!50!black] 
(1.34,1.63) -- +(0.05,-0.05) -- +(0.1,0) ; 

\draw[thick,yscale=-0.5,draw=white,double=green!50!black] 
(1.47,-3.8) arc (260:-80:0.2) ; 
\draw[thick,draw=green!50!black] 
(1.64,1.83) -- +(0.05,-0.05) -- +(0.1,0) ; 

\draw[ultra thick,double=white,white,yscale=0.5] (0.6,0.6) arc(0:-210:0.4) ; 
\draw[ultra thick,red,yscale=0.5] (0.6,0.6) arc(0:-210:0.4) ;
\draw[ultra thick,red,yscale=0.5] (0.6,0.6) arc(0:65:0.4) ;
\draw[ultra thick,red,yscale=0.5] (0.25,.97) arc(85:135:0.4) ;

\end{tikzpicture}

\caption[The magnetic 2-form in space.]{The magnetic 2-form $B$ in space as
  unbroken curves. Since the curves are unbroken this corresponds to
  the non divergence of $B$.  The
  external orientation is given (green). The circle (red) encloses some
  of the $B$-lines.}
\label{fig_EM_Mag_2-form}
\end{figure}


\begin{figure}[tb]
\centering
\begin{tikzpicture}[scale=2]
\draw [thick] (0,0) -- +(1,.3) (1.1,.4) node{$D$};
\draw [thick,-<] (0,0) -- +(1*.5,.3*.5) ;
\draw [thick] (0,0) -- +(.5,1) ;
\draw [thick,-<] (0,0) -- +(.5*.5,1*.5) ;
\draw [thick] (0,0) -- +(1,-.3) ;
\draw [thick,-<] (0,0) -- +(1*.5,-.3*.5) ;
\draw [thick] (0,0) -- +(.5,-1) ;
\draw [thick,-<] (0,0) -- +(.5*.5,-1*.5) ;
\draw [thick] (0,0) -- +(-1,-.3) ;
\draw [thick,-<] (0,0) -- +(-1*.5,-.3*.5) ;
\draw [thick] (0,0) -- +(-.5,-1) ;
\draw [thick,-<] (0,0) -- +(-.5*.5,-1*.5) ;
\draw [thick] (0,0) -- +(-1,.3) ;
\draw [thick,-<] (0,0) -- +(-1*.5,.3*.5) ;
\draw [thick] (0,0) -- +(-.5,1) ;
\draw [thick,-<] (0,0) -- +(-.5*.5,1*.5) ;
\shadedraw [ball color=red,draw=black] (0,0) circle (0.1) ; 
\shadedraw [ball color=green,draw=black,opacity=0.4] (0,0) circle (0.8) ;
\shadedraw [ball color=green,draw=black,opacity=0.4] (0,0) circle (0.4) ;
\shadedraw [ball color=red,draw=black,opacity=0.4] (0,0) circle (0.1) ; 
\draw [green!70!black] (0,0.9) node{$\boldsymbol E$} ;
\draw [ultra thick,green!70!black,->] (0,0.8) -- +(0,-.3) ;
\end{tikzpicture}

\caption[Gauss's Law in electrostatics.]{Gauss's Law in
  electrostatics. The positive charge, is given by the twisted 3-form
  $\varrho$ is the dot (red) in the centre. This is the exterior
  derivative of the twisted 2-form $D$ (black) lines. In the vacuum or
  similar constitutive relations, then $D$ is the Hodge dual of the
  untwisted 1-form $E$ (green spheres). }
\label{fig_EM_Gauss_ElectroStatics}
\end{figure}
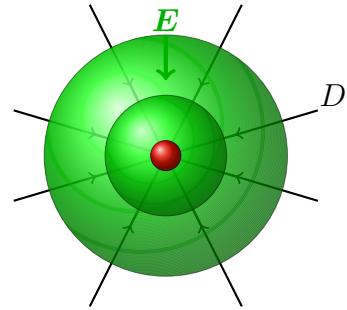

\begin{figure}[tb]
\centering
\begin{tikzpicture}
\fill [red,opacity=0.3] (0,0) -- +(0,4) -- +(2,5) -- +(2,1) -- cycle ;
\fill [red,opacity=0.3] (0,0) -- +(0,4) -- +(3,4) -- +(3,0) -- cycle ;
\fill [red,opacity=0.3] (0,0) -- +(0,4) -- +(1.5,3) -- +(1.5,-1) -- cycle ;
\fill [red,opacity=0.3,xscale=-1] (0,0) -- +(0,4) -- +(2,5) -- +(2,1) -- cycle ;
\fill [red,opacity=0.3,xscale=-1] (0,0) -- +(0,4) -- +(3,4) -- +(3,0) -- cycle ;
\fill [red,opacity=0.3,xscale=-1] (0,0) -- +(0,4) -- +(1.5,3) -- +(1.5,-1) -- cycle ;

\draw [ultra thick,double=black] (0,-1) -- (0,0.4) ;
\draw [ultra thick,double=black] (0,0.6) -- (0,2.4) ;
\draw [ultra thick,double=black] (0,2.6) -- (0,5) ;
\draw [ultra thick,double=black,->] (0,0.6) -- (0,2) ;

\draw [ultra thick,blue,yscale=0.5] (.2,3) arc (80:-260:1) ; 
\draw [ultra thick,blue!50!red,yscale=0.5] (.2,3) arc (80:30:1) ; 
\draw [ultra thick,blue!70!red,yscale=0.5] (0.89,2.52) arc (30:-30:1) ; 
\draw [ultra thick,blue!50!red,yscale=0.5] (-.15,3) arc (100:150:1) ; 
\draw [ultra thick,blue!70!red,yscale=0.5] (-0.84,2.52) arc (150:210:1) ;

\draw [ultra thick,blue,yscale=0.5,yshift=4cm] (.2,3) arc (80:-260:1) ; 
\draw [ultra thick,blue!50!red,yscale=0.5,yshift=4cm] (.2,3) arc (80:30:1) ; 
\draw [ultra thick,blue!70!red,yscale=0.5,yshift=4cm] (0.89,2.52) arc (30:-30:1) ; 
\draw [ultra thick,blue!50!red,yscale=0.5,yshift=4cm] (-.15,3) arc (100:150:1) ; 
\draw [ultra thick,blue!70!red,yscale=0.5,yshift=4cm] (-0.84,2.52) arc (150:210:1) ; 

\draw [very thick,green!80!black] (0.89,0.65) arc (-110:210:.3) ; 
\draw [very thick,green!80!black,-<] (0.89,0.65) arc (-110:0:.3) ; 

\draw [very thick,red!80!black,<-] (2.,0.65) arc (-110:210:.3) ; 

\end{tikzpicture}
\caption[The $B$ and $H$ fields from a wire.]{Ampere's law: The constant line current, the twisted
  2-form $J$ (Black), is the exterior derivative of the twisted 1-form
  $H$.  In the vacuum or similar constitutive relations, then $H$ is
  the Hodge dual of the untwisted 2-form $B$ (blue circles).  }
\label{fig_EM_Ampere_MagnetoStatics}

\end{figure}
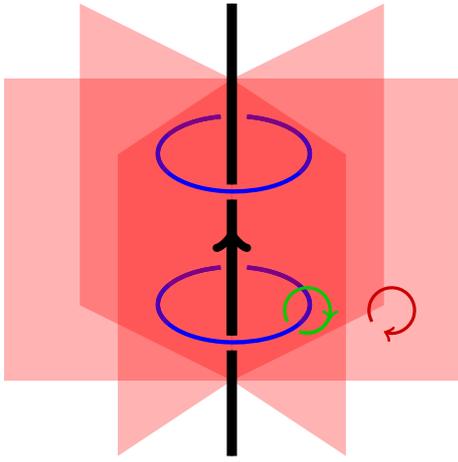


\begin{figure}[b!]
\begin{tikzpicture}[scale=0.8]

\filldraw [ultra thick,black,fill=black!50] (-1,-2) -- +(2,0) -- +(2,4) -- +(0,4) -- cycle ;
\draw [ultra thick,blue,xscale=0.3] (0.8,-2) arc (-138.3:138.3:3) 
  (0.8,-2) -- (0.8,2) (0.8,-2) -- (0.8,0) ; 
\draw [ultra thick,blue,xscale=-0.3] (0.8,-2) arc (-138.3:138.3:3) 
  (0.8,-2) -- (0.8,2) (0.8,-2) -- (0.8,0); 

\shade [ball color=red,opacity=0.4,yscale=0.35] (1,4.) arc (-163.5+90:163.5+90:3.5) 
  (1,4.6)--(-1,4.6) ;
\shade [ball color=red,opacity=0.4,yscale=-0.35] (1,4.) arc (-163.5+90:163.5+90:3.5) 
  (1,4.6)--(-1,4.6) ;
\shade [ball color=red,opacity=0.4,yscale=-0.35] (1,3.) arc (-168.5+90:168.5+90:5.) 
  (1,4.6)--(-1,4.6) ;
\shade [ball color=red,opacity=0.4,yscale=0.35] (1,3.) arc (-168.5+90:168.5+90:5.) 
  (1,4.6)--(-1,4.6) ;

\draw [line width=3pt,white,yscale=0.5] (1.85,-2.5) arc (70:-250:0.4) ; 
\draw [very thick,red!20!black,yscale=0.5] (1.85,-2.5) arc (70:-250:0.4) ; 
\draw [very thick,red!20!black,yscale=0.5,->] (1.85,-2.5) arc (70:-60:0.4) ; 

\draw [ultra thick,double=white,white,yscale=0.5] (0.4,-2.5) arc (70:-250:0.4) ; 
\draw [very thick,red!20!black,yscale=0.5] (0.4,-2.5) arc (70:-250:0.4) ; 
\draw [very thick,red!20!black,yscale=0.5,->] (0.4,-2.5) arc (70:-60:0.4) ; 
\draw (0,2.5) node {N} ;
\draw (0,-2.5) node {S} ;

\draw [ultra thick,double=white,white,yscale=0.5] (1.95,0) arc (70:-250:0.4) ; 
\draw [very thick,green!80!black,yscale=0.5] (1.95,0) arc (70:-250:0.4) ; 
\draw [very thick,green!80!black,yscale=0.5,->] (1.95,0) arc (70:-60:0.4) ; 

\draw [ultra thick,double=white,white,yscale=0.5] (0.4,0) arc (70:-250:0.4) ; 
\draw [very thick,green!80!black,yscale=0.5] (0.4,0) arc (70:-250:0.4) ; 
\draw [very thick,green!80!black,yscale=0.5,-<] (0.4,0) arc (70:-60:0.4) ; 

\end{tikzpicture}
\caption[The $B$ and $H$ fields from a magnet.]{A
  permanent magnet.  The untwisted magnetic 2-form $B$ (blue curves)
  and twisted 1-form $H$ (red ellipsoids) generated by a permanent
  magnet (grey rectangle). Observe that outside of the magnet $H$
  fields are perpendicular to the $B$ fields. Thus the Hodge dual,
  figure \ref{fig_Hodge_dual}, of $H$ coincides with $B$. Inside the
  magnet however the orientations of $H$ and $B$ are in opposite
  directions. This is due to the non vacuum constitutive relations of
  permanent magnets.}
\label{fig_Mag_Bulk}
\end{figure}
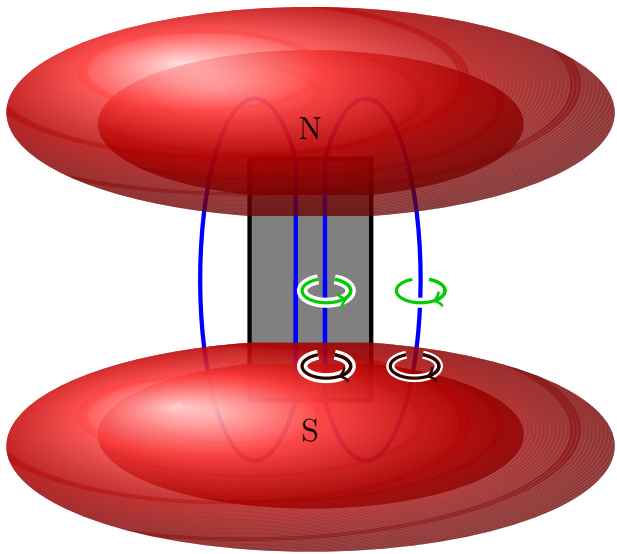



In this subsection we deal with electrostatics and
magnetostatics as 3-dimensional pictures are easier to draw.
One has that $\VB$ is divergence-free, Gauss's law which
states that the divergence of $\VD$ equals the charge, Faraday's law
which states that $\VE$ is curl-free and Ampere's law which states
that the curl of $\VH$ is the current.  In the language of forms this
corresponds to the fact that $B$ and $E$ are closed, whereas the
exterior derivative of $D$ and $H$ are $\rho$ and $J$ respectively.

In figure \ref{fig_EM_Mag_2-form} we see that the $B$ lines are
unbroken and therefore $B$ is closed. In figure
\ref{fig_EM_Gauss_ElectroStatics} we see Gauss's law, that a charge
twisted 3-form $\varrho$ is the exterior derivative of the twisted
2-form $D$. For the vacuum or for any constitutive relation where the
permittivity is constant or depends only on frequency, then the $D$ is
the permittivity times the Hodge dual of the electric 1-form $E$.  In
figure \ref{fig_EM_Ampere_MagnetoStatics} we see Ampere's law. The
current, a twisted 2-form $J$ in the wire, is the exterior derivative
of the twisted 1-form $H$, and $H$ is the Hodge dual of $B$.

In figure \ref{fig_Mag_Bulk} we show the fields for a permanent
magnet. This is an example of a slightly more complicated constitutive
relations. Outside of the magnet, the twisted 1-form $H$ is the Hodge
dual of the untwisted 2-form $B$. However in the magnet, the
orientations of $H$ and $B$ are in the opposite direction, due to the
constitutive relations for a permanent magnet.

\subsection{Maxwell's equations as four spacetime pictures}
\label{ch_Max_Maxwell}

\begin{figure}[tb]
\begin{tikzpicture}
\draw[->] (0,0) -- (0,5) ;
\draw (-0.3,2) node [rotate=90] {Time} ;
\draw[->] (0,0) -- (5,0) ;
\draw (3,-0.3) node [rotate=0] {Space} ;
\draw[->,rotate=-15] (0,0) -- (5,0) ;
\draw[rotate=-15]  (3,-0.3) node [rotate=-15] {Space} ;
\draw[->,rotate=15] (0,0) -- (5,0) ;
\draw[rotate=15]  (3,-0.3) node [rotate=15] {Space} ;

\draw[ultra thick,green!60!black,->-=.7] (2,0.2) 
to[out=100,in=-80] +(0,4)  ;
\draw[ultra thick,green!60!black,->-=.7] (2.5,0.) 
to[out=90,in=-90] +(0.3,4) ;
\draw[ultra thick,green!60!black,->-=.7] (3,0.2) 
to[out=90,in=-90] +(0.6,4) ;
\draw[ultra thick,green!60!black,->-=.7] (3.5,0.) 
to[out=90,in=-90] +(0.9,4) ;
\draw[ultra thick,green!60!black,->-=.7] (4,0.2) 
to[out=90,in=-90] +(1.15,4) ;
\draw[ultra thick,green!60!black,->-=.7] (4.5,0.) 
to[out=90,in=-90] +(1.35,4) ;
\draw[ultra thick,green!60!black,->-=.7] (5,0.2) 
to[out=90,in=-90] +(1.6,4) ;

\end{tikzpicture}
\caption[Conservation of charge]{Conservation of charge: There exists
  a closed twisted 3-form ${\cal J}$ called the electric current
  (green). The green line may be interpreted as charged particles,
  distributed over all of spacetime.}
\label{fig_Max_J}
\end{figure}
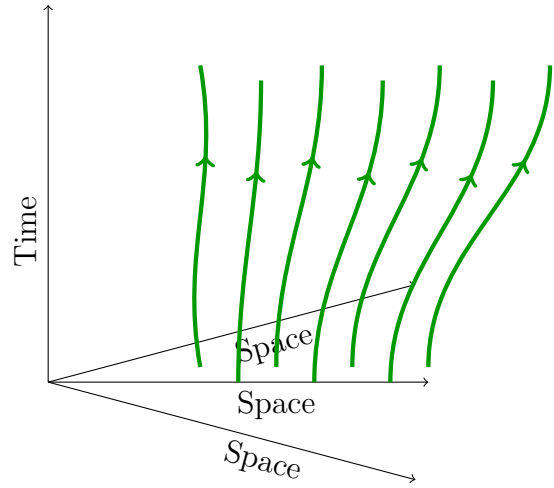

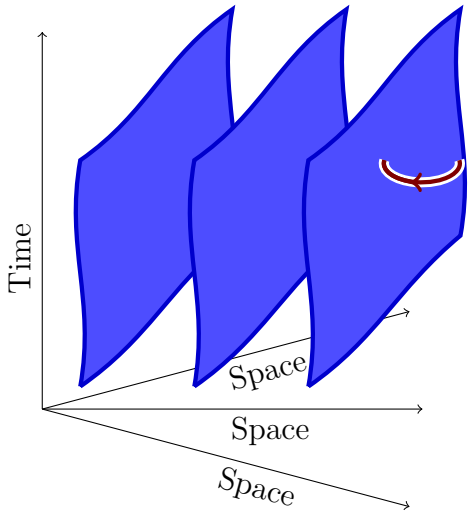
\begin{figure}[tb]
\begin{tikzpicture}
\draw[->] (0,0) -- (0,5) ;
\draw (-0.3,2) node [rotate=90] {Time} ;
\draw[->] (0,0) -- (5,0) ;
\draw (3,-0.3) node [rotate=0] {Space} ;
\draw[->,rotate=-15] (0,0) -- (5,0) ;
\draw[rotate=-15]  (3,-0.3) node [rotate=-15] {Space} ;
\draw[->,rotate=15] (0,0) -- (5,0) ;
\draw[rotate=15]  (3,-0.3) node [rotate=15] {Space} ;

\filldraw[ultra thick,draw=blue!80!black,fill=blue!70!white,shift={(.5,.3)}] 
(0,0) to[out=80,in=-100] (0,3) to [out=35,in=-145] (2,5) 
to [out=-100,in=80] (2,2) to [out=-140,in=35] (0,0); 
\filldraw[ultra thick,draw=blue!80!black,fill=blue!70!white,shift={(2,.3)}] 
(0,0) to[out=80,in=-100] (0,3) to [out=35,in=-145] (2,5) 
to [out=-100,in=80] (2,2) to [out=-140,in=35] (0,0); 
\filldraw[ultra thick,draw=blue!80!black,fill=blue!70!white,shift={(3.5,.3)}] 
(0,0) to[out=80,in=-100] (0,3) to [out=35,in=-145] (2,5) 
to [out=-100,in=80] (2,2) to [out=-140,in=35] (0,0); 
\draw[ultra thick,white,double=white,shift={(3.5,.3)},yscale=0.5] (1.,6) arc
(-190:10:0.5) ;
\draw[ultra thick,red!50!black,shift={(3.5,.3)},yscale=0.5] (1.,6) arc
(-190:10:0.5) ;
\draw[ultra thick,red!50!black,shift={(3.5,.3)},-<,yscale=0.5] (1.,6) arc
(-190:-90:0.5) ;
\end{tikzpicture}
\caption[Conservation of magnetic flux]
{Conservation of magnetic flux: There exists a closed untwisted
  2-form $F$ called the electromagnetic field (blue). This encodes the two
  fields $E$ and $B$. Since we are in 4-dimensions the orientation is
  a loop outside the 2-dim form-submanifolds.
 The closure of $F$ is equivalent to the two macroscopic
  Maxwell equations which contain $E$ and $B$.}
\label{fig_Max_F}
\end{figure}

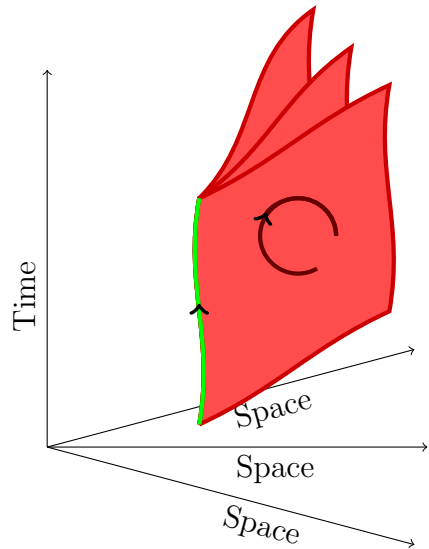
\begin{figure}[tb]
\begin{tikzpicture}
\draw[->] (0,0) -- (0,5) ;
\draw (-0.3,2) node [rotate=90] {Time} ;
\draw[->] (0,0) -- (5,0) ;
\draw (3,-0.3) node [rotate=0] {Space} ;
\draw[->,rotate=-15] (0,0) -- (5,0) ;
\draw[rotate=-15]  (3,-0.3) node [rotate=-15] {Space} ;
\draw[->,rotate=15] (0,0) -- (5,0) ;
\draw[rotate=15]  (3,-0.3) node [rotate=15] {Space} ;

\filldraw[ultra thick,draw=red!80!black,fill=red!70!white,shift={(2,.3)}] 
(0,0) to[out=80,in=-100] (0,3) to [out=45,in=-145] (1.5,5.5) 
to [out=-100,in=80] (1.5,2.5) to [out=-140,in=45] (0,0); 

\filldraw[ultra thick,draw=red!80!black,fill=red!70!white,shift={(2,.3)}] 
(0,0) to[out=80,in=-100] (0,3) to [out=35,in=-145] (2,5) 
to [out=-100,in=80] (2,2) to [out=-140,in=35] (0,0); 

\filldraw[ultra thick,draw=red!80!black,fill=red!70!white,shift={(2,.3)}] 
(0,0) to[out=80,in=-100] (0,3) to [out=25,in=-155] (2.5,4.5) 
to [out=-100,in=80] (2.5,1.5) to [out=-155,in=25] (0,0); 

\draw[ultra thick,draw=green,shift={(2,.3)}] 
(0,0) to[out=80,in=-100] (0,3); 
\draw[ultra thick,draw=green,shift={(2,1.8)},->] 
(0,0) -- +(0,0.1) ; 

\draw[ultra thick,draw=red!40!black,shift={(2,.3)}] (1.8,2.5) arc
(0:300:.5) ;
\draw[ultra thick,draw=red!40!black,shift={(2,.3)},-<] (1.8,2.5) arc
(0:160:.5) ;

\end{tikzpicture}
\caption[The excitation fields]{The excitation fields: There exists a
  twisted 2-form $\Hform$ called the excitation field (red). It
  encodes the fields $D$ and $H$. In general it is not closed and its
  exterior derivative is $\Jcurr$ (green). This is equivalent to the
  two macroscopic Maxwell equations which contain $D$ and $H$.  }
\label{fig_Max_H}
\end{figure}

\begin{figure}[b!]
\begin{tikzpicture}

\draw[shift={(-3,0)},->] (0,0) -- (0,5) ;
\draw[shift={(-3,0)}] (-0.3,2) node [rotate=90] {Time} ;
\draw[shift={(-3,0)},->] (0,0) -- (5,0) ;
\draw[shift={(-3,0)}] (3,-0.3) node [rotate=0] {Space} ;
\draw[shift={(-3,0)},->,rotate=-15] (0,0) -- (5,0) ;
\draw[shift={(-3,0)},rotate=-15]  (3,-0.3) node [rotate=-15] {Space} ;
\draw[shift={(-3,0)},->,rotate=15] (0,0) -- (5,0) ;
\draw[shift={(-3,0)},rotate=15]  (3,-0.3) node [rotate=15] {Space} ;

\filldraw[ultra thick,draw=red!80!black,fill=red!70!white,rotate=180] 
(-1,-2.5) -- +(2,-1) -- +(4,0) -- +(2,1) -- cycle ;
\filldraw[ultra thick,draw=blue!80!black,fill=blue!70!white] 
(0,0) to[out=80,in=-100] (0,3) to [out=35,in=-145] (2,5) 
to [out=-100,in=80] (2,2) to [out=-140,in=35] (0,0); 
\filldraw[ultra thick,draw=red!80!black,fill=red!70!white] 
(1,2.5) -- +(2,-1) -- +(4,0) -- +(2,1) -- cycle ;

\draw[ultra thick,white,double=white,shift={(0.9,3.7)},yscale=0.5] (0,0) arc
(-190:10:0.5) ;
\draw[ultra thick,red!50!black,shift={(0.9,3.7)},yscale=0.5] (0,0) arc
(-190:10:0.5) ;
\draw[ultra thick,red!50!black,shift={(0.9,3.7)},-<,yscale=0.5] (0,0) arc
(-190:-90:0.5) ;

\draw[ultra thick,draw=red!40!black,shift={(3.5,2.5)},yscale=0.5] (0,0) arc
(0:300:.5) ;
\draw[ultra thick,draw=red!40!black,shift={(3.5,2.5)},-<,yscale=0.5] (0,0) arc
(0:160:.5) ;

\draw[thick,shift={(1,2.5)}] (0.4,0) -- +(0,0.4) -- +(-0.4,0.6) ;

\end{tikzpicture}
\caption[The constitutive relations in the vacuum]
{The constitutive relations in the vacuum: The twisted 2-form
  ${\cal H}$ (red) is the Hodge dual of the untwisted 2-form ${F}$ (blue). Here we
  attempt to show that ${\cal H}$ and $F$ are orthogonal 2-dim
  form-submanifolds in
  4-dimensions. They have
  the same orientation.} 
\label{fig_Max_FH}
\end{figure}
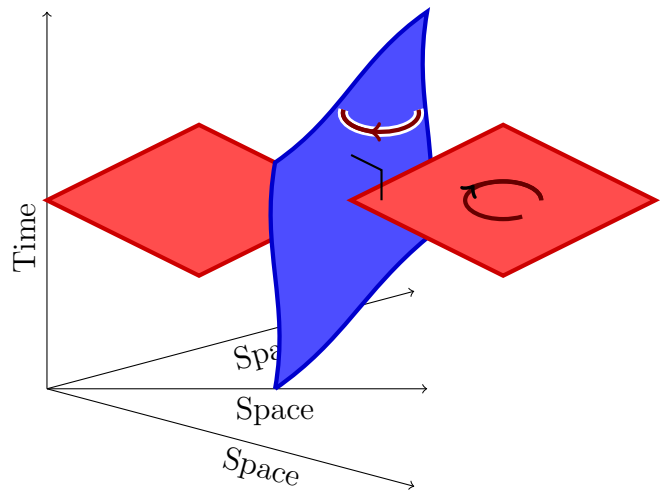

Figures \ref{fig_Max_J} to \ref{fig_Max_FH} represent conservation of
charge and Maxwell's equations in spacetime. Thus there are four axes,
one of time and three of space.



We start with conservation of charge, figure \ref{fig_Max_J}, which is
not one of Maxwell's equations but a simple consequence of them. In figure
\ref{fig_Max_J} we see the 1-dim form-submanifolds of $\Jcurr$ do not
have a boundary and therefore $\Jcurr$ is closed. One may think of the
form-submanifolds of $\Jcurr$ as worldlines of the individual
point charges such as electrons or ions.

Look again at figure \ref{fig_EM_Mag_2-form}, but now consider that the
fields are not static. We may imagine that the magnetic field lines
move in space.  As a result some of the field lines may cross the red
circle. This will generate an electromotive force as given by the
integral form of Faraday's law of induction.  In this case the
form-submanifolds of $B$ will map out a 2-dim form-submanifold in
4-dim spacetime. These form-submanifold will not have boundaries and
therefore it corresponds to the untwisted closed 2-form $F$. In figure
\ref{fig_Max_F} we depict the closed 2-form $F$ as
2-dim without boundary in 4-dim spacetime. However one should recall
that 2-forms in 4-dimensions can self intersect as indicated in
section \ref{ch_four-dim}. 

Repeating the same procedure for figure
\ref{fig_EM_Gauss_ElectroStatics}. If the charges $\varrho$ are not
static then they map out worldlines in spacetime. This is the closed 3-form
$\Jcurr$ as in figure \ref{fig_Max_J}. As a consequence the 1-dim
form-submanifolds for the 2-form $D$ map out the 2-dim
form-submanifolds for the 2-form $\Hform$. This gives figure
\ref{fig_Max_H}. 

As stated before, Maxwell's equations need to be augmented using the
constitutive relations. In the vacuum one has that $D$ and $H$ are the
Hodge duals of $E$ and $B$, as can be seen in figures
\ref{fig_EM_Gauss_ElectroStatics} and
\ref{fig_EM_Ampere_MagnetoStatics}. In figure \ref{fig_Max_FH} we see
that the twisted 2-form $\Hform$ is the Hodge dual of untwisted 2-form
$F$. In other media the relationship between $F$ and $\Hform$ is more
complicated.


\subsection{Lorentz force equation}
\label{ch_Max_Lorentz}

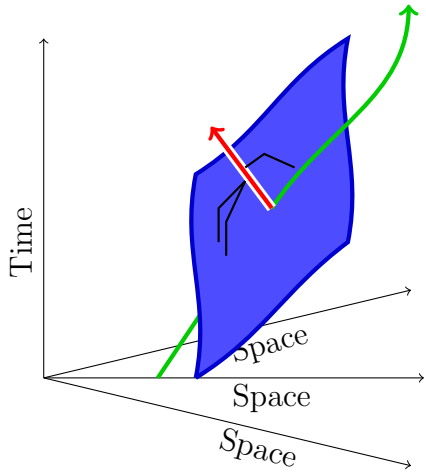
\begin{figure}[tb]
\centering
\begin{tikzpicture}[yscale=.9]
\draw[shift={(-2,0)},->] (0,0) -- (0,5) ;
\draw[shift={(-2,0)}] (-0.3,2) node [rotate=90] {Time} ;
\draw[shift={(-2,0)},->] (0,0) -- (5,0) ;
\draw[shift={(-2,0)}] (3,-0.3) node [rotate=0] {Space} ;
\draw[shift={(-2,0)},->,rotate=-15] (0,0) -- (5,0) ;
\draw[shift={(-2,0)},rotate=-15]  (3,-0.3) node [rotate=-15] {Space} ;
\draw[shift={(-2,0)},->,rotate=15] (0,0) -- (5,0) ;
\draw[shift={(-2,0)},rotate=15]  (3,-0.3) node [rotate=15] {Space} ;

\draw[ultra thick,green!80!black](-.5,0) -- (1,2.5) ;

\filldraw[ultra thick,draw=blue!80!black,fill=blue!70!white] 
(0,0) to[out=80,in=-100] (0,3) to [out=35,in=-145] (2,5) 
to [out=-100,in=80] (2,2) to [out=-140,in=35] (0,0); 

\draw[ultra thick,green!80!black,->]
(1,2.5) to[out=60,in=-90] +(1.5*1.2,2.5*1.2);

\draw[ultra thick,double=white,white] (1,2.5) -- +(-.8,1.2);
\draw[ultra thick,red!100,->] (1,2.5) -- +(-.8,1.2);


\draw[thick] (1.3,3.1) -- (0.9,3.3) -- (0.65,3.1) ;
\draw[thick] (0.65,2.9) -- (0.3,2.5) -- (0.3,2.) ;
\draw[thick] (0.65,2.9) -- (0.4,2.3) -- (0.4,1.8) ;
\end{tikzpicture}

\caption[The Lorentz force law.]{An attempt to draw the Lorentz force
  law on a particle. 
  Here the 2-form $F$ is the 2-dim form-submanifold (blue), the
  worldline of the particle is the curve (green) and the force is the
  straight arrow (red). The force is orthogonal to to both the
  2-form $F$ and the 4-velocity of the particle.}

\label{fig_Lorentz_force}
\end{figure}

In electrodynamics there is one remaining equation, the Lorentz force
law. As depicted in figure \ref{fig_Lorentz_force}. This states that
in spacetime the force is orthogonal to both the 4-velocity of the
particle and the electromagnetic 2-form $F$.



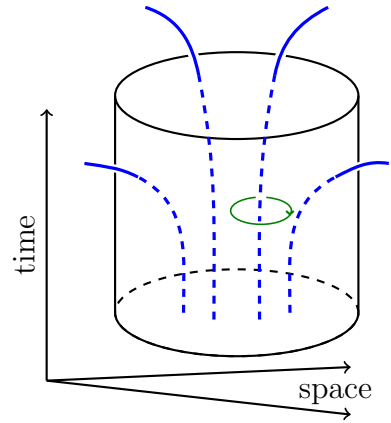
\begin{figure}[tb]
\centering
\begin{tikzpicture}[yscale=.9]
\draw[thick,yscale=0.4,dashed] (0,0) circle (1.6);
\draw[thick,yscale=0.4] (1.6,0) arc (0:-180:1.6);
\draw[thick,yscale=0.4] (0,8) circle (1.6);
\draw[thick] (-1.6,0) -- +(0,3.2) ;
\draw[thick] (1.6,0) -- +(0,3.2) ;

\draw[very thick, blue,dashed] (-.7,0) to[out=90,in=-30] (-1.3,2)
to[out=150,in=-10] +(-.7,.2);
\draw[very thick, double=blue,white] (-1.3,2) to[out=150,in=-10] +(-.7,.2);
\draw[very thick, blue] (-1.3,2) to[out=150,in=-10] +(-.7,.2);

\draw[very thick, blue,dashed] (.7,0) to[out=90,in=-150] (1.3,2)
to[out=30,in=170] +(.7,.2);
\draw[very thick, double=white,white] (1.3,2) to[out=30,in=170] +(.7,.2);
\draw[very thick, blue] (1.3,2) to[out=30,in=170] +(.7,.2);

\draw[very thick, blue,dashed] (-.3,-0.1) to[out=90,in=-80] (-.5,3.5)
to[out=100,in=-20] +(-.7,1);
\draw[ultra thick, double=blue,white] (-.5,3.5) to[out=100,in=-20] +(-.7,1);
\draw[very thick, blue] (-.5,3.5) to[out=100,in=-20] +(-.7,1);

\draw[very thick, blue,dashed] (.3,-0.1) to[out=90,in=-100] 
(.5,3.5) to[out=80,in=-150] +(.7,1);
\draw[ultra thick, double=blue,white] (.5,3.5) to[out=80,in=-150] +(.7,1);
\draw[very thick, blue] (.5,3.5) to[out=80,in=-150] +(.7,1);;


\draw [thick,->] (-2.5,-1) -- + (0,4) ;
\draw (-2.8,1) node {\rotatebox{90}{time}} ; 
\draw [thick,->] (-2.5,-1) -- + (4,.2) ;
\draw [thick,->] (-2.5,-1) -- + (4,-.5) ;
\draw (1.3,-1.2) node {{space}} ; 

\draw[thick,yscale=-0.5,draw=white,double=green!50!black] 
(0.25,-3.4) arc (260:-80:0.4) ; 
\draw[thick,draw=green!50!black] 
(0.65,1.5) -- +(0.05,-0.05) -- +(0.1,0) ; 

\end{tikzpicture}

\caption[Faraday's law of induction on the plane $z=0$.]{Faraday's law
  of induction on the plane $z=0$, integrated over time and
  space. The blue lines are pullback of the electromagnetic 2-form $F$
  onto the plane $z=0$. The blue lines are continuous, i.e. as $F$ is
  closed and therefore the pullback of $F$ onto the plane $z=0$ is
  closed. Thus the difference between the number of lines entering the
  bottom disc and the number leaving the top disc equals the number of
  lines leaving the sides of the cylinder.}
\label{fig_EM_Faraday}
\end{figure}


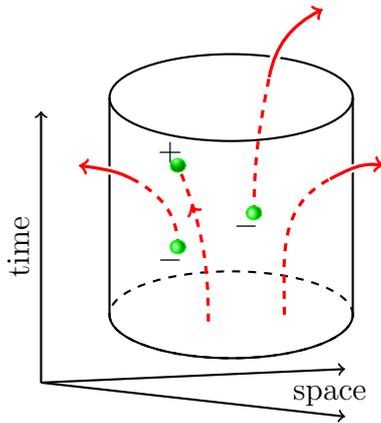
\begin{figure}[tb]
\centering
\begin{tikzpicture}[yscale=.9]
\draw[thick,yscale=0.4,dashed] (0,0) circle (1.6);
\draw[thick,yscale=0.4] (1.6,0) arc (0:-180:1.6);
\draw[thick,yscale=0.4] (0,8) circle (1.6);
\draw[thick] (-1.6,0) -- +(0,3.2) ;
\draw[thick] (1.6,0) -- +(0,3.2) ;

\draw[very thick, red,dashed] (-.7,1) to[out=90,in=-30] (-1.3,2)
to[out=150,in=-10] +(-.7,.2);
\draw[very thick, double=red,white] (-1.3,2) to[out=150,in=-10] +(-.7,.2);
\draw[very thick, red,->] (-1.3,2) to[out=150,in=-10] +(-.7,.2);
\shade [ball color=green,draw=green] (-.7,1) circle (0.1) +(-.1,-.2) node{$-$};

\draw[very thick, red,dashed] (.7,0) to[out=90,in=-150] (1.3,2)
to[out=30,in=170] +(.7,.2);
\draw[very thick, double=white,white] (1.3,2) to[out=30,in=170] +(.7,.2);
\draw[very thick, red,->] (1.3,2) to[out=30,in=170] +(.7,.2);

\draw[very thick, red,dashed] (-.3,-0.1) to[out=90,in=-70] (-.7,2.2);
\draw[very thick, red,->] (-.48,1.5) -- +(-.05,0.1);
\shade [ball color=green,draw=green] (-.7,2.2) circle (0.1) +(-.1,.2) node{+};

\draw[very thick, red,dashed] (.3,1.5) to[out=90,in=-100] (.5,3.5) to[out=80,in=-150] +(.7,1);
\draw[ultra thick, double=red,white] (.5,3.5) to[out=80,in=-150] +(.7,1);
\draw[very thick, red,->] (.5,3.5) to[out=80,in=-150] +(.7,1);;
\shade [ball color=green,draw=green] (.3,1.5) circle (0.1) +(-.1,-.2) node{$-$};

\draw [thick,->] (-2.5,-1) -- + (0,4) ;
\draw (-2.8,1) node {\rotatebox{90}{time}} ; 
\draw [thick,->] (-2.5,-1) -- + (4,.2) ;
\draw [thick,->] (-2.5,-1) -- + (4,-.5) ;
\draw (1.3,-1.2) node {{space}} ; 

\end{tikzpicture}

\caption[Maxwell Ampères law on the plane $z=0$.]{Maxwell Ampères law
  on the plane $z=0$, integrated over time and space. The green dots
  are the intersections of the positive and negative charges as they
  cross $z=0$. The red lines are the pullback of the excitation 2-form
  $\Hform$ onto the plane $z=0$. These terminate at charges. Thus the
  difference between the number of lines entering the bottom disc and
  the number leaving the top and sides of the cylinder equals the
  charge inside the cylinder.}
\label{fig_EM_AmpMax}
\end{figure}

\subsection{Constant space slice of Maxwell's equations.}
\label{ch_Max_3dim}

We may consider the static figures \ref{fig_EM_Mag_2-form} and
\ref{fig_EM_Gauss_ElectroStatics} correspond to the pullback, as
discussed in section \ref{ch_Maps}, of the
fields $F$, $\Hform$ and $\Jcurr$ on a 3-dimensional slice through
spacetime. 

An alternative slice of Maxwell's equations in spacetime is to
choose a plane in space and consider what the electromagnetic fields
look like when pulled back onto this plane for all time. We will call
this plane the $z=0$ plane. The result is a diagram with one time
coordinate and two space coordinates. 

The 2-form $F$ is a closed
untwisted 1-dimensional form-submanifold. However this time it
contains information about both $B$ and $E$. Figure
\ref{fig_EM_Faraday} shows Faraday's law of induction for a loop in
the plane $z=0$. Since $F$ is
closed the difference between the field at the top and bottom, which
is the difference in the $B$ field, must
correspond to field leaving the loop, which is the integral of the
electromotive force.
Similarly figure \ref{fig_EM_AmpMax} is the Maxwell-Ampères law. 

\section{Conclusion}
\label{ch_Conclusion}

We have presented a large proportion of differential geometry,
concentrating on differential forms, in order to arrive at figures
\ref{fig_Max_J}-\ref{fig_Lorentz_force} which are the diagrams that
describe electrodynamics.

In terms of pictorially presenting differential geometry, there are
still some open challenges. One is in the attempt to depict 2-forms in
4 dimensions, section \ref{ch_four-dim}. One question is whether for a
given 2-form, figure \ref{fig_2form_4dim_plane} or
\ref{fig_2form_4dim_moving_line} is correct and whether there are
alternatives which give more or different information. 

There are many phenomena in electromagnetism which may gain from a
pictorial representation. For example, how is it best to draw the the
electromagnetism potential, and can one shed light on the Aharonov–Bohm effect? 
How does one depict the electromagnetic fields as well as the flow of
electrons and ions in a plasma. In this case, ideally one would have 7
dimensional paper.


\section*{Acknowledgments}
I am grateful grateful for the support provided by STFC (the
Cockcroft Institute ST/G008248/1 and ST/P002056/1), EPSRC (the Alpha-X project
EP/J018171/1  and EP/N028694/1). I am grateful  
to the masters students who took my course Advanced
Electromagnetism and Gravity and 2017 and read the document, making
suggestions and to Taylor Boyd and Richard Dadhley who pointed
out some inconsistencies.

All figures are used with kind permission of Jonathan Gratus
\cite{Gratus-Internal-report}.

\bibliographystyle{unsrt}
\bibliography{Geometry}

\end{document}